\documentclass[12pt,leqno]{article}
\usepackage{amssymb}
\usepackage[mathscr]{eucal}
\usepackage{amsmath,amssymb,latexsym,theorem,bbm}
\usepackage{amsmath,amssymb,latexsym,theorem}
\usepackage{color,url}
\usepackage{appendix}

\newcommand{\NN}{\mathbb{N}}

\newcommand{\RR}{\mathbb{R}}

\newcommand{\ZZ}{\mathbb{Z}}

\newcommand{\bA}{{\boldsymbol{A}}}

\newcommand{\bB}{{\boldsymbol{B}}}

\newcommand{\tb}{\widetilde{b}}

\newcommand{\tbB}{{\widetilde{\bB}}}
\newcommand{\tB}{\widetilde{B}}

\newcommand{\bc}{{\boldsymbol{c}}}

\newcommand{\be}{{\boldsymbol{e}}}

\newcommand{\bI}{{\boldsymbol{I}}}

\newcommand{\bell}{{\boldsymbol{\ell}}}

\newcommand{\bQ}{{\boldsymbol{Q}}}

\newcommand{\br}{{\boldsymbol{r}}}

\newcommand{\bv}{{\boldsymbol{v}}}

\newcommand{\tbv}{\widetilde{\bv}}
\newcommand{\bx}{{\boldsymbol{x}}}
\newcommand{\bX}{{\boldsymbol{X}}}
\newcommand{\by}{{\boldsymbol{y}}}

\newcommand{\bz}{{\boldsymbol{z}}}

\newcommand{\bU}{{\boldsymbol{U}}}

\newcommand{\Bbeta}{{\boldsymbol{\beta}}}

\newcommand{\tBbeta}{\widetilde{\Bbeta}}

\newcommand{\blambda}{{\boldsymbol{\lambda}}}

\newcommand{\bxi}{{\boldsymbol{\xi}}}

\newcommand{\bmu}{{\boldsymbol{\mu}}}

\newcommand{\bbeta}{{\boldsymbol{\beta}}}

\newcommand{\bvarphi}{{\boldsymbol{\varphi}}}
\newcommand{\bzero}{{\boldsymbol{0}}}
\newcommand{\bone}{{\boldsymbol{1}}}

\newcommand{\cB}{{\mathcal B}}

\newcommand{\cF}{{\mathcal F}}

\newcommand{\cU}{{\mathcal U}}

\newcommand{\tw}{\widetilde{w}}

\newcommand{\biota}{\boldsymbol{\iota}}

\newcommand{\dd}{\mathrm{d}}
\newcommand{\ee}{\mathrm{e}}
\newcommand{\ii}{\mathrm{i}}

\newcommand{\EE}{\operatorname{\mathbb{E}}}

\newcommand{\PP}{\operatorname{\mathbb{P}}}

\newcommand{\card}{\operatorname{\mathrm{card}}}

\newcommand{\tN}{\widetilde{N}}

\newcommand{\tv}{\widetilde{v}}

\newcommand{\vare}{\varepsilon}

\renewcommand{\mid}{\,|\,}

\renewcommand{\leq}{\leqslant}
\renewcommand{\geq}{\geqslant}

\newcommand{\bbone}{\mathbbm{1}}

\newcommand{\proofend}{\hfill\mbox{$\Box$}}

\numberwithin{equation}{section}

\theoremstyle{change} \theorembodyfont{\em}
\newtheorem{Lem}{Lemma.}[section]
\newtheorem{Thm}[Lem]{Theorem.}
\newtheorem{Pro}[Lem]{Proposition.}
\newtheorem{Cor}[Lem]{Corollary.}
\newtheorem{Def}[Lem]{Definition.}

\theorembodyfont{\rm}
\newtheorem{Rem}[Lem]{Remark.}

\begin{document}

\begin{center}
 {\bfseries\Large Distributional properties of first jump times of\\[1mm]
   CBI processes with jump sizes in given Borel sets}

\vspace*{3mm}

 {\sc\large
 M\'aty\'as $\text{Barczy}^{*,\diamond}$,
  \ Sandra $\text{Palau}^{**}$,
  \ Yao $\text{Xue}^{***}$
  }

\end{center}

\vskip0.2cm

\noindent
 * HUN-REN--SZTE Analysis and Applications Research Group,
   Bolyai Institute, University of Szeged,
   Aradi v\'ertan\'uk tere 1, H--6720 Szeged, Hungary.

\noindent
 ** Instituto de Investigaciones en Matem\'aticas Aplicadas y en Sistemas,
    Universidad Nacional Aut\'onoma de M\'exico, CDMX, 04510, Ciudad de M\'exico, M\'exico.

\noindent
 *** School of Mathematical Sciences, Beijing Normal University, Beijing 100875, People's Republic of China.

\noindent E-mails: barczy@math.u-szeged.hu (M. Barczy),
                  sandra@sigma.iimas.unam.mx (S. Palau),
                  xuey1997@163.com (Y. Xue).

\noindent $\diamond$ Corresponding author.

\vskip0.2cm

%\centerline{\sl May 28, 2018.}

{\renewcommand{\thefootnote}{}
\footnote{Sandra Palau is supported by UNAM-PAPIIT (IN103924). Yao Xue is supported by the National Key R{\&}D Program of China (No.~2020YFA0712900).}
}

\vspace*{-10mm}

\begin{abstract}
We derive an expression for the joint distribution function of the first jump times
 of a continuous state and continuous time branching process with immigration (CBI process)
 with jump sizes in given Borel sets having finite total L\'evy measures,
 where by total L\'evy measure, we mean the sum of the measures appearing in the branching and
 immigration mechanisms of the CBI process in question.
Our result generalizes a corresponding result of He and Li (2016),
 who considered this problem in case of a single Borel set having finite total L\'evy measure.
 \end{abstract}

%\tableofcontents

\section{Introduction}
\label{section_intro}

Branching processes form an important class of stochastic processes,
 since they can be used for modeling random evolution of a population consisting of individuals,
 where individuals can generate offspring and immigrants may join the population as well. 
They can find applications in disease spread, biological reproduction, nuclear chain reactions, and computer virus propagation.
An important class of branching processes is the class of continuous state and continuous time branching processes without or with immigration (CB and CBI processes, respectively)
 that are special Markov processes with values in $[0,\infty)$. 
The study of CB processes was initiated by Feller \cite{Fel}, who proved that a diffusion
 process may arise in a limit theorem for (discrete time) Galton-Watson branching processes.
CBI processes cover the situation where immigrants may come to the population from outer sources. 
For a brief introduction to CBI processes, see Li \cite{Li3}. 
CBI processes can be also viewed as special measure-valued branching Markov processes, for a modern book in this field, 
 see Li \cite{Li}. 
 
In this paper, we study the distributional properties of some jump times of CBI processes.
Given a CBI process $(X_u)_{u\geq 0}$, see Definition \ref{Def_CBI} with $d=1$, and a Borel set $A$ of $(0,\infty)$,
 let us introduce
 \[
  \tau_A:= \inf\{u > 0: \Delta X_u \in A\} 
 \]
 with the convention $\inf(\emptyset) := \infty$, and
 \[
   J_t(A):= \begin{cases} 
                   \card(\{u \in (0, t] : \Delta X_u  \in A\}) & \text{if $t>0$,}\\
                    0 &\text{if $t=0$,}
                  \end{cases}  
 \]
  where $\Delta X_u := X_u - X_{u-}$, $u >0$,
 and $\card(H)$ denotes the cardinality of a set $H$.
One can call $\Delta X_u\in[0,\infty)$ the size of the jump of $X$ at time $u>0$.

He and Li \cite[Theorem 3.1]{HeLi} studied the distribution of $J_t(A)$ and $\tau_A$, namely,
 they showed that $\PP(J_t(A)<\infty)=1$, $t>0$, and derived an expression for
 $\PP(\tau_A>t \mid X_0=x)$, $t\geq 0$, $x\geq 0$, in terms of a solution of
 a deterministic differential equation provided that the set $A$ has a finite total L\'evy measure, which is defined
 as the sum of the measures appearing in the branching and immigration mechanisms of $(X_u)_{u\geq 0}$.
The above recalled result of He and Li \cite{HeLi} on the distributional properties of jumps of CBI processes
 is also contained in Section 10.3 in the new second edition of Li's book \cite{Li}.

Barczy and Palau \cite[Theorem 4.2]{BarPal} investigated the distributional properties of jumps of multi-type CBI processes.
They derived an expression for the distribution function of the first jump time of a multi-type
 CBI process with jump size in a given Borel set having finite total L\'evy measure,
 see also Theorem \ref{Thm_jump_times_dCBI} of the present paper.

In this paper, we give a generalization of the above recalled result of He and Li \cite[Theorem 3.1]{HeLi}
 on the distribution of $\tau_A$.
Let $(X_t)_{t\geq 0}$ be a CBI process, $k$ be a positive integer, and $A_1,\ldots,A_k$ be Borel sets of $(0,\infty)$
 such that each $A_i$, $i\in\{1,\ldots,k\}$, has a finite total L\'evy measure.
We derive an expression for the probability
 $\PP(\tau_{A_1}>t_1,\ldots, \tau_{A_k}>t_k \mid X_0=x)$, where $t_1,\ldots,t_k>0$ 
 and $x\geq 0$, in terms of a solution of a system of deterministic differential equations.
Next, we provide a possible motivation why it is worth studying such an extension.
CBI processes are sometimes used for modelling some non-negative financial quantities
  such as short interest rates or instantaneous variances.
In particular, the so-called $\alpha$-stable Cox-Ingersoll-Ross (CIR) processes, which are special CBI processes,
 are used for such a purpose (for the definition of $\alpha$-stable CIR processes, 
 see, e.g., Remark \ref{Rem_financial_application}).
When an $\alpha$-stable CIR process describes the evolution of short interest rates,
 an upward jump may be interpreted as an abrupt increase in funding costs. 
It is reasonable to study not only the occurrence of an individual large shock, but also the joint arrival
 of shocks of different sizes over different times (for more details, see Remark \ref{Rem_financial_application}).

Now we turn to summarize the content of our paper, it is structured as follows.
At the end of the Introduction, we collect all the notations that will be used throughout the paper.
In Section \ref{section_CBI}, we recall the definition of multi-type CBI processes, 
 a result on their representation as pathwise unique strong solutions
 of stochastic differential equations (see \eqref{SDE_atirasa_dimd}),
 and an expression for the distribution function of the first jump time of a multi-type
 CBI process with jump size in a given Borel set having finite total L\'evy measure 
 (see Theorem \ref{Thm_jump_times_dCBI}).

Section \ref{Section_CBI_results_1} is devoted to
 deriving a formula for the joint Laplace transform of a CBI process 
and its first jump times with jump sizes in given Borel sets having finite total L\'evy measures.
As a first step, given a CBI process $(X_t)_{t\geq 0}$,
 a positive integer $k$ and Borel sets $A_1,\ldots,A_k$ of $(0,\infty)$ having finite total L\'evy measures,
 in Theorem \ref{Thm_(k+1)CBI}, we prove that 
 $(X_t,J_t(A_1),\ldots,J_t(A_k))_{t\geq 0}$ is a $(k+1)$-type CBI process,
 and we determine its parameters.
Theorem \ref{Thm_(k+1)CBI} is a generalization of Theorem 10.11 in Li \cite{Li}, namely,
 in the special case $k=1$, Theorem \ref{Thm_(k+1)CBI} recovers Theorem 10.11 in Li \cite{Li}.
We also point out the fact that our proof of Theorem \ref{Thm_(k+1)CBI} is slightly different from that of 
 Theorem 10.11 in Li \cite{Li},
 we directly use a stochastic differential equation (SDE) representation of multi-type CBI processes (see \eqref{SDE_atirasa_dimd}).
Next, in Corollary \ref{Cor_Laplace}, having the parameters of the $(k+1)$-type
 CBI process in question at hand, we explicitly write out its branching and immigration mechanisms,
 and derive a formula for the Laplace transform of $(X_t,J_t(A_1),\ldots,J_t(A_k))$ 
 in terms of a solution of a deterministic differential equation, where $t>0$.
Corollary \ref{Cor_Laplace} with $k=1$ gives back Corollary 10.12 in Li \cite{Li}.
We close Section \ref{Section_CBI_results_1} with the auxiliary Lemma \ref{Lem_vjumps},
 which contains some monotonicity- and convergence-type properties of the unique locally bounded solution to the differential equation \eqref{def_v_k+1_CBI},
 that is related to the system of differential equations associated with the branching mechanism of
 $(X_t,J_t(A_1),\ldots,J_t(A_k))_{t\geq 0}$.
Lemma \ref{Lem_vjumps} is used in the proof of our main result Theorem \ref{Thm_main} (discussed below).
The proof of Lemma \ref{Lem_vjumps} is postponed to Appendix \ref{Appendix_Lemma_v} for better readability of the paper.

In Section \ref{Section_CBI_results_2}, we investigate the joint distribution function of a CBI process
 and its first jump times with jump sizes in given Borel sets with finite total L\'evy measures.
First, in Theorem \ref{Thm_single_jump_times_CBI}, we recall a known result due to He and Li \cite[Theorem 3.1]{HeLi}
 on the distribution function of $\tau_A$, where $A$ is a Borel set of $(0,\infty)$ having finite total L\'evy measure.
Theorem \ref{Thm_main} contains the main result of the paper: 
 given a positive integer $k$ and Borel sets $A_1,\ldots,A_k$ of $(0,\infty)$ having finite total L\'evy measures, 
 we derive an expression for the joint distribution function of $(\tau_{A_1},\ldots,\tau_{A_k})$ 
 at a point $(t_1,\ldots,t_k)\in\RR^k$ satisfying $0\leq t_1\leq t_2\leq\ldots\leq t_k$,
 in terms of the solution of a system of deterministic differential equations.
We close Section \ref{Section_CBI_results_2} with Corollary \ref{Cor_multiple_jump_times_CBI}, where we present 
 an expression for the joint distribution function of $(\tau_{A_1},\ldots,\tau_{A_k})$ 
 at a point $(t,\ldots,t)\in\RR^k$, where $t\geq 0$. 
We give two proofs of Corollary \ref{Cor_multiple_jump_times_CBI}.
The first one is a direct application of the known result of He and Li \cite[Theorem 3.1]{HeLi} (see also Theorem \ref{Thm_single_jump_times_CBI}),
 and the second one is a specialization of our main Theorem \ref{Thm_main} to the case $t_1=\cdots=t_k=t$.

We close the paper with two appendices.
In Appendix \ref{Appendix_multitype_CBI}, we derive an expression 
 for the joint Laplace transform of a multi-type CBI process,
 which may be interesting in itself as well, see Proposition \ref{Pro_CBI_joint_Laplace}.
This result is specialized and rewritten in case of a single-type CBI process,
 see Proposition \ref{Pro_Li_CB_joint_Laplace} and Remark \ref{Rem_AppA_2}.
In Proposition \ref{Relation_g_v}, we provide a connection 
 between the solutions of two differential equations 
 appearing in  Propositions \ref{Pro_CBI_joint_Laplace} and \ref{Pro_Li_CB_joint_Laplace}  
 for the joint Laplace transform of a single-type CB process,
 by giving a proof which does not use the theory of CB processes at all.
Appendix \ref{Appendix_Lemma_v} is devoted to the proof of the auxiliary Lemma \ref{Lem_vjumps}, 
 which plays a crucial role in the proof of our main result Theorem \ref{Thm_main}.

Finally, we introduce the notations that will be used throughout the paper.
Let $\ZZ_+$, $\NN$, $\RR$, $\RR_+$ and $\RR_{++}$ denote the
 set of non-negative integers, positive integers, real numbers, non-negative real
 numbers and positive real numbers, respectively.
For each $d\in\NN$, let $\cU_d:= \RR_+^d \setminus \{\bzero\}$.
The identity mapping on $\cU_1$ is denoted by $\mathrm{id}_{\cU_1}$.
For a subset $U$ of $\cU_d$, we will use the notation $U^{-1}:=\cU_d\setminus U$, and $\bbone_U$ denotes the indicator function of $U$.
For $k\in\NN$, \ $\bell=(\ell_1,\ldots,\ell_k)\in\{-1,1\}^k$ and $U_1,\dots, U_k\subset \cU_d$,
 we define $\bU^{\bell}:=U_1^{\ell_1}\cap\dots\cap U_k^{\ell_k}$.
For $x , y \in \RR$, we will use the notations
 $x \land y := \min \{x, y\}$ and $x^+ := \max \{0, x\}$.
For a function $g:\RR_+\times\RR_+^d\to \RR^d$, $\partial_1g$ denotes the partial derivative of $g$ with
 respect to its first variable (provided that it exists).
By $\|\bx\|$ and $\|\bQ\|$, we denote the norm and the induced norm of $\bx \in \RR^d$ and
 $\bQ \in \RR^{d\times d}$, respectively.
The null vector and the null matrix will be denoted by $\bzero$.
Moreover, $\bI_d \in \RR^{d\times d}$ denotes the $d\times d$ identity matrix,
 and $\be_1^{(d)}$, \ldots, $\be_d^{(d)}$ denotes the natural basis in $\RR^d$.
The Borel $\sigma$-algebra on a subset $U$ of $\RR^d$ is denoted by $\cB(U)$,
 and recall that $\cB(U) = U\cap \cB(\RR^d)$.
By a Borel measure on a Borel set $S\in\cB(\RR^d)$, we mean a measure on $(S, \cB(S))$.
Every random variable will be defined on an appropriate probability space $(\Omega,\cF,\PP)$.
Throughout this paper, we make the conventions $\int_a^b := \int_{(a,b]}$
 and $\int_a^\infty := \int_{(a,\infty)}$ for any $a, b \in \RR$ with $a \leq b$.

\section{Preliminaries on multi-type CBI processes}
\label{section_CBI}

In this section, we recall some preliminaries on multi-type CBI processes that we use in the paper.

\begin{Def}\label{Def_essentially_non-negative}
A matrix \ $\bQ = (q_{i,j})_{i,j\in\{1,\ldots,d\}} \in \RR^{d\times d}$ \ is
 called essentially non-negative if \ $q_{i,j} \in \RR_+$ \ whenever
 \ $i, j \in \{1,\ldots,d\}$ \ with \ $i \ne j$.
The set of essentially non-negative \ $d \times d$ \ matrices will be denoted by
 \ $\RR^{d\times d}_{(+)}$.
\end{Def}

\begin{Def}\label{Def_admissible}
A tuple \ $(d, \bc, \Bbeta, \bB, \nu, \bmu)$ \ is called a set of admissible
 parameters if
 \renewcommand{\labelenumi}{{\rm(\roman{enumi})}}
 \begin{enumerate}
  \item
   $d \in \NN$,
  \item
   $\bc = (c_i)_{i\in\{1,\ldots,d\}} \in \RR_+^d$,
  \item
   $\Bbeta = (\beta_i)_{i\in\{1,\ldots,d\}} \in \RR_+^d$,
  \item
   $\bB = (b_{i,j})_{i,j\in\{1,\ldots,d\}} \in \RR^{d \times d}_{(+)}$,
  \item
   $\nu$ \ is a Borel measure on \ $\cU_d$
    \ satisfying \ $\int_{\cU_d} (1 \land \|\br\|) \, \nu(\dd\br) < \infty$,
  \item
   $\bmu = (\mu_1, \ldots, \mu_d)$, \ where, for each
    \ $i \in \{1, \ldots, d\}$, \ $\mu_i$ \ is a Borel measure on
    \ $\cU_d$ \ satisfying
    \begin{align}\label{help_moment_mu}
      \int_{\cU_d}
       \biggl[\|\bz\| \land \|\bz\|^2
              + \sum_{j \in \{1, \ldots, d\} \setminus \{i\}} (1 \land z_j)\biggr]
       \mu_i(\dd\bz)
      < \infty .
    \end{align}
  \end{enumerate}
\end{Def}

Note that the measures \ $\nu$ \ and \ $\mu_i$, \ $i \in \{1, \ldots, d\}$, \ in Definition \ref{Def_admissible} are $\sigma$-finite,
 since the function \ $\cU_d\ni \br\mapsto 1 \land \Vert \br\Vert$ \ is strictly positive and measurable with a finite integral with respect to \ $\nu$,
 \ and, for each $i \in \{1, \ldots, d\}$, \ the function
 \ $\cU_d\ni \bz \mapsto \|\bz\| \land \|\bz\|^2 + \sum_{j \in \{1, \ldots, d\} \setminus \{i\}} (1 \land z_j)$ \
 is strictly positive and measurable with a finite integral with respect to \ $\mu_i$, \ see, e.g., Kallenberg \cite[Lemma 1.4]{K2}.

\begin{Thm}\label{CBI_exists}
Let \ $(d, \bc, \Bbeta, \bB, \nu, \bmu)$ \ be a set of admissible parameters.
Then there exists a unique conservative transition semigroup \ $(P_t)_{t\in\RR_+}$
 \ acting on the Banach space (endowed with the supremum norm) of real-valued
 bounded Borel measurable functions on the state space \ $\RR_+^d$ \ such that its
 Laplace transform has a representation
 \begin{equation*}%\label{Laplace_transform}
  \int_{\RR_+^d} \ee^{- \langle \blambda, \by \rangle} P_t(\bx, \dd \by)
  = \ee^{- \langle \bx, \bv(t, \blambda) \rangle
         - \int_0^t \psi(\bv(s, \blambda)) \, \dd s} , \qquad
  \bx \in \RR_+^d, \quad \blambda \in \RR_+^d , \quad t \in \RR_+ ,
 \end{equation*}
 where, for any \ $\blambda \in \RR_+^d$, \ the continuously differentiable
 function
 \ $\RR_+ \ni t \mapsto \bv(t, \blambda)
    = (v_1(t, \blambda), \ldots, v_d(t, \blambda)) \in \RR_+^d$
 \ is the unique locally bounded solution to the system 
 \begin{align}\label{help12}
  v_i(t, \blambda) = \lambda_i - \int_0^t \varphi_i(\bv(u, \blambda))\, \dd u ,  \qquad t\in\RR_+, \qquad i \in \{1, \ldots, d\} ,
 \end{align} 
 with
 \[
   \varphi_i(\blambda)
   := c_i \lambda_i^2 -  \langle \bB \be_i^{(d)}, \blambda \rangle
      + \int_{\cU_d}
         \bigl( \ee^{- \langle \blambda, \bz \rangle} - 1
                + \lambda_i (1 \land z_i) \bigr)
         \, \mu_i(\dd \bz)
 \]
 for \ $\blambda \in \RR_+^d$, \ $i \in \{1, \ldots, d\}$, \ and
 \begin{align}\label{dCBI_immigration}
   \psi(\blambda)
   := \langle \bbeta, \blambda \rangle
      + \int_{\cU_d}
         \bigl( 1 - \ee^{- \langle\blambda, \br\rangle} \bigr)
         \, \nu(\dd\br) , \qquad
   \blambda \in \RR_+^d .
 \end{align}
Further, the function $\RR_+\times \RR_+^d\ni(t,\blambda)\mapsto \bv(t, \blambda)$  is continuous.
\end{Thm}

Theorem \ref{CBI_exists} is a special case of Theorem 2.7 of Duffie et al.\ \cite{DufFilSch}
 with \ $m = d$, \ $n = 0$ \ and zero killing rate.
For the unique existence of a locally bounded solution to the system \eqref{help12}
 in Theorem \ref{CBI_exists}, see Li \cite[page 48]{Li} or Duffie et al.\ \cite[Proposition 6.4]{DufFilSch}.
The continuity of the function $\RR_+\times \RR_+^d\ni(t,\blambda)\mapsto \bv(t, \blambda)$ follows from 
 Proposition 6.4 in Duffie et al.\ \cite{DufFilSch}.
The system \eqref{help12} is usually written in differential form:
\begin{align*}
   \partial_1 v_i(t, \blambda) = - \varphi_i(\bv(t, \blambda)),  \qquad t\in\RR_+, \qquad
   v_i(0, \blambda) = \lambda_i , \qquad i \in \{1, \ldots, d\},
 \end{align*}
 and we follow this convention throughout the paper.

\begin{Def}\label{Def_CBI}
A conservative Markov process with state space \ $\RR_+^d$ \ and with transition
 semigroup \ $(P_t)_{t\in\RR_+}$ \ given in Theorem \ref{CBI_exists} is called a
 multi-type CBI process with parameters \ $(d, \bc, \Bbeta, \bB, \nu, \bmu)$.
\ The function
 \ $\RR_+^d \ni \blambda
    \mapsto \bvarphi(\blambda)= (\varphi_1(\blambda), \ldots, \varphi_d(\blambda)) \in \RR^d$
 \ is called its branching mechanism, and the function
 \ $\RR_+^d \ni \blambda \mapsto \psi(\blambda) \in \RR_+$ \ is called its
 immigration mechanism.
When there is no immigration, i.e., $\Bbeta = \bzero$ and $\nu = 0$,
 the process is simply called a multi-type CB process (a continuous state and continuous time branching process).
\end{Def}

Given a set of admissible parameters $(d, \bc, \Bbeta, \bB, \nu, \bmu)$, we get that
 $\nu$ and $\mu_i$, $i\in\{1,\ldots,d\}$, are L\'evy measures on $\cU_d$, since, by parts (v) and (vi) of Definition \ref{Def_admissible},
 we have
 \[
   \int_{\cU_d} (1 \land \|\br\|^2) \, \nu(\dd\br) \leq \int_{\cU_d}  (1 \land \|\br\|) \, \nu(\dd\br)< \infty,
 \]
 and
  \[
    \int_{\cU_d} (1\land \|\bz\|^2) \, \mu_i(\dd\bz)
       \leq \int_{\cU_d} (\|\bz\| \land \|\bz\|^2) \, \mu_i(\dd\bz)
      <\infty, \qquad i\in\{1,\ldots,d\}.
  \]
For this reason, we call $\nu+\sum_{i=1}^d \mu_i$ the total L\'evy measure corresponding to a multi-type CBI process
 with parameters $(d, \bc, \Bbeta, \bB, \nu, \bmu)$.

By Barczy et al.\ \cite[Remark 2.3 and (2.12)]{BarLiPap2}, for each $i\in\{1,\ldots,d\}$, the moment condition \eqref{help_moment_mu}
 is equivalent to
 \begin{align*}
      \int_{\cU_d}
       \biggl[\|\bz\| \land \|\bz\|^2
              + \sum_{j \in \{1, \ldots, d\} \setminus \{i\}} z_j\biggr]
       \mu_i(\dd\bz)
      < \infty ,
 \end{align*}
 which coincides with the moment condition in Example 2.5 in Li \cite[page 48]{Li} for multi-type CB processes,
 where multi-type CB processes are considered as special superprocesses.

It is known that a stochastically continuous multi-type CBI process $(\bX_t)_{t\in\RR_+}$ has a c\`{a}dl\`{a}g modification $(\widetilde\bX_t)_{t\in\RR_+}$
 (i.e.,  $\PP(\bX_t = \widetilde\bX_t)=1$, $t\in\RR_+$, and all the paths of $(\widetilde\bX_t)_{t\in\RR_+}$ are right continuous at every $t\in\RR_+$
 and possesses left limit at every $t\in\RR_{++}$), see, e.g., Keller-Ressel et al.\ \cite[Theorem 3.5]{KelSchTei} or He and Li \cite[page 1166]{HeLi}.
In what follows, it is always assumed that a multi-type CBI process has c\`{a}dl\`{a}g paths.

For a multi-type CBI process $(\bX_t)_{t\in\RR_+}$, $\bx\in\RR_+^d$, an event $A\in\sigma(\bX_t,t\in\RR_+)$ and an $\RR^p$-valued
 integrable random variable $\bxi$ which is $\sigma(\bX_t,t\in\RR_+)$-measurable (where $p\in\NN$), let $\PP_{\bx}(A):=\PP(A\mid \bX_0=\bx)$
 and $\EE_{\bx}(\bxi):=\EE(\bxi \mid \bX_0=\bx)$, respectively,
 where $\PP(A\mid \bX_0=\bx)$ is understood as $f(\bx)$, where $f$ is an appropriate function satisfying $\PP(A\mid \bX_0) = f(\bX_0)$ $\PP$-a.s.,
 and $\EE(\bxi \mid \bX_0=\bx)$ is understood in a similar way.
  
Let \ $(\bX_t)_{t\in\RR_+}$ \ be a multi-type CBI process with parameters
 \ $(d, \bc, \Bbeta, \bB, \nu, \bmu)$ \ such that the moment condition
 \begin{equation}\label{moment_condition_m_new}
  \int_{\cU_d} \|\br\| \bbone_{\{\|\br\|\geq1\}} \, \nu(\dd\br) < \infty
 \end{equation}
 holds.
Then, by formula (3.4) in Barczy et al.\ \cite{BarLiPap2} (see also formula (79) in Li \cite{Li3}),
 \begin{equation*}%\label{EXcond}
  \EE_{\bx}(\bX_t)
  = \ee^{t\tbB} \bx + \int_0^t \ee^{u\tbB} \tBbeta \, \dd u ,
  \qquad \bx \in \RR_+^d , \quad t \in \RR_+ ,
 \end{equation*}
 where
 \begin{gather}\label{help5}
  \tbB := (\tb_{i,j})_{i,j\in\{1,\ldots,d\}} , \qquad
  \tb_{i,j}
  := b_{i,j} + \int_{\cU_d} (z_i - \delta_{i,j})^+ \, \mu_j(\dd\bz) , \qquad
  \tBbeta := \Bbeta + \int_{\cU_d} \br \, \nu(\dd\br) ,
 \end{gather}
 with \ $\delta_{i,j}:=1$ \ if \ $i = j$, \ and \ $\delta_{i,j} := 0$ \ if
 \ $i \ne j$, where $i,j\in\{1,\ldots,d\}$.
\ Note that, for all \ $\bx \in \RR^d_+$, \ the function
 \ $\RR_+ \ni t \mapsto \EE_{\bx}(\bX_t)$ \ is continuous, and
 \ $\tbB \in \RR^{d \times d}_{(+)}$ \ and \ $\tBbeta \in \RR_+^d$.
Indeed, part (v) of Definition \ref{Def_admissible} together with the moment condition \eqref{moment_condition_m_new}
 and Barczy et al.\ \cite[Remark 2.3 and formulas (2.11) and (2.12)]{BarLiPap2} yield that
 \[
   \int_{\cU_d} \|\br\| \, \nu(\dd\br) < \infty , \qquad
   \int_{\cU_d} (z_i - \delta_{i,j})^+ \, \mu_j(\dd \bz) < \infty , \quad
   i, j \in \{1, \ldots, d\} .
 \]
Here we point out the fact that the matrix $\tbB$ belongs to $\RR^{d \times d}_{(+)}$ even if the moment condition
 \eqref{moment_condition_m_new} does not hold, however, the vector $\tBbeta$ belongs to $\RR_+^d$ if and only if
 the moment condition \eqref{moment_condition_m_new} holds.

Given a set of admissible parameters $(d, \bc, \Bbeta, \bB, \nu, \bmu)$ such that the moment condition
 \eqref{moment_condition_m_new} holds, let us consider the SDE
 \begin{align}\label{SDE_atirasa_dimd}
  \begin{split}
   \bX_t
   &=\bX_0
     + \int_0^t (\Bbeta + \tbB \bX_u) \, \dd u
     + \sum_{\ell=1}^d
        \int_0^t \sqrt{2 c_\ell X_{u,\ell}^+} \, \dd W_{u,\ell}
        \, \be_\ell^{(d)} \\
   &\quad
      + \sum_{\ell=1}^d
         \int_0^t \int_{\cU_d} \int_{\cU_1}
          \bz \bbone_{\{w\leq X_{u-,\ell}\}} \, \tN_\ell(\dd u, \dd\bz, \dd w)
      + \int_0^t \int_{\cU_d} \br \, M(\dd u, \dd\br)
  \end{split}
 \end{align}
 for $t \in\RR_+$, where
 $X_{t,\ell}$, $\ell \in \{1, \ldots, d\}$, denotes the $\ell^{\mathrm{th}}$ coordinate of $\bX_t$,
  $(W_{t,\ell})_{t\in\RR_+}$, $\ell\in\{1,\ldots,d\}$, are standard Wiener processes,
 \ $N_\ell$, \ $\ell \in \{1, \ldots, d\}$, \ and \ $M$ \ are Poisson
 random measures on \ $\cU_1 \times \cU_d \times \cU_1$ \ and on
 \ $\cU_1 \times \cU_d$ \ with intensity measures
 \ $\dd u \, \mu_\ell(\dd\bz) \, \dd w$, \ $\ell \in \{1, \ldots, d\}$, \ and
 \ $\dd u \, \nu(\dd\br)$, \ respectively, and
 \ $\tN_\ell(\dd u, \dd\bz, \dd w)
    := N_\ell(\dd u, \dd\bz, \dd w) - \dd u \, \mu_\ell(\dd\bz) \, \dd w$,
 \ $\ell \in \{1, \ldots, d\}$.
\ We suppose that $\EE(\Vert \bX_0\Vert)<\infty$ and that \ $\bX_0$, \ $(W_{t,1})_{t\in\RR_+}$, \ldots, \ $(W_{t,d})_{t\in\RR_+}$, \ $N_1$,
 \ldots, \ $N_d$ \ and \ $M$ \ are mutually independent.
The SDE \eqref{SDE_atirasa_dimd} has a pathwise unique $\RR_+^d$-valued c\`{a}dl\`{a}g strong solution,
 and the solution is a multi-type CBI process with parameters \ $(d, \bc, \Bbeta, \bB, \nu, \bmu)$,
 \ see Theorem 4.6 and Section 5 in
 Barczy et al.~\cite{BarLiPap2}, where \eqref{SDE_atirasa_dimd} was proved only for
 \ $d \in \{1, 2\}$, \ but their method clearly works for all \ $d \in\NN$.
Consequently, given a c\`{a}dl\`{a}g multi-type CBI process $(\bX_t)_{t\in\RR_+}$ with parameters $(d, \bc, \Bbeta, \bB, \nu, \bmu)$
 such that $\EE(\Vert \bX_0\Vert)<\infty$ and the moment condition \eqref{moment_condition_m_new} hold,
 its law on the space of $\RR^d$-valued c\`{a}dl\`{a}g functions defined on $\RR_+$ coincides with
 the law of the pathwise unique c\`{a}dl\`{a}g strong solution of the SDE \eqref{SDE_atirasa_dimd}.
In the remaining part of the paper, when we refer to a multi-type CBI process $(\bX_t)_{t\in\RR_+}$ with parameters \ $(d, \bc, \Bbeta, \bB, \nu, \bmu)$
 such that $\EE(\Vert \bX_0\Vert)<\infty$ and the moment condition \eqref{moment_condition_m_new} hold, we consider it
 as a pathwise unique c\`{a}dl\`{a}g strong solution of the SDE \eqref{SDE_atirasa_dimd}.

Recall that \ $\cB(\cU_d)$ \ denotes the set of Borel sets of \ $\cU_d$.
\ For all \ $A \in\cB(\cU_d)$, \ let
\begin{align}\label{tau_def}
  \tau_A
    := \inf\{u \in \RR_{++}
              : \Delta \bX_u \in A\}
 \end{align}
 with the convention \ $\inf(\emptyset) := \infty$, \ and 
 \begin{align*}
  J_t(A) := \begin{cases}
             \card(\{u \in (0, t] : \Delta \bX_u  \in A\}) & \text{if $t\in\RR_{++}$,}\\
             0 & \text{if $t=0$,}
             \end{cases}
 \end{align*}
 where \ $\Delta \bX_u := \bX_u - \bX_{u-}$, \ $u \in \RR_{++}$, \
 and \ $\card(H)$ \ denotes the cardinality of a set \ $H$.

In the forthcoming results, given a set \ $A\in\cB(\cU_d)$, \ the condition that \ $\nu(A)+\sum_{\ell=1}^d \mu_\ell(A)<\infty$
 \ will come into play.
In the next remark, we give a sufficient condition under which it holds.

\begin{Rem}\label{Rem_A_condition}
For any $\vare>0$, let  $K_\vare := \{\by\in\RR_+^d : \Vert \by\Vert<\vare\}$.
If \ $A\in\cB(\cU_d)$ \ is such that there exists \ $\vare\in(0,1)$ \ with \ $A\subseteq \RR_+^d\setminus K_\vare$
 \ (or equivalently, \ $A\in\cB(\cU_d)$ \ is such that \ $\bzero$ \ is not contained in the closure of $A$), \ then
 \ $\nu(A)+\sum_{\ell=1}^d \mu_\ell(A)<\infty$.
\ Indeed, if \ $A\in \cB(\cU_d)$ \ is such that \ $A\subseteq \RR_+^d\setminus K_\vare$ \ with some $\vare\in(0,1)$,
 \ then, by parts (v) and (vi) of Definition \ref{Def_admissible}, we have
 \begin{align*}
    \nu(A) \leq \nu(\RR_+^d\setminus K_\vare)
           \leq \frac{1}{\vare} \int_{\{ \br\in\cU_d\,:\, \vare\leq \Vert \br\Vert \leq 1 \}} \Vert\br\Vert \,\nu(\dd\br)
             + \int_{\{ \br\in\cU_d\,:\, \Vert \br\Vert >1 \}}1 \,\nu(\dd\br)
           <\infty,
 \end{align*}
 and, for each \ $\ell\in\{1,\ldots,d\}$,
 \begin{align*}
    \mu_\ell(A) \leq \mu_\ell(\RR_+^d\setminus K_\vare)
           \leq \frac{1}{\vare^2} \int_{\{ \bz\in\cU_d\,:\, \vare\leq \Vert \bz\Vert \leq 1 \}} \Vert\bz\Vert^2 \,\mu_\ell(\dd\bz)
             + \int_{\{ \bz\in\cU_d\,:\, \Vert \bz\Vert >1 \}} \Vert \bz\Vert \,\mu_\ell(\dd\bz)
           <\infty,
 \end{align*}
 as desired.

Note also that, given Borel measures $\mu$ and $\nu$ on $(0,\infty)$, $k\in\NN$ 
 and $A_1,\ldots,A_k\in\cB((0,\infty))$, the condition $\nu(A_i)+\mu(A_i) <\infty$, $i\in\{1,\ldots,k\}$,
 is equivalent to $\nu\left( \bigcup_{i=1}^k A_i\right) + \mu\left( \bigcup_{i=1}^k A_i\right) <\infty$,
 since the measures $\mu$ and $\nu$ are $\sigma$-subadditive. 
\proofend
\end{Rem}

Given a set of admissible parameters \ $(d, \bc, \Bbeta, \bB, \nu, \bmu)$ \ and a Borel set
 $A\in\cB(\cU_d)$ such that $\nu(A)+\sum_{\ell=1}^d \mu_\ell(A)<\infty$, let us introduce
 \ $\bB^{(A)}\in\RR^{d\times d}_{(+)}$, \ $\nu^{(A^{-1})}$ \ and $\bmu^{(A^{-1})}$ \ by
  \begin{align}\label{mod_parameters}
   \begin{split}
             &(\bB^{(A)})_{i,j}
                 := b_{i,j} + \delta_{i,j}\int_A \big( (z_i - 1)^+ - z_i \big)\,\mu_i(\dd\bz),\qquad i,j\in\{1,\ldots,d\},\\
             &\nu^{( A^{-1})}:=\nu\vert_{A^{-1}},\\
             &\bmu^{(A^{-1})}:=\Big(\mu^{(A^{-1})}_1,\ldots,\mu^{(A^{-1})}_d\Big)  \quad \text{with}\quad
               \mu^{(A^{-1})}_\ell:=\mu_\ell\vert_{A^{-1}}, \quad \ell\in\{1,\ldots,d\},
   \end{split}
  \end{align}
  where we recall $A^{-1}=\cU_d\setminus A$.
Then one can easily see that \ $\big(d, \bc, \Bbeta, \bB^{(A)}, \nu^{(A^{-1})}, \bmu^{(A^{-1})}\big)$ \
 is a set of admissible parameters.
Indeed, for each $i\in\{1,\ldots,d\}$, we have 
 \[
  \int_A \vert (z_i - 1)^+ - z_i \vert \,\mu_i(\dd\bz)
    = \int_A  z_i \bbone_{\{z_i\leq1\}}\,\mu_i(\dd\bz) + \int_A  \bbone_{\{z_i>1\}}\,\mu_i(\dd\bz) \leq \mu_i(A)<\infty.
 \]  
If $(\bX_t)_{t\in\RR_+}$ is a  multi-type CBI process with parameters \ $(d, \bc, \Bbeta, \bB, \nu, \bmu)$, \
 $A\in\cB(\cU_d)$ is such that $\nu(A)+\sum_{\ell=1}^d \mu_\ell(A)<\infty$, and 
 \ $(\bX^{(A)}_t)_{t\in\RR_+}$ \ is a multi-type CBI process with parameters
 \ $\big(d, \bc, \Bbeta, \bB^{(A)}, \nu^{(A^{-1})}, \bmu^{(A^{-1})}\big)$ \ such that \ $\bX_0^{(A)} = \bX_0$, then,
 intuitively, the process $(\bX^{(A)}_t)_{t\in\RR_+}$ can be obtained by removing from $(\bX_t)_{t\in\RR_+}$ all the
 masses produced by the jumps with size vectors in the set $A$.
This argument has been made precise in mathematical terms in the proof of Theorem 4.2 in Barczy and Palau \cite{BarPal}.
Note that if the moment condition \eqref{moment_condition_m_new} holds for $(\bX_t)_{t\in\RR_+}$, then it also holds for
 $(\bX^{(A)}_t)_{t\in\RR_+}$, since
 \begin{align*}%\label{moment_condition_m_new_A}
  \int_{\cU_d} \|\br\| \bbone_{\{\|\br\|\geq1\}} \, \nu^{(A^{-1})}(\dd\br)
     =  \int_{A^{-1}} \|\br\| \bbone_{\{\|\br\|\geq1\}} \, \nu(\dd\br) <\infty.
 \end{align*}
Barczy and Palau \cite[Theorem 4.2]{BarPal} proved that 
 the branching mechanism $\bvarphi^{(A)}$ and the immigration mechanism $\psi^{(A)}$ of $(\bX^{(A)}_t)_{t\in\RR_+}$ take 
 the following forms
 $\RR_+^d\ni \blambda \mapsto \bvarphi^{(A)}(\blambda):= (\varphi^{(A)}_1(\blambda),\ldots,\varphi^{(A)}_d(\blambda)) \in\RR^d$ with
 \begin{align}\label{branching_X_A}
  \varphi^{(A)}_\ell(\blambda)
                    := \varphi_\ell(\blambda) + \int_A  (1  - \ee^{-\langle \blambda, \bz \rangle}) \, \mu_\ell(\dd \bz),
                      \qquad \blambda\in\RR_+^d, \qquad \ell\in\{1,\ldots,d\},
 \end{align}
 and  
 \begin{align}\label{immigration_X_A}
  \RR_+^d \ni \blambda \mapsto
   \psi^{(A)}(\blambda)
    : = \langle \bbeta, \blambda \rangle
     + \int_{A^{-1}}
     \bigl( 1 - \ee^{- \langle\blambda, \br\rangle} \bigr)
     \, \nu(\dd\br)\in \RR_+,
 \end{align}
 respectively.
Concerning our notations, we note that \ $(A)$ \ in the superscript of a formula (e.g. in that of \ $\varphi^{(A)}_\ell$) \ means 
 only that the corresponding expression depends on \ $A$.
Nonetheless, the restrictions of the measures \ $\nu$ \ and \ $\bmu$ \ onto \ $A^{-1}$
 \ are denoted by \ $\nu^{(A^{-1})}$ \ and \ $\bmu^{(A^{-1})}$, \ respectively,
 in order to avoid some possible confusion.  
 
The next result was also showed in Barczy and Palau \cite[Theorem 4.2]{BarPal}.

\begin{Thm}\label{Thm_jump_times_dCBI}
Let \ $(\bX_t)_{t\in\RR_+}$ \ be a multi-type CBI process with parameters \  $(d, \bc, \Bbeta, \bB, \nu, \bmu)$ \ such that
 \ $\EE(\|\bX_0\|) < \infty$ \ and the moment condition \eqref{moment_condition_m_new} hold.
Then for all \ $A\in\cB(\cU_d)$ \ such that
 \ $\nu(A)+\sum_{\ell=1}^d \mu_\ell(A)<\infty$, \ we have
  \begin{itemize}
    \item[(i)] $\EE(J_t(A))<\infty$, \ $t \in \RR_{++}$, \ which, in particular, implies that
                \ $\PP(J_t(A) < \infty) = 1$, \ $t\in\RR_{++}$;
    \item[(ii)] for all \ $t\in\RR_+$ \ and \ $\bx=(x_1,\ldots,x_d)\in\RR_+^d$,
       \begin{align*}%\label{formula3}
         \PP_{\bx}(\tau_A>t)
          = \exp\left\{ -\nu(A)t - \sum_{\ell=1}^d x_\ell\,\tv_\ell^{(A)}(t,\bmu(A))
                           - \int_0^t \psi^{(A)}\big(\tbv^{(A)}(s,\bmu(A))\big)\,\dd s \right\},
        \end{align*}
        where
          \begin{itemize}
           \item[$\bullet$] $\bmu(A):=(\mu_1(A),\ldots,\mu_d(A))$,
           \item[$\bullet$] $\psi^{(A)}$  is the immigration mechanism of $(\bX^{(A)}_t)_{t\in\RR_+}$ given in \eqref{immigration_X_A},
           \item[$\bullet$] the continuously differentiable function
\begin{align*}%\label{tvA}
\RR_+\ni t \mapsto \tbv^{(A)}(t,\bmu(A)):= \big(\tv^{(A)}_1(t,\bmu(A)),\ldots,\tv^{(A)}_d(t,\bmu(A))\big) \in\RR_+^d
\end{align*}
is the unique locally bounded solution to the system of differential equations
\begin{align}\label{help_DE_tv}
 \begin{split}
 &\partial_1 \tv_i^{(A)}(t, \bmu(A)) = \mu_i(A) - \varphi_i^{(A)}( \tbv^{(A)}(t, \bmu(A)) ), \qquad i \in \{1, \ldots, d\} ,\\
                    & \tv_i^{(A)}(0, \bmu(A)) = 0 , \qquad i \in \{1, \ldots, d\},
 \end{split}
 \end{align}
    where $\bvarphi^{(A)}= (\varphi^{(A)}_1,\ldots,\varphi^{(A)}_d)$ is the branching mechanism of $(\bX^{(A)}_t)_{t\in\RR_+}$ given in \eqref{branching_X_A}.
          \end{itemize}
     \end{itemize}
\end{Thm}

In the next remark, we recall some ingredients of the proof of Theorem 4.2 in Barczy and Palau \cite{BarPal}
 that will be used later on.

\begin{Rem}\label{REM_Prelimin}
In Step 2 of the proof of Theorem 4.2 in Barczy and Palau \cite{BarPal}, we showed that
 \begin{align}\label{help4}
  J_t(A)
  = \sum_{\ell=1}^d
      \int_0^t \int_{\cU_d} \int_{\cU_1} \bbone_A(\bz)
       \bbone_{\{w\leq X_{u-,\ell}\}}
       \, N_\ell(\dd u, \dd\bz, \dd w)
    + \int_0^t \int_{\cU_d} \bbone_A(\br) \, M(\dd u, \dd\br)
 \end{align}
 almost surely for all $t\in\RR_{++}$.
Since $J_0(A)$ is defined to be $0$, \eqref{help4} holds for $t=0$ as well.
\proofend
\end{Rem}

In the next remark, in case of $d=1$,
 we rewrite the branching mechanism $\varphi^{(A)}$ and the immigration mechanism $\psi^{(A)}$
 defined in Theorem \ref{Thm_jump_times_dCBI}.

\begin{Rem}\label{rem:VarphiA}
Let  $(1, c, \beta, B, \nu, \mu)$ be a set of admissible parameters such that the moment condition \eqref{moment_condition_m_new} with $d=1$ holds.
Let $A\in\cB((0,\infty))$ be such that $\nu(A) + \mu(A)<\infty$. If $(X^{(A)}_t)_{t\in\RR_+}$ is a CBI process with parameters
 $\big(1, c, \beta, B^{(A)}, \nu^{(A^{-1})}, \mu^{(A^{-1})}\big)$, where $B^{(A)}, \nu^{(A^{-1})}$ and $\mu^{(A^{-1})}$ are given \eqref{mod_parameters},
 then its branching mechanism $\varphi^{(A)}:\RR_+\to\RR$ and immigration mechanism $\psi^{(A)}:\RR_+\to\RR_+$ take the form
\begin{align}\label{help13}
  \begin{split}
    \psi^{(A)}(\lambda)  =
 \beta\lambda + \int_{\cU_1\setminus A}(1-\ee^{-\lambda r})\,\nu(\dd r),
        \qquad \lambda\in\RR_+,
  \end{split}
 \end{align}
and 
 \begin{align}\label{eq:psiA}
   \begin{split}
   \varphi^{(A)}(\lambda)
    & = \varphi(\lambda) + \int_{A} ( 1 - \ee^{-\lambda z} )\,\mu(\dd z)
      = \varphi(\lambda) + \mu(A) - \int_{A} \ee^{-\lambda z} \,\mu(\dd z) \\ 
    &= c\lambda^2 - B\lambda
       + \int_{\cU_1} ( \ee^{-\lambda z} - 1 + \lambda(1\wedge z) )\,\mu(\dd z)
       + \int_{A} ( 1 - \ee^{-\lambda z} )\,\mu(\dd z) \\
    & = c\lambda^2 - \left( B - \int_{A}(1\wedge z)\,\mu(\dd z)\right)\lambda
        + \int_{\cU_1\setminus A} ( \ee^{-\lambda z} - 1 + \lambda(1\wedge z) )\,\mu(\dd z)
  \end{split}
 \end{align}
 for all $\lambda\in\RR_+$, where we used that $\int_{A}(1\wedge z)\,\mu(\dd z)\leq \int_{A}1\,\mu(\dd z)=\mu(A)<\infty$.
\proofend 
\end{Rem}

\section{Joint Laplace transform of a CBI process and its first jump times with jump sizes in given Borel sets}\label{Section_CBI_results_1}
 
First, we recall the notations $A^{-1} = \Omega\setminus A$ for any subset $A$ of $(0,\infty)$ and
 $\bA^{\bell}  = A_1^{\ell_1} \cap \dots \cap A_k^{\ell_k}$ for any $k\in\NN$, subsets $A_1,\ldots,A_k$ of $(0,\infty)$,
 and $\bell=(\ell_1,\dots,\ell_k) \in \{-1,1\}^k$.
Note also that, given $z\in(0,\infty)$ and $\ell\in\{-1,1\}$, the inclusion $z\in A^{\ell}$ holds
 if and only if $\bbone_{A}(z)=\frac{1}{2}(\ell+1)$.
Further, let $\bone$ denote the $k$-dimensional vector having all the coordinates $1$.

The next Theorem \ref{Thm_(k+1)CBI} is a generalization of Theorem 10.11 in Li \cite{Li}, namely,
 in the special case $k=1$, Theorem \ref{Thm_(k+1)CBI} gives back Theorem 10.11 in Li \cite{Li}.
We also point out the fact that our proof of Theorem \ref{Thm_(k+1)CBI} is slightly different from that of 
 Theorem 10.11 in Li \cite{Li}, we directly use an SDE representation of multi-type CBI processes (see \eqref{SDE_atirasa_dimd}).

\begin{Thm}\label{Thm_(k+1)CBI}
Let $(X_t)_{t\in\RR_+}$ be a CBI process with parameters $(1, c, \beta, B, \nu, \mu)$ such that
 $\EE(X_0) < \infty$ and the moment condition \eqref{moment_condition_m_new} with $d=1$ hold.
Let $k\in\NN$ and $A_1,\ldots,A_k\in\cB((0,\infty))$ be such that $\mu(A_i) + \nu(A_i)<\infty$, $i\in\{1,\ldots,k\}$.
Then $(X_t,J_t(A_1),\ldots,J_t(A_k))_{t\in\RR_+}$ is a multi-type CBI process with parameters
 \ $(k+1, \bc^{(k)}, \Bbeta^{(k)}, \bB^{(k)}, \nu^{(k)}, \bmu^{(k)})$ starting from $(X_0,\bzero)$, where
\begin{itemize}
\item $\bc^{(k)}:=c\, \be_1^{(k+1)}\in\RR_+^{k+1}$, \ $\Bbeta^{(k)}:= \beta\, \be_1^{(k+1)}\in\RR_+^{k+1}$,
\item $\bB^{(k)}:= 
     (b^{(k)}_{i,j})_{i,j\in\{1,\ldots,k+1\}}\in\RR^{(k+1)\times(k+1)}_{(+)}$, where $b^{(k)}_{i,j}:= B\bbone_{\{i=j=1\}}$,
 \item $\nu^{(k)}$ is a  Borel measure on $\cU_{k+1}$ defined by
 \begin{align}\label{def_nu1}
  \begin{split}
  \nu^{(k)}(B)
  := \big( (\mathrm{id}_{\cU_1}, \bbone_{A_1},\ldots,\bbone_{A_k})_{*}(\nu)\big)(B)
  =\sum_{\bell\in\{-1,1\}^k}\int_{\bA^{\bell}}\bbone_{B}(r,\tfrac{1}{2}(\bell+\bone))\,\nu(\dd r)
  \end{split}
 \end{align}
 for $B\in\cB(\cU_{k+1})$, where $(\mathrm{id}_{\cU_1}, \bbone_{A_1},\ldots,\bbone_{A_k})_{*}(\nu)$ denotes the pushforward measure of $\nu$
  by the mapping $(\mathrm{id}_{\cU_1}, \bbone_{A_1},\ldots,\bbone_{A_k})$ defined on $\cU_1$ with range in $\cU_{k+1}$,
 \item $\bmu^{(k)}:=(\mu_1^{(k)},\dots, \mu_{k+1}^{(k)})$ is a vector of Borel measures on $\cU_{k+1}$ defined by
       \begin{align}\label{def_mu1}
       \begin{split}
         \mu_1^{(k)}(B)
         := \big( (\mathrm{id}_{\cU_1}, \bbone_{A_1},\ldots,\bbone_{A_k})_{*}(\mu)\big)(B)
         =\sum_{\bell\in\{-1,1\}^k}\int_{\bA^{\bell}}\bbone_{B}(z,\tfrac{1}{2}(\bell+\bone))\,\mu(\dd z)
        \end{split} 
       \end{align}
      for $B\in\cB(\cU_{k+1})$, and $\mu_i^{(k)}(B):=0$, $B\in\cB(\cU_{k+1})$, $i\in\{2,\dots,k+1\}$.
 \end{itemize}
\end{Thm}

\noindent{\bf Proof.}
We will use some results with $d=1$  recalled in Section \ref{section_CBI} for the single-type CBI process $(X_t)_{t\in\RR_+}$,
 and, for simplicity, in this case, in the SDE representation \eqref{SDE_atirasa_dimd},
 instead of $(W_{t,1})_{t\in\RR_+}$ and $N_1$ we will write $(W_t)_{t\in\RR_+}$ and $N$, respectively.
The (second) equalities in \eqref{def_nu1} and \eqref{def_mu1} follow from the definition of pushforward measure.
Indeed, in case of $\nu^{(k)}$, for all $B\in\cB(\cU_{k+1})$, using that  $\{\bA^\bell : \bell\in\{-1,1\}^k\}$ is a partition of $\cU_1$,
 we get that
 \begin{align*}
   \nu^{(k)}(B)
   & = \nu\big( (\mathrm{id}_{\cU_1}, \bbone_{A_1},\ldots,\bbone_{A_k})^{-1}(B) \big)
      = \nu\big( \{ z\in \cU_1 : (\mathrm{id}_{\cU_1}(z),\bbone_{A_1}(z),\ldots,\bbone_{A_k}(z))\in B \} \big) 
\end{align*}      
\begin{align*}
    & = \sum_{\bell\in\{-1,1\}^k} \nu\big( \{ z\in \bA^{\bell} : (z,\tfrac{1}{2}(\bell+\bone))\in B \} \big)  
      = \sum_{\bell\in\{-1,1\}^k} \int_{\cU_1} \bbone_{\left\{ z \in  \bA^{\bell} : \left( z, \tfrac{1}{2} (\bell+\bone) \right) \in B \right\} }\,\nu(\dd z) \\
    & = \sum_{ \bell \in \{-1,1\}^k } \int_{\bA^{\bell}} \bbone_{B}(z, \tfrac{1}{2}(\bell+\bone))\,\nu(\dd z)  
     = \nu_1^{(k)}(B),
  \end{align*}
 as desired.

{\sl Step 1.}
By Remark \ref{REM_Prelimin}, for each $i\in\{1,\ldots,k\}$, we have that
 \begin{align*}
  J_t(A_i)
  = \int_0^t \int_{\cU_1} \int_{\cU_1} \bbone_{A_i}(z)
       \bbone_{\{w\leq X_{u-}\}}
       \, N(\dd u, \dd z, \dd w)
    + \int_0^t \int_{\cU_1} \bbone_{A_i}(r) \, M(\dd u, \dd r)
 \end{align*}
 almost surely for all $t\in\RR_+$.
We rewrite the first term on the right hand side of the previous equation in a way that it becomes the sum of $\mu(A_i) \int_0^t X_u\,\dd u$ and
 an integral with respect to the compensated Poisson random measure $\tN$.
Since $\mu(A_i)<\infty$, $i\in\{1,\ldots,k\}$, and $(X_t)_{t\in\RR_+}$ has c\`{a}dl\`{a}g paths almost surely,
 for each $i\in\{1,\ldots,k\}$, we get
 \begin{align*}
   &\int_0^t \int_{\cU_1} \int_{\cU_1} \bbone_{A_i}(z)
       \bbone_{\{w\leq X_{u-}\}}
       \, N(\dd u, \dd z, \dd w) \\
  &\qquad = \int_0^t \int_{\cU_1} \int_{\cU_1} \bbone_{A_i}(z)
       \bbone_{\{w\leq X_{u-}\}}
       \,\tN(\dd u, \dd z, \dd w)
      + \int_0^t \int_{\cU_1} \int_{\cU_1} \bbone_{A_i}(z)
        \bbone_{\{w\leq X_{u-}\}}
        \, \dd u \,\mu(\dd z)\,\dd w \\
  &\qquad = \int_0^t \int_{\cU_1} \int_{\cU_1} \bbone_{A_i}(z)
       \bbone_{\{w\leq X_{u-}\}}
       \,\tN(\dd u, \dd z, \dd w)
       + \mu(A_i) \int_0^t X_u\,\dd u
   \end{align*}
 almost surely for all $t\in\RR_+$.
Consequently,  for each $i\in\{1,\ldots,k\}$,  we have that
 \begin{align*}
  J_t(A_i)
  = \int_0^t \int_{\cU_1} \int_{\cU_1} \bbone_{A_i}(z)
       \bbone_{\{w\leq X_{u-}\}}
       \,\tN(\dd u, \dd z, \dd w)
    + \int_0^t \int_{\cU_1} \bbone_{A_i}(r) \, M(\dd u, \dd r)
   + \mu(A_i) \int_0^t X_u\,\dd u
 \end{align*}
 almost surely for all $t\in\RR_+$.

Using the SDE \eqref{SDE_atirasa_dimd} for $(X_t)_{t\in\RR_+}$, we have
 \begin{align}\label{help_SDE_joint}
 \begin{split}
  \begin{bmatrix}
    X_t \\
    J_t(A_1) \\
\vdots\\
J_t(A_k)\\
  \end{bmatrix}
  & = \begin{bmatrix}
        X_0 \\
        0 \\
       \vdots\\
       0 \\
     \end{bmatrix}
     + \int_0^t \left(
                      \begin{bmatrix}
                          \beta \\
                          0 \\
                          \vdots\\
                          0\\
                       \end{bmatrix}
                     + \begin{bmatrix}
                         \tB & 0& \cdots& 0 \\
                         \mu(A_1) & 0& \cdots& 0 \\
                         \vdots & \vdots& \ddots& \vdots \\
                         \mu(A_k) & 0& \cdots& 0 \\
                       \end{bmatrix}
                       \begin{bmatrix}
                         X_u \\
                         J_u(A_1) \\
                        \vdots\\
                        J_u(A_k)\\
                       \end{bmatrix}
                      \right) \, \dd u \\
  &\phantom{=\;} + 
    \begin{bmatrix}
        \int_0^t \sqrt{2 c X_u^+ } \, \dd W_u &
         0 &
        \cdots &
         0 \\
       \end{bmatrix}^\top\\
  &\phantom{=\;} + \begin{bmatrix}
              \int_0^t \int_{\cU_1} \int_{\cU_1}
                   z \bbone_{\{w\leq X_{u-}\}} \, \tN(\dd u, \dd z, \dd w) \\
               \int_0^t \int_{\cU_1} \int_{\cU_1} \bbone_{A_1}(z) \bbone_{\{w\leq X_{u-}\}} \, \tN(\dd u, \dd z, \dd w)  \\
               \vdots\\
               \int_0^t \int_{\cU_1} \int_{\cU_1} \bbone_{A_k}(z) \bbone_{\{w\leq X_{u-}\}} \, \tN(\dd u, \dd z, \dd w)  \\
           \end{bmatrix}
         + \begin{bmatrix}
             \int_0^t \int_{\cU_1} r \, M(\dd u, \dd r) \\
             \int_0^t \int_{\cU_1} \bbone_{A_1}(r) \, M(\dd u, \dd r)  \\
               \vdots\\
             \int_0^t \int_{\cU_1} \bbone_{A_k}(r) \, M(\dd u, \dd r)  \\
           \end{bmatrix}
 \end{split}
 \end{align}
 for $t\in\RR_+$.

{\sl Step 2.}
Let $(W^{(1)}_t)_{t\in\RR_+}, \dots, (W^{(k)}_t)_{t\in\RR_+}$ be independent  standard Wiener processes such that they are
 independent of $(W_t)_{t\in\RR_+}$, $N$ and $M$.

 Let $h:\cU_1\times\cU_1\times\cU_1 \to \cU_1\times\cU_{k+1}\times\cU_1$,%
 \[
   h(u,z,w):=(u,z,\bbone_{A_1}(z),\ldots,\bbone_{A_k}(z),w), \qquad (u,z,w)\in \cU_1\times\cU_1\times\cU_1.
 \]
Here the range of $h$ is indeed a subset of $\cU_1\times\cU_{k+1}\times\cU_1$, since $(z,\bbone_{A_1}(z),\ldots,\bbone_{A_k}(z))\neq \bzero$ for all $z\in\cU_1$ and
 $(z,\bbone_{A_1}(z),\ldots,\bbone_{A_k}(z))\in\cU_1\times \{0,1\}^k\subset \cU_{k+1}$.
Let $N^{(1)}$ be the random point measure $N\circ h^{-1}$, which is concentrated on $\cU_1\times\cU_{k+1}\times\cU_1$.
We check that $N^{(1)}$ is a Poisson random measure on $\cU_1\times\cU_{k+1}\times\cU_1$ with intensity measure
 $\dd u\,\mu_1^{(k)}(\dd \bz)\,\dd w$, where the Borel measure $\mu_1^{(k)}$ on $\cU_{k+1}$ is given by \eqref{def_mu1}.
By the mapping theorem for Poisson point processes (see, e.g., Kingman \cite[Section 2.3]{Kin}),
 it is enough to check that
 $(\lambda_1\otimes \mu\otimes\lambda_1)\circ h^{-1} = \lambda_1\otimes \mu_1^{(k)}\otimes\lambda_1$,
 where $\lambda_1$ denotes the Lebesgue measure on $\cU_1=(0,\infty)$.
For all $C\in \cB(\cU_1)$, $B\in\cB(\cU_{k+1})$ and $D\in\cB(\cU_1)$, we have that
 \begin{align*}
 (&(\lambda_1\otimes \mu\otimes\lambda_1)\circ h^{-1})(C\times B\times D)
   = (\lambda_1\otimes \mu\otimes\lambda_1) (h^{-1}(C\times B\times D))\\
  &\quad = (\lambda_1\otimes \mu\otimes\lambda_1)
     \big( \{ (u,z,w)\in \cU_1\times\cU_1\times\cU_1 : u\in C, (z,\bbone_{A_1}(z),\ldots,\bbone_{A_k}(z))\in B, w\in D  \} \big)\\
  &\quad = (\lambda_1\otimes \mu\otimes\lambda_1)
     \big( C \times \{z\in\cU_1 : (z,\bbone_{A_1}(z),\ldots,\bbone_{A_k}(z))\in B\} \times D   \big)\\
  &\quad = \lambda_1(C)  \mu\big( \{ z\in \cU_1 : (z,\bbone_{A_1}(z),\ldots,\bbone_{A_k}(z))\in B \} \big) \lambda_1(D)\\
  &\quad = \lambda_1(C) \mu_1^{(k)}(B)  \lambda_1(D),
 \end{align*}
 as desired, where the last equality is due to the definition of pushforward meausure.

Let $N^{(2)}, \dots , N^{(k+1)}$ be zero Poisson random measures on $ \cU_1\times\cU_{k+1}\times\cU_1$
 with intensity measure $\dd u \,\mu_i^{(k)}(\dd\bz)\,\dd w$, where $\mu_i^{(k)}(B):=0$, $B\in \cB(\cU_{k+1})$, for any $i\in\{2,\dots,k+1\}$.

Let $\widetilde h:\cU_1\times\cU_1 \to \cU_1\times\cU_{k+1}$,
 \[
   \widetilde h(u,z):=(u,z,\bbone_{A_1}(z),\ldots,\bbone_{A_k}(z)), \qquad (u,z)\in \cU_1\times\cU_1.
 \]
Let $M^{(1)}$ be the random point measure $M\circ \widetilde h^{-1}$, which is concentrated on $\cU_1\times\cU_{k+1}$.
By the mapping theorem for Poisson point processes,
 as we have seen above in case of $N^{(1)}$,
 one can check that $M^{(1)}$ is a Poisson random measure on $\cU_1\times\cU_{k+1}$ with intensity measure
 $\dd u\,\nu^{(k)}(\dd \bz)$, where the Borel measure $\nu^{(k)}$ on $\cU_{k+1}$ is given by
 \eqref{def_nu1}.

Since $(W_t)_{t\in\RR_+}$, $(W^{(1)}_t)_{t\in\RR_+}, \dots, (W^{(k)}_t)_{t\in\RR_+}$, $N$ and $M$ are mutually independent
 and $h$ and $\widetilde h$ are Borel measurable functions,
 we have that $(W_t)_{t\in\RR_+}$, $(W^{(1)}_t)_{t\in\RR_+}, \dots, (W^{(k)}_t)_{t\in\RR_+}$, $N^{(1)}=N\circ h^{-1}$, $N^{(2)}=\cdots=N^{(k+1)}=0$,
 and $M^{(1)}=M\circ \widetilde h^{-1}$ are mutually independent as well.

{\sl Step 3.}
  Using \eqref{help_SDE_joint}, Step 2 and Applebaum \cite[part (2) of Theorem 4.3.4]{App}, we have
 \begin{align}\label{SDE_X_J}
 \begin{split}
\begin{bmatrix}
    X_t \\
    J_t(A_1) \\
\vdots\\
J_t(A_k)\\
  \end{bmatrix}
  & = \begin{bmatrix}
        X_0 \\
       0 \\
       \vdots\\
        0\\
     \end{bmatrix}
     + \int_0^t \left(
                      \begin{bmatrix}
                          \beta \\
                          0 \\
                          \vdots\\
                          0\\
                       \end{bmatrix}
                     + \begin{bmatrix}
                         \tB & 0& \cdots& 0 \\
                         \mu(A_1) & 0& \cdots& 0 \\
                         \vdots & \vdots& \ddots& \vdots \\
                         \mu(A_k) & 0& \cdots& 0 \\
                       \end{bmatrix}
                       \begin{bmatrix}
                         X_u \\
                         J_u(A_1) \\
                        \vdots\\
                        J_u(A_k)\\
                       \end{bmatrix}
                      \right) \, \dd u \\
  &\quad + \int_0^t \sqrt{2 c X_u^+ }  \, \dd W_u \cdot \be_1^{(k+1)}
     + \underset{i=1}{\overset{k}{\sum}} \int_0^t 0  \, \dd W^{(i)}_u \cdot \be_{i+1}^{(k+1)}  \\
  &\quad
     + \int_0^t \int_{\cU_{k+1}} \int_{\cU_1}
            \bz \bbone_{\{w\leq X_{u-}\}} \, \widetilde{N^{(1)}}   (\dd u, \dd \bz, \dd w)  \\
   &\quad
      + \underset{i=1}{\overset{k}{\sum}} \int_0^t \int_{\cU_{k+1}} \int_{\cU_1}
         \bz
           \bbone_{\{w\leq J_{u-}(A_i)\}} \, \widetilde{N^{(i+1)}}  (\dd u, \dd \bz, \dd w) \\
    &\quad   + \int_0^t \int_{\cU_{k+1}}
      \br \, M^{(1)}(\dd u, \dd \br),
      \qquad t\in\RR_+.
 \end{split}
 \end{align}

Note that for any Borel measurable function \ $f:\cU_{k+1}\to\RR_+$ \ being integrable with respect to \ $\nu^{(k)}$,
 \ we have that
 \begin{align}\label{help_int1}
  \begin{split}
  \int_{\cU_{k+1}} f(r_1,\dots,r_{k+1})\,\nu^{(k)}(\dd r_1, \dots, \dd r_{k+1})
      &= \int_{\cU_1} f(r,\bbone_{A_1}(r),\ldots,\bbone_{A_k}(r))\,\nu(\dd r)\\ 
      &=  \sum_{\bell\in\{-1,1\}^k}  \int_{\bA^{\bell}} f(r, \tfrac{1}{2}(\bell+\bone) )\,\nu(\dd r),
 \end{split}       
 \end{align}
 as a consequence of the construction of the Lebesgue integral.
Similarly, for any Borel measurable function \ $f:\cU_{k+1}\to\RR_+$ \ being integrable with respect to $\mu_1^{(k)}$, we have that
 \begin{align}\label{help_int2}
  \begin{split}
  \int_{\cU_{k+1}} f(z_1, \dots, z_{k+1})\,\mu_1^{(k)}(\dd z_1, \dots, \dd z_{k+1})
       &= \int_{\cU_1} f(r,\bbone_{A_1}(r),\ldots,\bbone_{A_k}(r))\,\mu(\dd r)\\
       &= \sum_{\bell\in\{-1,1\}^k}  \int_{\bA^{\bell}} f(z, \tfrac{1}{2}(\bell+\bone))\,\mu(\dd z).
  \end{split}       
 \end{align}

Now, we show that the matrix that appears
 in the drift part of the SDE \eqref{SDE_X_J} is $\widetilde{\bB^{(k)}}$,
 where $\widetilde{\bB^{(k)}}$ is the matrix appearing in the corresponding SDE representation \eqref{SDE_atirasa_dimd} of
 $(X_t,J_t(A_1),\ldots,J_t(A_k))_{t\in\RR_+}$.
Since $\mu_n^{(k)}=0$ for any $n\in\{2,\dots,k+1\}$, by using \eqref{help5}, \eqref{help_int2},
 the facts that $\{ \bA^{\bell} : \bell \in \{-1,1\}^k\}$ is a partition of $\cU_1$, and
 $\{ \bA^{\bell} : \bell \in \{-1,1\}^k,\, \ell_{i-1}=1\}$ is a partition of $A_{i-1}$ for each $i\in\{2,\ldots,k+1\}$,
 we obtain that for any $i,j\in\{1,\dots, k+1\}$, it holds that
\begin{align*}
\widetilde{b^{(k)}}_{i,j}
  &= b^{(k)}_{i,j} + \int_{\cU_{k+1}}(z_i - \delta_{i,j})^+\,\mu_j^{(k)}(\dd \bz)\\
  & = B \bbone_{\{i=j=1\}}
    + \bbone_{\{i=j=1\}}
      \int_{\cU_1} (z-1)^+ \,\mu(\dd z)
    + \bbone_{\{j=1\}}
      \bbone_{\{i\neq 1\}} \int_{\cU_1} \bbone_{A_{i-1}}(z) \,\mu(\dd z) \\
 &= B \bbone_{\{i=j=1\}}
    + \bbone_{\{i=j=1\}}
      \int_{\cU_1} (z-1)^+ \,\mu(\dd z)
    + \bbone_{\{j=1\}}
      \bbone_{\{i\neq 1\}} \mu(A_{i-1})\\      
  &  = \bbone_{\{i=j=1\}} \tB + \bbone_{\{j=1\}}
      \bbone_{\{i\neq 1\}} \mu(A_{i-1}),
\end{align*}
 which shows that $\widetilde{\bB^{(k)}}$ is the matrix that appears
 in the drift part of the SDE \eqref{SDE_X_J}.

{\sl Step 4.} 
We check that the moment conditions in part (v) of Definition \ref{Def_admissible} and \eqref{moment_condition_m_new} 
 hold for $\nu^{(k)}$, and the moment conditions in part (vi) of Definition \ref{Def_admissible} 
 hold for $\mu^{(k)}_i$, $i\in\{1,\dots,k+1\}$.
By \eqref{help_int1}, we have that
\[
 \int_{\cU_{k+1}} (1 \land \|\br\|) \, \nu^{(k)}(\dd \br)
  \leq  \int_{\cU_{k+1}} \|\br\| \, \nu^{(k)}(\dd \br) 
  = \sum_{\bell\in\{-1,1\}^k}  \int_{\bA^{\bell}}  \|(r,\tfrac{1}{2}\bell+\tfrac{1}{2}\bone)\| \,\nu(\dd r).
\]
We will decompose the previous sum according to the cases  $\bell = -\bone$ or $\bell \neq -\bone$.
In the first case, i.e., when  $\bell = -\bone$, we have that $\bA^{-\bone}=\cU_1\setminus (\cup_{i=1}^k A_i)$,
 $\tfrac{1}{2} (-\bone + \bone) = \bzero$ and $\|(r, \bzero )\|=r$ for all $r\in\cU_1$.
In the second case, i.e., when $\bell \neq -\bone$, we have that
\begin{align}\label{help_bounded}
1<\|(r,\tfrac{1}{2}\bell+\tfrac{1}{2}\bone)\|
 =\sqrt{r^2+\tfrac{1}{4}(\ell_1+1)^2+\dots+\tfrac{1}{4}(\ell_k+1)^2}
 \leq \sqrt{r^2+k}\leq r+ \sqrt{k}
\end{align}
 for all \ $r\in\cU_1$.
Then, using \eqref{moment_condition_m_new} with $d=1$, the fact that
 $\{ \bA^{\bell} : \bell \in \{-1,1\}^k\}$ is a partition of $\cU_1$,
 and the moment condition for $\nu$ given in part (v) of Definition \ref{Def_admissible} with $d=1$,
 we have that
 \begin{align*}
   \int_{\cU_{k+1}} (1 \land \|\br\|) \, \nu^{(k)}(\dd \br)
    &\leq \int_{\cU_{k+1}} \|\br\| \, \nu^{(k)}(\dd \br) \\
   &\leq \int_{\bA^{-\bone}}  r\, \nu(\dd r)+
     \sum_{\{ \bell\in\{-1,1\}^k \,:\,\bell\neq -\bone\}}  \int_{\bA^{\bell}} (r+\sqrt{k}) \,\nu(\dd r)\\
   &= \int_{\cU_1} r\,\nu(\dd r)+
                        \sqrt{k}\sum_{\{\bell\in\{-1,1\}^k \,:\, \bell\neq -\bone\}}  \nu(\bA^{\bell}) 
       <\infty,
  \end{align*}
 since for any $\bell\in\{-1,1\}^k$ such that  $\bell\neq -\bone$,
 there exists an $i\in\left\{1,\dots, k\right\}$ that satisfies $\ell_i=1$ and  $\bA^{\bell}\subset A_{\ell_i}$
 yielding that $\nu(\bA^{\bell})\leq \nu(A_{\ell_i})<\infty$.
It means that $\nu^{(k)}$ satisfies the moment condition in part (v) of Definition \ref{Def_admissible}.
In particular, the previous calculation also shows that
 the moment condition \eqref{moment_condition_m_new} holds for the measure \ $\nu^{(k)}$.

Furthermore, by \eqref{help_int2}, we get that
 \begin{align*}
    &\int_{\cU_{k+1}}
       \biggl[\|\bz\| \land \|\bz\|^2 +\sum_{i=2}^{k+1} (1 \land z_i)\biggr] \mu_1^{(k)}(\dd \bz) \\
   &\qquad  =  \sum_{\bell\in\{-1,1\}^k}
    \int_{\bA^{\bell}}\biggl[
\|(z,\tfrac{1}{2}\bell+\tfrac{1}{2}\bone)\| \land \|(z,\tfrac{1}{2}\bell+\tfrac{1}{2}\bone)\|^2 + \sum_{i=1}^{k} (1 \land \tfrac{1}{2}(\ell_i+1))\biggr] \mu(\dd z).
 \end{align*}
Again, we decompose the previous sum according to the cases
 $\bell=-\bone$ or $\bell\neq -\bone$.
If $\bell = -\bone$, then $\|(z,\tfrac{1}{2}\bell+\tfrac{1}{2}\bone)\| = z$, $z\in\cU_1$, and $\tfrac{1}{2}(\ell_i+1)=0$, $i\in\{1,\ldots,k\}$.
If $\bell\neq -\bone$, then, by \eqref{help_bounded}, we have that $1<\|(z,\tfrac{1}{2}\bell+\tfrac{1}{2}\bone)\|\leq z+\sqrt{k}$, $z\in\cU_1$.
Hence using also the moment condition for $\mu$ given in part (vi) of Definition \ref{Def_admissible} with $d=1$ and the facts that
 $1\land \tfrac{1}{2}(\ell_i+1)\leq 1$, $i\in\{1,\ldots,k\}$, and $z\bbone_{\{z\geq 1\}}\leq z\wedge z^2$, $z\in\cU_1$, we get that
 \begin{align*}
    &\int_{\cU_{k+1}}
       \biggl[\|\bz\| \land \|\bz\|^2 +\sum_{i=2}^{k+1} (1 \land z_i)\biggr] \mu_1^{(k)}(\dd \bz) \\
   &\qquad \leq
       \int_{\bA^{-\bone}}
      (z\land z^2)\mu(\dd z)
    + \sum_{ \{ \bell\in\{-1,1\}^k \,:\, \bell\neq -\bone\} }  \int_{\bA^{\bell}} \biggl[ z+ \sqrt{k} + k \biggr]\, \mu(\dd z) \\
  &\qquad {= \int_{\bA^{-\bone}} 
       (z\land z^2)\mu(\dd z) }\\
  &\phantom{\qquad =\,}     
        + \sum_{ \{ \bell\in\{-1,1\}^k \,:\, \bell\neq -\bone\} } \left(  \int_{\bA^{\bell}}z\bbone_{\{z\geq 1\}} \,\mu(\dd z)
            +\int_{ \bA^{\bell}} z\bbone_{\{z<1\}}\, \mu(\dd z) + (\sqrt{k}+k)\mu(\bA^\bell)\right)\\
   &\qquad  \leq
       \int_{\cU_1}
       (z\land z^2)\mu(\dd z)
        + \sum_{ \{ \bell\in\{-1,1\}^k \,:\, \bell\neq -\bone\} } \left(  \int_{\cU_1}z\bbone_{\{z\geq 1\}}\, \mu(\dd z)
            +\int_{ \bA^{\bell}}\bbone_{\{z< 1\} } \,\mu(\dd z) +2k\mu(\bA^\bell)\right)\\
  &\qquad \leq \int_{\cU_1} (z\land z^2)\mu(\dd z)
               + \sum_{ \{ \bell\in\{-1,1\}^k \,:\, \bell\neq -\bone\} } \left( \int_{\cU_1} (z\wedge z^2)\,\mu(\dd z)  + \mu(\bA^\bell) + 2k\mu(\bA^\bell) \right) \\
  & \qquad = 2^k  \int_{\cU_1}
      (z\land z^2)\mu(\dd z) + (1+2k) \sum_{ \{ \bell\in\{-1,1\}^k \,:\, \bell\neq -\bone\} } \mu(\bA^\bell)<\infty,
 \end{align*}
 since for any $\bell\in\{-1,1\}^k$ such that  $\bell\neq -\bone$, there exists
 an $i\in\left\{1,\dots, k\right\}$ that satisfies $\ell_i=1$ and  $\bA^{\bell}\subset A_{\ell_i}$ yielding that $\mu(\bA^{\bell})\leq \mu(A_{\ell_i})<\infty$.
 Further, since $\mu_{i}^{(k)}=0$ for any $i=2,\dots, k+1$, it trivially holds that
 \[
 \int_{\cU_{k+1}}
       \biggl[\|\bz\| \land \|\bz \|^2 + \sum_{j \in \{1, \ldots, k+1\} \setminus \{i\}}  (1 \land z_j)\biggr] \mu_i^{(k)}(\dd \bz)
         =0, \qquad i=2,\dots, k+1.
 \]
It means that $\mu^{(k)}_i$, $i\in\{1,\dots,k+1\}$, satisfy the moment conditions in part (vi) of Definition \ref{Def_admissible}.

{\sl Step 5.} Using Steps 1--4 and Theorem 4.6 in Barczy et al.\ \cite{BarLiPap2}, we get that the stochastic process
 \ $(X_t,J_t(A_1),\dots, J_t(A_k))_{t\in\RR_+}$ \ is a CBI process
 with parameters \ $(k+1, \bc^{(k)}, \Bbeta^{(k)}, \bB^{(k)}, \nu^{(k)}, \bmu^{(k)})$
 starting from $(X_0,\bzero)$.
\proofend

In Theorem \ref{Thm_(k+1)CBI}, we proved that 
$(X_t,J_t(A_1),\ldots,J_t(A_k))_{t\in\RR_+}$ is a $(k+1)$-type CBI process starting from $(X_0,\bzero)$,
 by deriving an SDE for $(X_t,J_t(A_1),\ldots,J_t(A_k))_{t\in\RR_+}$ 
and then using a general representation of multi-type CBI processes
 via stochastic equations due to Barczy et al.\ \cite[Theorem 4.6]{BarLiPap2}.
We also determined the parameters of $(X_t,J_t(A_1),\ldots,J_t(A_k))_{t\in\RR_+}$.
In the next Corollary \ref{Cor_Laplace}, having the parameters of the multi-type
 CBI process in question at hand, using Theorem \ref{CBI_exists}, we explicitly write out the branching and immigration mechanisms,
 and derive a formula for the Laplace transform of $(X_t,J_t(A_1),\ldots,J_t(A_k))$ 
 in terms of a solution of a deterministic differential equation, where $t\in\RR_{++}$.
This corollary is a generalization of Corollary 10.12 in Li \cite{Li}, namely,
 in the special case $k=1$, our Corollary \ref{Cor_Laplace} gives back Corollary 10.12 in Li \cite{Li}.
For completeness, we mention that, in case of $k=1$, in the formulae $\Phi_A$ and $\Psi_A$ for the branching and immigration mechanisms 
 of $(X_t,J_t(A_1))_{t\in\RR_+}$ in Li \cite[page 293]{Li}
 (using the notations of Li \cite{Li}), one should write $-z\lambda_1 - \lambda_2 \bone_A(z)$
 instead of $-z(\lambda_1+\lambda_2\bone_A(z))$.

\begin{Cor}\label{Cor_Laplace}
Let $(X_t)_{t\in\RR_+}$ be a CBI process with parameters $(1, c, \beta, B, \nu, \mu)$ such that
 $\EE(X_0) < \infty$ and the moment condition \eqref{moment_condition_m_new} with $d=1$ hold.
Let $k\in\NN$ and $A_1,\ldots,A_k\in\cB((0,\infty))$ be such that $\nu(A_i) + \mu(A_i)<\infty$, $i\in\{1,\ldots,k\}$.
Then
\begin{itemize}
\item[(i)] the immigration mechanism of \ $(X_t,J_t(A_1),\ldots,J_t(A_k))_{t\in\RR_+}$ is the function $ \psi^{(k)}:\RR_+^{k+1}\rightarrow \RR_+$ given by 
    \begin{align}\label{2CBI_immigration}
 \psi^{(k)}(\lambda_0,\lambda_1,\ldots,\lambda_k):= \beta\lambda_0 +  \int_{\cU_1}\bigl( 1 - \ee^{-\lambda_0 r - \sum_{i=1}^k \lambda_i \bbone_{A_i}(r)} \bigr) \, \nu(\dd r)
   \end{align}
 for all $(\lambda_0,\lambda_1,\ldots,\lambda_k)\in\RR_+^{k+1}$,
\item[(ii)] the branching mechanism of $(X_t,J_t(A_1),\ldots,J_t(A_k))_{t\in\RR_+}$ is the function 
 \[ 
   \bvarphi^{(k)}=(\varphi^{(k)}_1, \ldots, \varphi^{(k)}_{k+1}): \RR_+^{k+1}\rightarrow \RR^{k+1}
 \]  
  given by 
\begin{align}\label{2CBI_branching}
 \begin{split}
&\varphi_1^{(k)}(\lambda_0,\lambda_1,\ldots,\lambda_k)
:= c\lambda_0^2 - B\lambda_0\\
&\phantom{\varphi_1^{(k)}(\lambda_0,\lambda_1,\ldots,\lambda_k):=}
+ \int_{\cU_1}\bigl( \ee^{- \lambda_0 z - \sum_{i=1}^k \lambda_i \bbone_{A_i}(z)} - 1+ \lambda_0 (1 \land z) \bigr) \, \mu(\dd z),\\
&\varphi_i^{(k)}(\lambda_0,\lambda_1,\ldots,\lambda_k):=0,\qquad i\in\{2,\ldots,k+1\}
\end{split}
\end{align}
 for all $(\lambda_0,\lambda_1,\ldots,\lambda_k)\in\RR_+^{k+1}$,

\item[(iii)] for all  $(\lambda_0,\lambda_1,\ldots,\lambda_k)\in\RR_+^{k+1}$, the unique locally bounded solution
 \[
  \RR_+\ni t\mapsto \bv^{(k)}(t,\lambda_0,\lambda_1,\dots,\lambda_{k})\in\RR_+^{k+1} 
 \]
 to the system \eqref{help12} associated with the branching mechanism \eqref{2CBI_branching} takes the form
\begin{align}\label{eq:vJumps}
  \bv^{(k)}(t,\lambda_0,\lambda_1,\dots,\lambda_{k}) =v^{(k)}(t,\lambda_0,\lambda_1,\dots,\lambda_{k})\be^{(k+1)}_1+\sum_{i=1}^{k}\lambda_i\be^{(k+1)}_{i+1}, 
    \qquad {t\in\RR_+,}
\end{align}
where $\RR_+\ni t\mapsto v^{(k)}(t,\lambda_0,\lambda_1,\dots,\lambda_{k})$ 
is the unique locally bounded non-negative solution to the differential equation
\begin{align}\label{def_v_k+1_CBI}
\begin{split}
&\partial_1 v^{(k)}(t,\lambda_0,\lambda_1,\ldots,\lambda_k)
=  - \varphi_1^{(k)}\big( v^{(k)}(t,\lambda_0,\lambda_1,\ldots,\lambda_k),\lambda_1,\ldots,\lambda_k\big),\\
&v^{(k)}(0, \lambda_0,\lambda_1,\ldots,\lambda_k) = \lambda_0,
\end{split}
\end{align}

\item[(iv)] for all $t\in\RR_+$, the Laplace transform of $(X_t,J_t(A_1),\ldots,J_t(A_k))$  is 
 \begin{align}\label{2CBI_Laplace}
  \begin{split}
    &\EE\left( \ee^{-\lambda_0 X_t - \sum_{i=1}^k \lambda_i J_t(A_i)}  \mid X_0 = x\right) \\
    &= \exp\Bigg\{ -x v^{(k)}(t,\lambda_0,\lambda_1,\ldots,\lambda_k) 
                 - \int_0^t \psi^{(k)}\big( v^{(k)}(s,\lambda_0,\lambda_1,\ldots,\lambda_k), \lambda_1,\ldots,\lambda_k) \,\dd s\Bigg\}
  \end{split}
 \end{align}
 for all $(\lambda_0,\lambda_1,\ldots,\lambda_k)\in\RR_+^{k+1}$ and $x\in\RR_+$.
\end{itemize}
\end{Cor}

\noindent{\bf Proof.}
First, recall that in Theorem \ref{Thm_(k+1)CBI} we proved that 
  $(X_t,J_t(A_1),\ldots,J_t(A_k))_{t\in\RR_+}$ is a multi-type CBI process with parameters
  \ $(k+1, \bc^{(k)}, \Bbeta^{(k)}, \bB^{(k)}, \nu^{(k)}, \bmu^{(k)})$ starting from $(X_0,\bzero)$.
Having the set of admissible parameters at hand, we will show the assertion using Theorem \ref{CBI_exists}.

(i).
We determine the immigration mechanism of the multi-type CBI process in question.
For all $(\lambda_0,\lambda_1,\ldots,\lambda_k)\in\RR_+^{k+1}$,
 using \eqref{help_int1} with
 \[
  f(r_0,r_1,\ldots,r_k):= 1 - \ee^{- \big\langle (\lambda_j)_{j=0}^k, (r_j)_{j=0}^k \big\rangle}
                        = 1 - \ee^{-\lambda_0r_0 - \sum_{i=1}^k \lambda_i r_i},
   \qquad  (r_0,r_1,\ldots,r_k)\in\cU_{k+1},
 \]
 we get that
 \begin{align*}
    \int_{\cU_{k+1}}
         \bigl( 1 - \ee^{-\lambda_0r_0 - \sum_{i=1}^k \lambda_i r_i} \bigr)
         \, \nu^{(k)}(\dd r_0,\dd r_1,\ldots,\dd r_k) 
    = \int_{\cU_1}  \left( 1 -\ee^{-\lambda_0 r - \sum_{i=1}^k \lambda_i \bbone_{A_i}(r) }  \right) \, \nu(\dd r).
 \end{align*}
Consequently, by Theorems \ref{CBI_exists} and \ref{Thm_(k+1)CBI}, we get that
 the immigration mechanism of the $(k+1)$-type CBI process \ $(X_t,J_t(A_1),\ldots,J_t(A_k))_{t\in\RR_+}$ 
 starting from $(X_0,\bzero)$ is given by \eqref{2CBI_immigration}, as desired.

(ii).
 We are going to determine the branching mechanism of the multi-type CBI process in question.
For all \ $(\lambda_0,\lambda_1,\ldots,\lambda_k)\in\RR_+^{k+1}$, we have
 \begin{align*}
  \left\langle \bB^{(k)} \be_i^{(k+1)}, (\lambda_j)_{j=0}^k \right\rangle
    & = \left\langle
               (B\bbone_{\{\ell=j=1\}})_{\ell,j=1}^{k+1}
                  (\delta_{j,i})_{j=1}^{k+1},
              (\lambda_j)_{j=0}^k
        \right\rangle 
       =  \left\langle (B\bbone_{\{i=\ell=1\}})_{\ell=1}^{k+1},  (\lambda_j)_{j=0}^k  \right\rangle    \\
    & = \begin{cases}
         \big\langle
              (B \bbone_{\{\ell=1\}})_{\ell=1}^{k+1},
              (\lambda_j)_{j=0}^k
        \big\rangle
         = B\lambda_0  & \text{if \ $i=1$,}\\[1mm]
       \left\langle
              \bzero,
              (\lambda_j)_{j=0}^k
        \right\rangle
        = 0  & \text{if \ $i\in\{2,\ldots,k+1\}$.}
       \end{cases} 
 \end{align*}
Furthermore, for all \ $(\lambda_0,\lambda_1,\ldots,\lambda_k)\in\RR_+^{k+1}$, using \eqref{help_int2} with
 \[
   f(z_0,z_1,\ldots,z_k):= \ee^{- \lambda_0 z_0 - \sum_{i=1}^k \lambda_i z_i} - 1 + \lambda_0 (1 \land z_0),
             \qquad (z_0,z_1,\ldots,z_k)\in\cU_{k+1},
  \]
 we get that
 \begin{align*}
  &\int_{\cU_{k+1}}
         \bigl( \ee^{-\lambda_0 z_0 - \sum_{i=1}^k \lambda_i z_i} - 1 + \lambda_0 (1 \land z_0) \bigr)
         \, \mu_1^{(k)}(\dd z_0,\dd z_1,\ldots,\dd z_k)\\
 &\qquad = \int_{\cU_1}\left[ \ee^{-\lambda_0 z - \sum_{i=1}^k \lambda_i \bbone_{A_i}(z) }
  - 1 + \lambda_0 (1 \land z) \right] \,\mu(\dd z).
 \end{align*}
Since $\mu^{(k)}_j=0$, $j\in\{2,\ldots,k+1\}$, it trivially holds that
 \[
  \int_{\cU_{k+1}}
         \bigl( \ee^{-\lambda_0 z_0 - \sum_{i=1}^k \lambda_i z_i} - 1 + \lambda_j (1 \land z_j) \bigr)
         \, \mu_j^{(k)}(\dd z_0,\dd z_1,\ldots,\dd z_k) = 0,\qquad j\in\{2,\ldots,k+1\}.
 \]
Consequently, by Theorems \ref{CBI_exists} and \ref{Thm_(k+1)CBI}, we get that
 the branching mechanism of the $(k+1)$-type CBI process \ $(X_t,J_t(A_1),\ldots,J_t(A_k))_{t\in\RR_+}$
 starting from $(X_0,\bzero)$ is given by \eqref{2CBI_branching}, as desired.

(iii).
For all $(\lambda_0,\lambda_1,\ldots,\lambda_k)\in\RR_+^{k+1}$, let
 \begin{align*}
    \RR_+\ni t \mapsto 
      \bv^{(k)}(t,\lambda_0,\lambda_1,\ldots,\lambda_k)=
         \big(v^{(k)}_0(t,\lambda_0,\lambda_1,\ldots,\lambda_k),\ldots,v^{(k)}_k(t,\lambda_0,\lambda_1,\ldots,\lambda_k)\big)\in\RR_+^{k+1}
 \end{align*}
 denote the unique locally bounded solution to the system \eqref{help12} associated with the branching mechanism \eqref{2CBI_branching}, i.e.,
\begin{align*}
  \begin{cases}
      \partial_1 v_i^{(k)}(t,\lambda_0,\lambda_1,\ldots,\lambda_k) 
         = -\varphi_{i+1}^{(k)}\big( v_0^{(k)}(t,\lambda_0,\lambda_1,\ldots,\lambda_k), \ldots, v_k^{(k)}(t,\lambda_0,\lambda_1,\ldots,\lambda_k)\big), & \text{}\\
      v_i^{(k)}(0, \lambda_0,\lambda_1,\ldots,\lambda_k) = \lambda_i, \qquad i\in\{0,1,\ldots,k\}, & \text{}
   \end{cases}   
 \end{align*} 
 where the functions $\varphi^{(k)}_i$, $i\in\{1,\ldots,k+1\}$, are given by \eqref{2CBI_branching}. Since $\varphi^{(k)}_{i+1}=0$, $i\in\{1,\ldots,k\}$,
 we have $v_i^{(k)}(t,\lambda_0,\lambda_1,\ldots,\lambda_k)=\lambda_i$, $t\in\RR_+$, $i\in\{1,\ldots,k\}$.  
Consequently, with $v^{(k)}:=v_0^{(k)}$, we get that 
 $\RR_+\ni t\mapsto v^{(k)}(t,\lambda_0,\lambda_1,\ldots,\lambda_k)$ is the unique locally bounded non-negative solution 
 to the differential equation \eqref{def_v_k+1_CBI}.
Moreover, 
 \[
 \bv^{(k)}(t,\lambda_0,\lambda_1,\ldots,\lambda_k)
    = \big(v^{(k)}(t,\lambda_0,\lambda_1,\ldots,\lambda_k),\lambda_1,\ldots,\lambda_k\big),\qquad t\in\RR_+, 
 \] 
 which implies \eqref{eq:vJumps}.

(iv).
Finally, we determine the joint Laplace transform of the multi-type CBI process in question.
Using Theorems \ref{CBI_exists} and \ref{Thm_(k+1)CBI}, equation \eqref{eq:vJumps} and that $J_0(A_i)=0$, $i\in\{1,\ldots,k\}$, we have that
\begin{align*}
 \begin{split}
&\EE\left( \ee^{-\lambda_0 X_t - \sum_{i=1}^k \lambda_i J_t(A_i)}  \mid X_0 = x\right) \\
&\qquad = \exp\Bigg\{ - \Big\langle 
(x,\bzero), \bv^{(k)}(t, \lambda_0,\lambda_1,\dots,\lambda_k) \Big\rangle - \int_0^t \psi^{(k)}\big( \bv^{(k)}(s, \lambda_0,\lambda_1,\dots,\lambda_k) \big) \,\dd s\Bigg\}\\
   &\qquad = \exp\Bigg\{ -  x v^{(k)}(t,\lambda_0,\lambda_1,\ldots,\lambda_k)
 - \int_0^t \psi^{(k)}\big(v^{(k)}(s,\lambda_0,\lambda_1,\ldots,\lambda_k), \lambda_1,\ldots,\lambda_k\big) \,\dd s\Bigg\}                                          
 \end{split}
 \end{align*}
 for all $(\lambda_0,\lambda_1,\ldots,\lambda_k)\in\RR_+^{k+1}$, $x\in\RR_+$, and $t\in\RR_+$. 
Hence \eqref{2CBI_Laplace} holds.
\proofend

\smallskip 

For all $k\in\NN$, $A_1,\ldots,A_k\in\cB((0,\infty))$ and $\blambda=(\lambda_1,\dots,\lambda_k)\in\RR_+^k$, let us introduce the notation
 \begin{equation}\label{Alambda}
  A_\blambda:=\bigcup_{i\in\{1,\dots,k\}\,:\,\lambda_i\neq 0}A_i.
 \end{equation}

In the next lemma, we provide some useful properties of the function $v^{(k)}$ defined in \eqref{def_v_k+1_CBI}, 
 that will be applied in Section \ref{Section_CBI_results_2}.

\begin{Lem}\label{Lem_vjumps}
Let $(1, c,0, B, 0,\mu)$ be a set of admissible parameters.
Let $k\in\NN$ and $A_1,\ldots,A_k\in\cB((0,\infty))$ be such that $\mu(A_i)<\infty$, $i\in\{1,\ldots,k\}$. 
For all $(\lambda_0,\lambda_1,\dots,\lambda_k)\in\RR_+^{k+1}$, let us consider the function
 \[
  \RR_+ \ni t\mapsto  v^{(k)}(t,\lambda_0, \lambda_1,\dots,\lambda_k)\in\RR_+,
 \] 
  which is the unique locally bounded solution to the differential equation \eqref{def_v_k+1_CBI}. 
Then the following assertions hold:
 \begin{enumerate}
\item[(i)] For all $i\in\{0,1,\dots,k\}$ and  $t,\lambda_j \in\RR_+$, $j\in \{0,1,\dots,k\}\setminus\{i\}$, the function
 \[
  \RR_+\ni\lambda_i \mapsto v^{(k)}(t,\lambda_0,\lambda_1,\dots,\lambda_k)
 \]
 is monotone increasing. 
\item[(ii)] 
Let $t,\lambda_0\in\RR_+$ and $\blambda=(\lambda_1,\dots,\lambda_k)\in\RR_+^k$. 
Suppose that $\RR_+\ni s\mapsto \eta_s\in\RR_+$ is a monotone increasing function such that $\eta_s\uparrow\lambda_0$ as $s\uparrow \infty$.
Then we have that
 \[
   v^{(k)}(t,\eta_s,s\lambda_1,\dots,s\lambda_k) \uparrow  v^{(A_\blambda)}(t,\lambda_0) \qquad \text{as $s\uparrow\infty$,} 
 \] 
 where $A_\blambda$ is given by \eqref{Alambda} and the function  $\RR_+\ni t \mapsto v^{(A_\blambda)}(t,\lambda_0){\in\RR_+}$ is the unique locally bounded non-negative solution
 to the differential equation
  \begin{align}\label{eq_vAlambda}
  \partial_1 v^{(A_\blambda)}(t,\lambda_0) =\mu(A_\blambda)- \varphi^{(A_\blambda)}(v^{(A_\blambda)}(t,\lambda_0)),\quad t\in\RR_+,\qquad  v^{(A_\blambda)}(0,\lambda_0) = \lambda_0,
  \end{align}
  where the function $\varphi^{(A_\blambda)}:\RR_+\to\RR$ is given by \eqref{eq:psiA} with $A:=A_\blambda$. 
\end{enumerate}
\end{Lem}

We call the attention that the limit $ v^{(A_\blambda)}(t,\lambda_0)$ in part (ii) of Lemma \ref{Lem_vjumps} does not depend
 on the choice of the function $\RR_+\ni s\mapsto \eta_s$ (only on its limit $\lambda_0$ as $s\uparrow \infty$).
Further, if $\lambda_0=0$ in part (ii) of Lemma \ref{Lem_vjumps}, then $\eta_s=0$, $s\in\RR_+$.
The proof of Lemma \ref{Lem_vjumps} can be found in Appendix \ref{Appendix_Lemma_v}.

In the next remark, we specialize part (ii) of Lemma \ref{Lem_vjumps} to the case $\lambda_i=0$, $i\in\{1,\ldots,k\}$
 by giving a direct proof (in the considered special case) as well.

\begin{Rem}
Let $(1, c,0, B, 0,\mu)$ be a set of admissible parameters.
Let $k\in\NN$ and $A_1,\ldots,A_k\in\cB((0,\infty))$ be such that $\mu(A_i)<\infty$, $i\in\{1,\ldots,k\}$. 
Let $\lambda_i:=0$, $i\in\{1,\ldots,k\}$, i.e., $\blambda:=\bzero\in\RR^k$.
Then 
 \begin{align*}
  \varphi_1^{(k)}(\lambda_0,\lambda_1,\ldots,\lambda_k)
    &= \varphi_1^{(k)}(\lambda_0,0,\ldots,0)\\
    &= c\lambda_0^2 - B\lambda_0
      + \int_{\cU_1}\bigl( \ee^{- \lambda_0 z } - 1+ \lambda_0 (1 \land z) \bigr) \, \mu(\dd z),
      \qquad \lambda_0\in\RR_+,
 \end{align*}
 where $\varphi_1^{(k)}$ is given by \eqref{2CBI_branching}, and hence the function $\RR_+\ni\lambda_0 \mapsto \varphi_1^{(k)}(\lambda_0,0,\ldots,0)$
 is nothing else but the branching mechanism of a CB process with parameters $(1, c,0, B, 0,\mu)$.
Therefore, for all $\lambda_0\in\RR_+$, the function 
 $\RR_+\ni t\mapsto v^{(k)}(t,\lambda_0,0,\ldots,0)$ given by \eqref{def_v_k+1_CBI} coincides with the function 
 $\RR_+\ni t\mapsto v_1(t,\lambda_0)$ given by \eqref{help12} with $d=1$.
Hence, using also part (i) of Lemma \ref{Lem_vjumps}, for all $t\in\RR_+$, the function 
 $\RR_+\ni \lambda_0 \mapsto v^{(k)}(t,\lambda_0,0,\ldots,0) = v_1(t,\lambda_0)$ is monotone increasing.
Recall also that, by Theorem \ref{CBI_exists}, the function $\RR_+\times\RR_+\ni(t,\lambda_0)\mapsto v_1(t,\lambda_0)$ is continuous.
Therefore, if $\lambda_0\in\RR_+$ and $\RR_+\ni s\mapsto \eta_s \in\RR_+$ is a function such that $\eta_s\uparrow \lambda_0$ as $s\uparrow \infty$, 
 then, for all $t\in\RR_+$, we have that
 \[
   v^{(k)}(t,\eta_s,s\lambda_1,\ldots,s\lambda_k)
      = v^{(k)}(t,\eta_s,0,\ldots,0)
      =v_1(t,\eta_s) \uparrow v_1(t,\lambda_0)\qquad \text{as \ $s\uparrow \infty$.}
 \] 
This is in accordance with part (ii) of Lemma \ref{Lem_vjumps}, since, using that $\lambda_i=0$, $i\in\{1,\ldots,k\}$,
 we have that $A_\blambda = \emptyset$, $\mu(A_\blambda)=0$, and $\varphi^{(A_\blambda)}$ (given by \eqref{eq:psiA} with $A=A_\blambda$) coincides with the branching mechanism of a CB process 
 with parameters $(1, c,0, B, 0,\mu)$, and therefore, for all $\lambda_0\in\RR_+$,
 the function $\RR_+\ni t \mapsto v^{(A_\blambda)}(t,\lambda_0)$ given by \eqref{eq_vAlambda}
 coincides with the function $\RR_+\ni t \mapsto v_1(t,\lambda_0)$ given by \eqref{help12} with $d=1$.
\proofend
\end{Rem}

\section{Joint distribution function of a CBI process 
 and its first jump times with jump sizes in given Borel sets}\label{Section_CBI_results_2}
  
First, we recall a result on the distribution function of $\tau_A$ (introduced in \eqref{tau_def})
 for single-type CBI processes, where $A\in\cB((0,\infty))$ is such that
 its total L\'evy measure is finite, see Theorem \ref{Thm_single_jump_times_CBI}.
This result is due to He and Li \cite[Theorem 3.1]{HeLi},
 see also Theorem 10.13 in the new second edition of Li's book \cite{Li}.
For a corresponding result for multi-type CBI processes, 
see Barczy and Palau \cite[Theorem 4.2]{BarPal} (or Theorem \ref{Thm_jump_times_dCBI}).
The proof techniques of He and Li \cite[Theorem 3.1]{HeLi} and Barczy and Palau \cite[Theorem 4.2]{BarPal} 
 are the same, both use an SDE representation of (multi-type) CBI processes.
The proof technique of Li \cite[Theorem 10.13]{Li} is different,
 it is based on the Laplace transform of the $2$-type CBI process $(X_t,J_t(A))_{t\in\RR_+}$ 
 starting from $(X_0,0)$ (see Corollary 3.12 in Li \cite{Li} or our Corollary \ref{Cor_Laplace} with $k=1$). 

\begin{Thm}\label{Thm_single_jump_times_CBI}
Let $(X_t)_{t\in\RR_+}$ be a CBI process with parameters $(1, c, \beta, B, \nu, \mu)$ such that
 $\EE(X_0) < \infty$ and the moment condition \eqref{moment_condition_m_new} with $d=1$ hold.
Let $A\in\cB((0,\infty))$ be such that $\nu(A) + \mu(A)<\infty$.
Then we have
 \begin{align}\label{formula3_multiple}
  \begin{split}
  \PP_x(\tau_{A}>t) = \exp\left\{ -\nu(A) t - x \tv^{(A)}(t)
- \int_0^t \psi^{(A)}(\tv^{(A)}(s) )\,\dd s \right\}
 \end{split}
 \end{align}
 for all $t\in\RR_+$ and $x\in\RR_+$, where the continuously differentiable function
 $\RR_+\ni t \mapsto \tv^{(A)}(t) \in\RR_+$ is the unique locally bounded solution to the differential
 equation
 \begin{align}\label{DE_jump_time}
     (\tv^{(A)})'(t)  = \mu(A) - \varphi^{(A)}(\tv^{(A)}(t) ), \qquad t\in\RR_+,\qquad
     \tv^{(A)}(0)=0.
 \end{align}
Here the functions $\psi^{(A)}$ and $\varphi^{(A)}$ are given by \eqref{help13} and \eqref{eq:psiA}, respectively.
\end{Thm}

We call the attention that the differential equation \eqref{DE_jump_time} coincides with the differential equation  \eqref{help_DE_tv} with $d=1$, and for simplicity instead of 
 $\tv_1^{(A)}(t,\mu(A))$ we write  $\tv^{(A)}(t)$ for the solution.
Further, note that,  using \eqref{eq:psiA}, the differential equation \eqref{DE_jump_time} can also be written in the form 
 \[
  (\tv^{(A)})'(t) = -\varphi(\tv^{(A)}(t)) + \int_A\ee^{-\tv^{(A)}(t) z}\,\mu(\dd z), \qquad t\in\RR_+,\qquad \tv^{(A)}(0)=0.
 \]

Next, we present our main result, we derive an expression for the joint distribution function of $(\tau_{A_1},\ldots,\tau_{A_k})$ 
 at a point $(t_1,\ldots,t_k)\in\RR_+^k$ satisfying $0\leq t_1\leq t_2\leq\ldots\leq t_k$,
 in terms of the solution of a system of deterministic differential equations,
 where $k\in\NN$ and $A_1,\ldots,A_k\in\cB((0,\infty))$ having finite total L\'evy measures.
 
\begin{Thm}\label{Thm_main}
Let $(X_t)_{t\in\RR_+}$ be a CBI process with parameters $(1, c, \beta, B, \nu, \mu)$ such that
 $\EE(X_0) < \infty$ and the moment condition \eqref{moment_condition_m_new} with $d=1$ hold.
Let $k\in\NN$ and $A_1,\ldots,A_k\in\cB((0,\infty))$ be such that $\nu(A_i) + \mu(A_i)<\infty$, $i\in\{1,\ldots,k\}$.
Let $C_i:=\bigcup_{j=i}^k A_j$, $i\in\{1,\ldots,k\}$. 
Then, for all $0 =:t_0\leq t_1\leq t_2\leq\ldots\leq t_k$ and $x\in\RR_+$, we have that
 \begin{align*}%\label{formula_multiple_diff_times}
  %\begin{split}
  &\PP_x(\tau_{A_1}>t_1,\ldots,\tau_{A_k}>t_k) 
    = \exp\left\{    -x \tw^{(k)}_1(t_1)-\sum_{i=1}^{k}\int_{0}^{t_i-t_{i-1}}
        \widetilde{\psi}^{(i)}\big( \tw^{(k)}_i(u)\big)\,\dd u   \right\},
 %\end{split}
 \end{align*}
 where
 \begin{itemize}
  \item $\tw^{(k)}_k(t):=v^{(C_k)}(t,0)$, $t\in\RR_+$, and
       \begin{align}\label{help_new_01}
 \tw^{(k)}_i(t):=v^{(C_i)}(t, \tw^{(k)}_{i+1}(t_{i+1}-t_{i})), \qquad t\in\RR_+,
                \;\; i\in\{1,\ldots,k-1\},
       \end{align}
   where, for all $\lambda_0\in\RR_+$ and $i\in\{1,\ldots,k\}$, the continuously differentiable function $\RR_+\ni t \mapsto v^{(C_i)}(t,\lambda_0){\in\RR_+}$
    is the unique locally bounded non-negative solution to the differential equation 
  \begin{align}\label{help_DE_main_result}
  \partial_1 v^{(C_i)}(t,\lambda_0) =\mu(C_i)- \varphi^{(C_i)}(v^{(C_i)}(t,\lambda_0)),\quad t\in\RR_+,\qquad  v^{(C_i)}(0,\lambda_0) = \lambda_0,
  \end{align}
  where the function $\varphi^{(C_i)}:\RR_+\to\RR$ is given by
  $$
  \varphi^{(C_i)}(\lambda_0)
   = c\lambda_0^2 - \left(B - \int_{C_i} (1\wedge z)\,\mu(\dd z)\right)\lambda_0+ \int_{\cU_1\setminus C_i} \bigl( \ee^{- \lambda_0 z } - 1 + \lambda_0 (1 \land z) \bigr) \, \mu(\dd z)
  $$
  for all $\lambda_0\in\RR_+$,
 \item $\widetilde{\psi}^{(i)}(\lambda_0) :=\psi^{(C_i)}(\lambda_0) + \nu(C_i)$, $\lambda_0\in\RR_+$, $i\in\{1,\ldots,k\}$,
       where the function $\psi^{(C_i)}:\RR_+\to\RR_+$  is given by
       $$
       \psi^{(C_i)}(\lambda_0)
= \beta\lambda_0 + \int_{\cU_1\setminus C_i}\bigl( 1 - \ee^{-\lambda_0 r} \bigr) \, \nu(\dd r),
\qquad \lambda_0\in\RR_+.
       $$
 \end{itemize} 
\end{Thm}

Note that, for any $i\in\{1,\dots, k\}$, the function $\varphi^{(C_i)}$ is nothing else but the 
 function given in \eqref{eq:psiA} with $A=C_i$,
 and the function $\psi^{(C_i)}$ is nothing else but the function given in \eqref{help13} with $A=C_i$. 
Moreover, for any $i\in\{1,\dots, k\}$, the differential equations \eqref{help_DE_main_result} and \eqref{eq_vAlambda} 
 with $\blambda =\sum_{j=i}^{k}\be^{(k)}_{j}=: \biota_i$ coincide.
Using the notation \eqref{Alambda}, in fact, we have that $C_i=A_{\biota_i}$, $i\in\{1,\ldots,k\}$,
 and part (ii) of Lemma \ref{Lem_vjumps} yields that, for all $\lambda_0\in\RR_+$,
 there is a unique locally bounded non-negative solution to the differential equation \eqref{help_DE_main_result}.
We mention that, using \eqref{eq:psiA}, the differential equation \eqref{help_DE_main_result} can also be written in the form
  \[
   \partial_1 v^{(C_i)}(t,\lambda_0) = - \varphi(v^{(C_i)}(t,\lambda_0)) 
                                       + \int_{C_i} \ee^{- v^{(C_i)}(t,\lambda_0)\, z}\,\mu(\dd z), 
                                       \qquad t\in\RR_+,
  \]
 with initial condition $v^{(C_i)}(0,\lambda_0) = \lambda_0$.

\noindent{\bf Proof of Theorem \ref{Thm_main}.} 
 By part (i) of Theorem \ref{Thm_jump_times_dCBI} (in the proof of which only the equality \eqref{help4} was used,
 see Step 3 of Barczy and Palau \cite[Theorem 4.2]{BarPal}) and the fact that $J_0(A_i)=0$, $i\in\{1,\ldots,k\}$,
 we have $\EE(J_t(A_i))<\infty$, $t \in \RR_+$, $i\in\{1,\ldots,k\}$, which, in particular,
 implies that $\PP(J_t(A_i) < \infty) = 1$, $t\in\RR_+$, $i\in\{1,\ldots,k\}$.

Let $0=t_0\leq t_1\leq t_2\leq\ldots\leq t_k$ be arbitrarily fixed.
We divide the proof into several steps.

{\sl Step 1.}
We check that 
 \begin{align}\label{help202_versionk}
   \PP( \tau_{A_1}>t_1,\ldots,\tau_{A_k}>t_k \mid X_0 = x)
       = \PP(J_{t_1}(A_1)=0,\ldots,J_{t_k}(A_k)=0 \mid X_0 = x)
 \end{align}
 for all $x\in\RR_+$ and $t_1,\ldots,t_k\in\RR_+$.
Without loss of generality, we can assume that $t_1,\ldots,t_k\in\RR_{++}$, since if $t_i=0$ for some 
 $i\in\{1,\ldots,k\}$, then $\PP_x(\tau_{A_i}>0)=1$, $x\in\RR_+$ (following from \eqref{formula3_multiple})
 and $\PP_x(J_0(A_i)=0)=1$, $x\in\RR_+$ (due to the definition of $J_0(A_i)$).
If $t_i\in\RR_{++}$ and $\tau_{A_i}>t_i$, $i\in\{1,\ldots,k\}$, 
 then, by the definition of $J_{t_i}(A_i)$, $i\in\{1,\ldots,k\}$, we have that $J_{t_i}(A_i)=0$, $i\in\{1,\ldots,k\}$, 
 i.e., $\{\tau_{A_i}>t_i\}\subset \{J_{t_i}(A_i)=0\}$, $i\in\{1,\ldots,k\}$, 
 which, using the monotonicity of the probability, yields that
 \begin{align}\label{help203_versionk}
   \PP( \tau_{A_1}>t_1,\ldots,  \tau_{A_k}>t_k\mid X_0 = x)
     \leq \PP(J_{t_1}(A_1)=0,\ldots, J_{t_k}(A_k)=0\mid X_0 = x)
 \end{align}
 for all $x\in\RR_+$ and $t_1,\ldots,t_k\in\RR_{++}$.
If  $t_i\in\RR_{++}$ and $J_{t_i}(A_i)=0$, $i\in\{1,\ldots,k\}$, 
 then, by the definition of $\tau_{A_i}$, $i\in\{1,\ldots,k\}$, we have that
 $\tau_{A_i}\geq t_i$, $i\in\{1,\ldots,k\}$, 
 i.e., $\{J_{t_i}(A_i)=0\}\subset \{\tau_{A_i}\geq t_i\}$, $i\in\{1,\ldots,k\}$, 
 which, using again the monotonicity of the probability, yields that
 \begin{align*}
   &\PP(J_{t_1}(A_1)=0,\ldots, J_{t_k}(A_k)=0 \mid X_0 = x)\\
   &\qquad\leq \PP( \tau_{A_1}\geq t_1,\ldots, \tau_{A_k}\geq t_k, J_{t_1}(A_1)=0,\ldots, J_{t_k}(A_k)=0\mid X_0 = x)
  \end{align*} 
 for all $x\in\RR_+$ and $t_1,\ldots,t_k\in\RR_{++}$.
Hence, using that 
  \[
    \{ \tau_{A_1}\geq t_1,\ldots, \tau_{A_k}\geq t_k \}
      \subset  \{ \tau_{A_1}> t_1,\ldots, \tau_{A_k}> t_k \} \cup \bigcup_{i=1}^k \{\tau_{A_i}= t_i\},
  \]
 we get that
 \begin{align}\label{help204_versionka}
  \begin{split}
   &\PP(J_{t_1}(A_1)=0,\ldots, J_{t_k}(A_k)=0 \mid X_0 = x)\\
   &\qquad  \leq \PP( \tau_{A_1}>t_1,\ldots, \tau_{A_k}>t_k, J_{t_1}(A_1)=0,\ldots, J_{t_k}(A_k)=0 \mid X_0 = x)\\
      &\phantom{\qquad=\;}  + \PP\left(  \bigcup_{i=1}^k \big\{ \tau_{A_i}=t_i, J_{t_1}(A_1)=0,\ldots, J_{t_k}(A_k)=0 \big\} \,\Big\vert\, X_0 = x\right)\\
      &\qquad\leq \PP( \tau_{A_1}> t_1,\ldots, \tau_{A_k}> t_k, J_{t_1}(A_1)=0,\ldots, J_{t_k}(A_k)=0  \mid X_0 = x)\\
      &\phantom{\qquad=\;}  + \sum_{i=1}^k \PP( \tau_{A_i}=t_i, J_{t_1}(A_1)=0,\ldots, J_{t_k}(A_k)=0 \mid X_0 = x)\\
      &\qquad\leq \PP( \tau_{A_1}> t_1,\ldots, \tau_{A_k}> t_k \mid X_0 = x)
          + \sum_{i=1}^k \PP( \tau_{A_i}=t_i, J_{t_i}(A_i)=0\mid X_0 = x)
 \end{split}         
 \end{align}
 for all $x\in\RR_+$ and $t_1,\ldots,t_k\in\RR_{++}$.
 
 Now, we verify that, for all $A\in\cB((0,\infty))$ with $\nu(A)+\mu(A)<\infty$, we have that
 \begin{align}\label{help205}
    \PP( \tau_{A} = t, J_t(A)=0 \mid X_0 = x) =0 
 \end{align}
  for all $x\in\RR_+$ and $t\in\RR_{++}$.
Using that $\PP(J_s(A) < \infty)=1$, $s\in\RR_+$ (see the beginning of the proof),
 we have that, for all $t\in\RR_{++}$, the probability of the event $S\in \cF$ that there exists a sequence $(t_n)_{n\in\NN}$
 of pairwise distinct positive real numbers such that $t_n>t$, $n\in\NN$,
 $t_n\to t$ as $n\to\infty$ and $\Delta X_{t_n}\in A$, $n\in\NN$, is zero
 (since if such a sequence $(t_n)_{n\in\NN}$ exists, then $J_{t+\vare}(A) = \infty$ for any $\vare>0$).
Hence, for all $t\in\RR_{++}$, we have $\{\tau_{A} = t\} \subseteq \{ \Delta X_t\in A\}\cup S$, and therefore
 \begin{align*}
    \PP( \tau_{A} = t, J_t(A)=0 \mid X_0 = x) 
    &\leq \PP( \Delta X_t\in A, J_t(A)=0 \mid X_0 = x)\\
    &\leq \PP( J_t(A)\geq 1, J_t(A)=0 \mid X_0 = x)=0
 \end{align*}
  for all $x\in\RR_+$ and $t\in\RR_{++}$, which yields \eqref{help205}.
 
Taking into account \eqref{help204_versionka} and \eqref{help205}, we get that
 \[
  \PP(J_{t_1}(A_1)=0,\ldots, J_{t_k}(A_k)=0 \mid X_0 = x) 
   \leq  \PP( \tau_{A_1}> t_1,\ldots, \tau_{A_k}> t_k \mid X_0 = x)
 \]
 for all $x\in\RR_+$ and $t_1,\ldots,t_k\in\RR_{++}$.
Therefore, using also \eqref{help203_versionk}, we get \eqref{help202_versionk}, as desired.

{\sl Step 2.}
We express the right hand side of \eqref{help202_versionk} in terms of the limit of some Laplace transforms.
By the dominated convergence theorem and  the fact that $\sum_{i=1}^kJ_{t_i}(A_i)$ is a non-negative integer-valued random variable,
 we have that
 \begin{align}\label{help201_versionk}
  \begin{split}
  &\lim_{\lambda\to\infty} \EE\left( \ee^{-\lambda \sum_{i=1}^k J_{t_i}(A_i)}  \mid X_0 = x \right)\\
%   & = \lim_{\lambda\to\infty}
%      \EE\left( \ee^{- \lambda \sum_{i=1}^k J_{t_i}(A_i) }\bbone_{\{ \sum_{i=1}^k J_{t_i}(A_i) \ne 0\}} 
%           + \ee^{-\lambda\sum_{i=1}^k J_{t_i}(A_i) }\bbone_{\{ \sum_{i=1}^k J_{t_i}(A_i) = 0\}}
%                 \mid X_0 = x\right)\\
   & = \lim_{\lambda\to\infty}
              \EE\left( \ee^{- \lambda \sum_{i=1}^k J_{t_i}(A_i) }\bbone_{\{ \sum_{i=1}^k J_{t_i}(A_i)\ne 0\}} \mid X_0 = x\right)
             + \PP\left(\sum_{i=1}^k J_{t_i}(A_i)=0 \mid X_0 = x\right)\\
   & = \PP\left(\sum_{i=1}^k J_{t_i}(A_i)=0 \mid X_0 = x\right)\\
   &= \PP\left( J_{t_1}(A_1)=0,\ldots, J_{t_k}(A_k)=0 \mid X_0 = x\right),
        \qquad x\in\RR_+.
 \end{split}
 \end{align}

By \eqref{help202_versionk} and \eqref{help201_versionk}, we get that 
 \begin{equation}\label{eq:limittaus}
  \PP( \tau_{A_1}>t_1,\ldots,\tau_{A_k}>t_k \mid X_0 = x)
   = \lim_{\lambda\to\infty} \EE\left( \ee^{-\lambda \sum_{i=1}^k J_{t_i}(A_i)}  \mid X_0 = x \right),
   \qquad  x\in\RR_+.
 \end{equation} 
 
{\sl Step 3.}
We derive an expression for the joint Laplace transform
 \begin{equation*}\label{eq:Laplace}
   \EE\left( \ee^{-\lambda \sum_{i=1}^k J_{t_i}(A_i)}  \mid X_0 = x \right),\qquad \lambda\in\RR_+,\;\; x\in\RR_+.
 \end{equation*}  
We recall that, by Theorem \ref{Thm_(k+1)CBI} and parts (i) and (ii) of Corollary \ref{Cor_Laplace},  
 $\bX_t:=(X_t,J_t(A_1),\ldots,J_t(A_k))$, $t\in\RR_+$, is a $(k+1)$-type CBI process starting from $(X_0,\bzero)$,
  with branching mechanism $\bvarphi^{(k)}$ given by \eqref{2CBI_branching} and immigration mechanism $\psi^{(k)}$ given by \eqref{2CBI_immigration}.
Since
  \[
    J_{t_i}(A_i)  = \left\langle \be^{(k+1)}_{i+1},
                                   \begin{pmatrix}
                                     X_{t_i},
                                     J_{t_i}(A_1),
                                     \ldots,
                                     J_{t_i}(A_k) 
                                   \end{pmatrix}
                                   \right\rangle, \qquad i\in\{1,\ldots,k\},
 \]
  we have that
 \begin{align*}
   \EE\big( \ee^{-\lambda \sum_{i=1}^k J_{t_i}(A_i)}  \mid X_0 = x \big)
     = \EE( \ee^{-\sum_{i=1}^k \langle \blambda_i,\bX_{t_i} \rangle }  \mid \bX_0={(x,\bzero)}),
     \qquad  \lambda\in\RR_+,\;\; x\in\RR_+,
 \end{align*}  
 where $\blambda_i:=\lambda \be^{(k+1)}_{i+1}\in\RR_+^{k+1}$, $i\in\{1,\ldots,k\}$ 
 for any $\lambda\in\RR_+$. 
Then, by Proposition \ref{Pro_CBI_joint_Laplace}, we get that
 \begin{equation} \label{eq:laplace_with_v}
 \begin{split}
 & \EE\big( \ee^{-\lambda \sum_{i=1}^k J_{t_i}(A_i)}  \mid X_0 = x \big)\\
 &= \exp\Bigg\{-\big\langle  (x,\bzero),  \bv^{(k)}\big(t_1, \by^{(k)}_1\big) \big\rangle  - \sum_{i=1}^k \int_0^{t_i-t_{i-1}} \psi^{(k)}\big(  \bv^{(k)}(u,\by^{(k)}_{i}) \big)\,\dd u\Bigg\},
 \end{split}
 \end{equation}
 where
  \begin{itemize}
   \item for any $\blambda\in\RR_+^{k+1}$,  the function
         $\RR_+\ni t\mapsto \bv^{(k)}(t, \blambda)\in \RR_+^{k+1}$ is the unique locally bounded solution to the system 
           \eqref{help12} associated with the branching mechanism $\bvarphi^{(k)}$ given in \eqref{2CBI_branching},
    \item $\by^{(k)}_k:=\blambda_k \in \RR_+^{k+1}$ and the vectors $\by^{(k)}_i\in\RR_+^{k+1}$, 
          $i\in\{1,\ldots,k-1\}$, are given by \eqref{help_y_recursion},
    \item $\psi^{(k)}$ is given by \eqref{2CBI_immigration}.   
  \end{itemize}
Note that the vectors $\by^{(k)}_i$, $i\in\{1,\ldots,k\}$ depend on $\lambda\in\RR_+$ as well, but for simplifying the notations,
 we do not denote this dependence.
In what follows, by a kind of backward induction argument, 
 we show that 
 \begin{equation}\label{eq: def y_i^k} 
  \by_k^{(k)}=(0,\lambda \biota_k), \qquad \qquad 
  \by_{i}^{(k)}= (v^{(k)}(t_{i+1}-t_{i},\by_{i+1}^{(k)}),\lambda \biota_i),  \quad i\in\{1,\ldots,k-1\}
 \end{equation}
 for any $\lambda\in\RR_+$, where $\biota_i:=\sum_{j=i}^{k} \be_{j}^{(k)}\in \RR_+^{k}$, $i\in\{1,\ldots,k\}$, and, for any $\blambda\in \RR_+^{k+1}$, 
 the function $\RR_+\ni t\mapsto v^{(k)}(t,\blambda)$ is the unique locally bounded non-negative solution to the differential equation
 \eqref{def_v_k+1_CBI}. 
By the definitions of $\blambda_k$ and $\biota_k$, we have that 
 \[
   \by_k^{(k)}= \blambda_k=\lambda\be^{(k+1)}_{k+1}=(0,\lambda\biota_k)\in\RR_+^{k+1},
    \qquad \lambda\in\RR_+.
 \]  
Therefore, using also \eqref{help_y_recursion} and part (iii) of Corollary \ref{Cor_Laplace}
 with $(\lambda_0,\ldots,\lambda_k)=\by^{(k)}_k$, we have that
\begin{align*}
\by_{k-1}^{(k)}&=\blambda_{k-1}+\bv^{(k)}(t_k-t_{k-1},\by_{k}^{(k)})=\lambda\be_{k}^{(k+1)}+v^{(k)}(t_k-t_{k-1},\by_{k}^{(k)})\be_{1}^{(k+1)}+ \lambda \be_{k+1}^{(k+1)}\\
               &=(v^{(k)}(t_k-t_{k-1},\by_{k}^{(k)}), \lambda\biota_{k-1}), \qquad \lambda\in\RR_+.
\end{align*}
Thus the second formula in \eqref{eq: def y_i^k} holds for $i=k-1$. 
Now, suppose that the second formula in \eqref{eq: def y_i^k} holds for $i=m$, where $m\in\{2,\ldots,k-1\}$.
Using again \eqref{help_y_recursion}, part (iii) of Corollary \ref{Cor_Laplace} with $(\lambda_0,\ldots,\lambda_k)=\by^{(k)}_m$,
 and that the last $k-m+1$ coordinates of $\by_m^{(k)}$ are $\lambda$, and the other coordinates are zero except
 the first coordinate, we have that
\begin{align*}
\by_{m-1}^{(k)}&=\blambda_{m-1}+\bv^{(k)}(t_m-t_{m-1},\by_{m}^{(k)})
                =\lambda\be_{m}^{(k+1)} + v^{(k)}(t_{m}-t_{m-1},\by_{m}^{(k)})\be_{1}^{(k+1)}+ \sum_{j=m}^{k} \lambda \be_{j+1}^{(k+1)}\\
               &=v^{(k)}(t_{m}-t_{m-1},\by_{m}^{(k)})\be_{1}^{(k+1)} + \sum_{j=m-1}^{k} \lambda \be_{j+1}^{(k+1)}
                      =(v^{(k)}(t_{m}-t_{m-1},\by_{m}^{(k)}), \lambda\biota_{m-1}),
\end{align*}
 i.e., the second formula in \eqref{eq: def y_i^k} holds for $i=m-1$ as well.
This implies that the second formula in \eqref{eq: def y_i^k} holds for every $i\in\{1,\dots,k-1\}$, as desired.

Furthermore, using \eqref{eq: def y_i^k} and part (iii) of Corollary \ref{Cor_Laplace} 
 with $(\lambda_0,\ldots,\lambda_k)=\by^{(k)}_i$, $i\in\{1,\ldots,k\}$,  we have that
 \begin{equation}\label{eq:vaux}
 \bv^{(k)}(t,\by_{i}^{(k)})=(v^{(k)}(t,\by_{i}^{(k)}), \lambda \biota_i),
   \qquad  i\in\{1,\dots, k\}, \qquad t,\lambda\in\RR_+.
\end{equation}
  
{\sl Step 4.} By a kind of backward induction, we prove
 that, for all $t\in\RR_+$ and $i\in\{1,\dots,k\}$, it holds that
\begin{equation}\label{eq:relvw}
  v^{(k)}(t,\by_{i}^{(k)})\uparrow \tw^{(k)}_i(t) \qquad \text{as $\lambda\uparrow\infty$,} 
\end{equation}
 where $\tw^{(k)}_k(t)=v^{(C_k)}(t,0)$ and $\tw^{(k)}_i(t)$ is given by \eqref{help_new_01} for $i\in\{1,\dots,k-1\}$,
 and we remind the readers that the vectors $\by_{i}^{(k)}$, $i\in\{1,\ldots,k\}$, depend on $\lambda$ as well.

First, we check that \eqref{eq:relvw} holds for $i=k$.
Since $\by_k^{(k)}=(0,\lambda \biota_k)$ (see \eqref{eq: def y_i^k}), 
 part (ii) of Lemma \ref{Lem_vjumps} with the function $\RR_+\ni \lambda \mapsto \eta_\lambda:=0$, $\lambda_0:=0$, and $\blambda:=\lambda \biota_k$ yields that 
 \[
    v^{(k)}(t,\by_k^{(k)})
      = v^{(k)}(t,0,\lambda \biota_k) \uparrow  v^{(A_{\biota_k})}(t,0)
      \qquad \text{as \ $\lambda\uparrow \infty$}
 \]   
 for all $t\in\RR_+$, where $A_{\biota_k} = A_k = C_k$ and the function $\RR_+\ni t\mapsto v^{(A_{\biota_k})}(t,0)$ is the unique
 locally bounded non-negative solution to the differential equation \eqref{eq_vAlambda}
 (with $\blambda = \biota_k$ and $\lambda_0=0$).
Therefore, since $\tw^{(k)}_k(t) = v^{(C_k)}(t,0)$, $t\in\RR_+$ (by definition),
 we obtain that \eqref{eq:relvw} holds for $i=k$.

Now suppose that \eqref{eq:relvw} holds for $m\in\{2,\dots,k\}$, that is, for all $t\in\RR_+$, it holds that
 \[
   v^{(k)}(t,\by_{m}^{(k)}) \uparrow \tw^{(k)}_m(t) \qquad \text{as \ $\lambda\uparrow \infty$.} 
 \]
In particular, for $m\in\{2,\dots,k\}$, it holds that 
 \[
    v^{(k)}(t_m-t_{m-1},\by_{m}^{(k)}) \uparrow \tw^{(k)}_m(t_m-t_{m-1}) \qquad \text{as \ $\lambda\uparrow \infty$.} 
 \]
We prove that \eqref{eq:relvw} holds for $i=m-1$ as well.
Then, by \eqref{eq: def y_i^k} and part (ii) of Lemma \ref{Lem_vjumps}  
 with the function $\RR_+\ni\lambda \mapsto \eta_\lambda:=v^{(k)}(t_m-t_{m-1},\by_m^{(k)})$, $\lambda_0:=\tw^{(k)}_m(t_m-t_{m-1})$
 and $\blambda:=\biota_{m-1}$, we get that
\begin{align}\label{help_main_2}
 \begin{split}
  v^{(k)}(t,\by_{m-1}^{(k)})&=v^{(k)}(t,v^{(k)}(t_{m}-t_{m-1},\by_{m}^{(k)}),\lambda \biota_{m-1}) \\
 &\; \uparrow v^{(A_{\biota_{m-1}})}(t,\tw^{(k)}_m(t_m-t_{m-1}))
    \qquad \text{as \ $\lambda\uparrow \infty$} 
\end{split}    
\end{align}
 for all $t\in\RR_+$.
Since $A_{\biota_{m-1}} = C_{m-1}$, using also \eqref{help_new_01}, we have that 
 \[
   v^{(A_{\biota_{m-1}})}(t,\tw^{(k)}_m(t_m-t_{m-1}))
      = \tw^{(k)}_{m-1}(t),\qquad t\in\RR_+,
 \] 
 which together with \eqref{help_main_2} yield that \eqref{eq:relvw} also holds for $i=m-1$, as desired.

{\sl Step 5.}
By \eqref{eq:limittaus} and  \eqref{eq:laplace_with_v}, we get that 
\begin{equation}\label{eq:S51}
\begin{split}
&\PP( \tau_{A_1}>t_1,\ldots,\tau_{A_k}>t_k \mid X_0 = x)
   = \lim_{\lambda\uparrow\infty} \EE\left( \ee^{-\lambda \sum_{i=1}^k J_{t_i}(A_i)}  \mid X_0 = x \right)\\
&\qquad = \lim_{\lambda\uparrow \infty}
 \exp\Bigg\{-\big\langle  (x,\bzero),  \bv^{(k)}\big(t_1, \by^{(k)}_1\big) \big\rangle  - \sum_{i=1}^k \int_0^{t_i-t_{i-1}} \psi^{(k)}\big(  \bv^{(k)}(u,\by^{(k)}_{i}) \big)\,\dd u\Bigg\} 
\end{split}
\end{equation}
  for all $x\in\RR_+$, where $\psi^{(k)}$ is given by \eqref{2CBI_immigration}.
By \eqref{eq:vaux}  and \eqref{eq:relvw}, for all $x\in\RR_+$, we have that 
 \begin{align}\label{help_main_1}
 \begin{split}
    \big\langle  (x,\bzero),  \bv^{(k)}\big(t_1, \by^{(k)}_1\big) \big\rangle
     & = \big\langle  (x,\bzero),  (v^{(k)}(t_1,\by^{(k)}_1),\lambda \biota_1 )  \big\rangle\\
     & = xv^{(k)}(t_1,\by^{(k)}_1)\uparrow x\tw^{(k)}_1(t_1)
      \qquad \text{as \ $\lambda\uparrow \infty$.}
 \end{split}     
 \end{align}
In what follows, let $i\in\{1,\dots,k\}$ be fixed, and recall the notation
 $\biota_i=\sum_{j=i}^{k} \be_{j}^{(k)} $. 
By \eqref{2CBI_immigration} and  \eqref{eq:vaux}, we have that 
\begin{align*}
 \psi^{(k)}\big(  \bv^{(k)}(u,\by^{(k)}_{i}) \big)
 = \beta v^{(k)}(u,\by^{(k)}_{i}) + \int_{\cU_1}\bigl( 1 - \ee^{-v^{(k)}(u,\by^{(k)}_{i})  r - \sum_{j=i}^k \lambda \bbone_{A_j}(r)} \bigr)\, \nu(\dd r),
  \qquad u\in\RR_+.
\end{align*}
By \eqref{eq:relvw}, we get that $v^{(k)}(t,\by_{i}^{(k)})\uparrow \tw^{(k)}_i(t)$ as $\lambda\uparrow \infty$ for all $t\in\RR_+$, 
 which yields that the non-negative function $\RR_+\ni\lambda\mapsto\psi^{(k)}\big(  \bv^{(k)}(u,\by^{(k)}_{i}) \big)\in\RR_+$ 
 is non-decreasing for all $u\in\RR_+$. 
Consequently, by the monotone convergence theorem and the non-negativity of $\psi^{(k)}$,
 we obtain that
\begin{align}\label{help_change_limit_psi}
  \begin{split}
&\lim_{\lambda\uparrow\infty}\int_{0}^{t_i-t_{i-1}}\psi^{(k)}\big(  \bv^{(k)}(u,\by^{(k)}_{i}) \big)\,\dd u 
  =\int_{0}^{t_i-t_{i-1}}\lim_{\lambda\uparrow\infty}\psi^{(k)}\big(  \bv^{(k)}(u,\by^{(k)}_{i}) \big)\,\dd u\\
&=\int_{0}^{t_i-t_{i-1}}\beta \tw^{(k)}_i(u) \,\dd u+\int_{0}^{t_i-t_{i-1}} \int_{\cU_1\setminus A_{\biota_i}}\bigl( 1 - \ee^{-\tw^{(k)}_i(u) r }\bigr) \nu(\dd r)\,\dd u
    +\nu(A_{\biota_i})(t_i-t_{i-1})\\
&= \int_{0}^{t_i-t_{i-1}} \psi^{(C_i)}(\tw^{(k)}_i(u))\,\dd u  + \nu(C_i)(t_i-t_{i-1}) 
       = \int_{0}^{t_i-t_{i-1}}   \widetilde\psi^{(i)}(\tw^{(k)}_i(u))\,\dd u,
  \end{split}
 \end{align}
 where, at the second step, we used that $\sum_{j=i}^k \lambda \bbone_{A_j}(r)=0$ for all $r\in\cU_1\setminus A_{\biota_i}$, 
 and $\sum_{j=i}^k \lambda \bbone_{A_j}(r) \to\infty$ as $\lambda \uparrow \infty$ for all $r\in A_{\biota_i}$,
 and, at the third step, $A_{\biota_i}=C_i$, the definition of $\psi^{(C_i)}$ and the non-negativity of $\tw^{(k)}_i(u)$, $u\in\RR_+$.
Therefore, by \eqref{eq:S51}, \eqref{help_main_1} and \eqref{help_change_limit_psi}, for all $x\in\RR_+$, we have that
 \begin{align*}
   &\PP( \tau_{A_1}>t_1,\ldots,\tau_{A_k}>t_k \mid X_0 = x)
 =\exp\left\{  -x \tw^{(k)}_1(t_1)-\sum_{i=1}^{k}\int_{0}^{t_i-t_{i-1}}\widetilde{\psi}^{(i)}\big(\tw^{(k)}_i(u)\big)\,\dd u \right\},
 \end{align*}
 as desired.
\proofend

In the next Corollary \ref{Cor_multiple_jump_times_CBI}, we present an expression for the joint distribution function of $(\tau_{A_1},\ldots,\tau_{A_k})$ 
 at a point $(t,\ldots,t)\in\RR_+^k$.
Note also that this corollary in the special case $k=1$ gives back Theorem \ref{Thm_single_jump_times_CBI} as well.
We will give two proofs of Corollary \ref{Cor_multiple_jump_times_CBI}.
The first one is a direct application of the known result of He and Li \cite[Theorem 3.1]{HeLi} (see also Theorem \ref{Thm_single_jump_times_CBI}),
 and the second one is a specialization of Theorem \ref{Thm_main} to the case $t_1=\cdots=t_k=t\in\RR_+$.
The first proof of Corollary \ref{Cor_multiple_jump_times_CBI} also demonstrates that 
 Theorem \ref{Thm_single_jump_times_CBI} directly implies a formula for the joint distribution function 
 of $(\tau_{A_1},\ldots,\tau_{A_k})$ at a point $(t,\ldots,t)\in\RR_+^k$, for this we do not need our 
 general Theorem \ref{Thm_main}.

\begin{Cor}\label{Cor_multiple_jump_times_CBI}
Let $(X_t)_{t\in\RR_+}$ be a CBI process with parameters $(1, c, \beta, B, \nu, \mu)$ such that
 $\EE(X_0) < \infty$ and the moment condition \eqref{moment_condition_m_new} with $d=1$ hold.
Let $k\in\NN$ and $A_1,\ldots,A_k\in\cB((0,\infty))$ be such that $\nu(A_i) + \mu(A_i)<\infty$, $i\in\{1,\ldots,k\}$.
Then we have
 \begin{align}\label{formula_multiple_same_times}
  \begin{split}
  \PP_x(\tau_{A_1}>t,\ldots,\tau_{A_k}>t) 
         = \exp\left\{ -\nu(A) t - x \tv^{(A)}(t)
                       - \int_0^t \psi^{(A)}( \tv^{(A)}(s) )\,\dd s \right\}
 \end{split}
 \end{align}
 for all $t\in\RR_+$ and $x\in\RR_+$, where $A:=\bigcup_{i=1}^k A_i$,
 and the continuously differentiable function
 $\RR_+\ni t \mapsto \tv^{(A)}(t) \in\RR_+$ is the unique locally bounded solution to the differential equation \eqref{DE_jump_time}.
\end{Cor}

\noindent{\bf First proof.}
It is based on Theorem \ref{Thm_single_jump_times_CBI}.
We check that
\begin{equation}\label{eq:taus}
\{\tau_{A_1}>t,\ldots,\tau_{A_k}>t\}= \{\tau_{A}>t\},\qquad t\in\RR_+.
\end{equation}
Let $t\in\RR_+$ be fixed.
First, suppose that $\tau_{A_i}>t$ for all $i\in\{1,\dots,k\}$. 
Let $h>0$ be such that $\tau_{A_i}>t+h$ for all $i\in\{1,\dots,k\}$. 
Then, by the definition of $\tau_{A_i}$, $i\in\{1,\ldots,k\}$, we have that 
$J_{t+h}(A_i)=0$ for all $i\in\{1,\dots,k\}$ implying that $J_{t+h}(A)=0$. 
Then $\tau_{A}\geq t+h>t$, yielding that $\{\tau_{A_1}>t,\ldots,\tau_{A_k}>t\} \subset \{\tau_{A}>t\}$.
   
Suppose now that $\tau_{A}>t$, and let $h>0$ be such that $\tau_{A}>t+h$. 
Then, by the definition of $\tau_A$, we have that $J_{t+h}(A)=0$. 
This implies that $J_{t+h}(A_i)=0$ for all $i\in\{1,\dots,k\}$. 
Hence $\tau_{A_i}\geq t+h>t$ for each $i\in\{1,\dots,k\}$,
 yielding that $\{\tau_{A}>t\}\subset \{\tau_{A_1}>t,\ldots,\tau_{A_k}>t\}$.
Then \eqref{eq:taus} holds, and the result directly follows by Theorem \ref{Thm_single_jump_times_CBI}.

\smallskip

\noindent{\bf Second proof.}
It is based on Theorem \ref{Thm_main}.
By choosing $t_1:=\cdots:=t_k:=t\in\RR_+$ in Theorem \ref{Thm_main}, we get that
 \begin{align}\label{formula_multiple_diff_times_spec}
  \begin{split}
  \PP_x(\tau_{A_1}>t,\ldots,\tau_{A_k}>t) 
  &= \exp\left\{ -x \tw^{(k)}_1(t) - \int_{0}^t\widetilde{\psi}^{(1)}\big( \tw^{(k)}_1(u)\big)\dd u   \right\}\\
  &= \exp\left\{ -\nu(C_1)t -x\tw^{(k)}_1(t) - \int_{0}^t \psi^{(C_1)}\big( \tw^{(k)}_1(u)\big)\dd u   \right\}
 \end{split}
 \end{align} 
 for all $x,t\in\RR_+$, where $C_1=\bigcup_{j=1}^k A_j =A$.
Therefore, in order to see that the right hand side of \eqref{formula_multiple_diff_times_spec} indeed coincides 
 with the right hand side of \eqref{formula_multiple_same_times}, it is enough to check that 
 \[
   \tw^{(k)}_1(t) = \tv^{(C_1)}(t),\qquad t\in\RR_+.
 \]
Here $\tw^{(k)}_k(t) = v^{(C_k)}(t,0)$, $t\in\RR_+$ (by definition), and, for $i\in\{1,\ldots,k-1\}$, using also \eqref{help_new_01} and \eqref{help_DE_main_result},
 we get that
 \begin{align}\label{help_multiple_diff_times_spec_1}
 \begin{split}
  \tw^{(k)}_1(t)
  & = v^{(C_1)}(t, \tw^{(k)}_2(0))
   = v^{(C_1)}\big(t, v^{(C_2)}(0, \tw^{(k)}_3(0)) \big)\\
  & = v^{(C_1)}(t,\tw^{(k)}_3(0))
   = \cdots
   = v^{(C_1)}(t,\tw^{(k)}_k(0))\\
  & = v^{(C_1)}(t,v^{(C_k)}(0,0))
    = v^{(C_1)}(t,0),\qquad t\in\RR_+.
  \end{split}  
 \end{align}
Further, by the uniqueness of a locally bounded non-negative solution to the differential equations
 \eqref{DE_jump_time} with $A=C_1$ and \eqref{help_DE_main_result} with $i=1$ and $\lambda_0=0$,
 we get that $v^{(C_1)}(t,0) = \tv^{(C_1)}(t)$, $t\in\RR_+$.
This together with \eqref{help_multiple_diff_times_spec_1} yield that
 $\tw^{(k)}_1(t) = \tv^{(C_1)}(t)$, $t\in\RR_+$, as desired.
\proofend

In the next remark, we mention a possible application of Theorem \ref{Thm_main} in financial mathematics.

\begin{Rem}\label{Rem_financial_application}
Suppose that a treasury desk monitors interest-rate modelled by an $\alpha$-stable CIR process $(X_t)_{t\in\RR_+}$, 
 which is given by the pathwise unique non-negative strong solution of the SDE
 \begin{align}\label{stable_CIR}
  \dd X_t = (a - b X_t) \, \dd t + \sigma \sqrt{X_t} \, \dd W_t
            + \delta \sqrt[\alpha]{X_{t-}} \, \dd L_t ,
  \qquad  t \in\RR_+,
 \end{align}
 with an almost surely non-negative initial value $X_0$ satisfying $\EE(X_0)<\infty$,
 where $a \in\RR_+$, $b \in \RR$, $\sigma \in \RR_+$, $\delta \in (0, \infty)$, $\alpha \in (1, 2)$,
 $(W_t)_{t \in\RR_+}$ is a standard Wiener process, and
 $(L_t)_{t\in\RR_+}$ is a spectrally positive $\alpha$-stable L\'evy
 process such that the characteristic function of $L_1$ takes the form
 \begin{align*}
  \EE(\ee^{\ii\theta L_1})
  = \exp\left\{(\alpha \Gamma(-\alpha))^{-1} \int_0^\infty
                (\ee^{\ii \theta z} - 1 - \ii \theta z)  z^{-1-\alpha}
                \, \dd z\right\} ,
  \qquad \theta \in \RR ,
 \end{align*}
 where $\Gamma$ denotes the Gamma function, and we suppose that $X_0$, $(W_t)_{t\in\RR_+}$ and $(L_t)_{t\in\RR_+}$ are independent
 (see Fu and Li \cite[Corollary 6.3]{FuLi} and Jiao et al.\ \cite{JiaoMaScotti}).
It is known that the $\alpha$-stable CIR process $(X_t)_{t\in\RR_+}$ given by the SDE \eqref{stable_CIR} is a
 CBI process with branching mechanism
 \[
      \varphi(\lambda) = \frac{\sigma^2}{2} \lambda^2  + b\lambda + \frac{\delta^\alpha}{\alpha}\lambda^\alpha,
      \qquad \lambda \in \RR_+ ,
    \]
  and immigration mechanism $\psi(\lambda) = a \lambda$, $\lambda \in \RR_+$,
 see, e.g., part (ii) of Proposition 2.1 in Barczy et al.\ \cite{BarBenKebPap}.
Suppose that a treasury desk's short-term liquidity arrangements may require that no
 jump larger than $25$ basis points occur before the next monthly review,
 while its long-term arrangements entail that no jump larger than $100$
 basis points occur during the next year. 
If the short interest rates are measured as a decimal annual rate and time is measured in years, 
 then these requirements correspond to $A_1=(0.0025,\infty)$, $A_2=(0.01,\infty)$, $t_1=1/12$ and $t_2=1$
 (with the notations of Theorem \ref{Thm_main}).
The probability that both requirements are satisfied is
 $\PP(\tau_{A_1}>t_1,\tau_{A_2}>t_2\mid X_0=x)$, where $x$ is the (initial) short interest rate at time $0$.
Here we emphasize that the thresholds $0.0025$ and $0.01$ concern the jump sizes of the short interest rates,
 rather than the value of the short interest rate itself.
Theorem \ref{Thm_main} can be used to calculate the probability $\PP(\tau_{A_1}>t_1,\tau_{A_2}>t_2\mid X_0=x)$.
\proofend
\end{Rem}
    
\vspace*{5mm}

\appendix

\vspace*{5mm}

\noindent{\bf\Large Appendix}

\section{Joint Laplace transform of multi-type CBI processes}\label{Appendix_multitype_CBI}

In this appendix, we derive a formula for the joint Laplace transform of a multi-type CBI process.

\begin{Pro}\label{Pro_CBI_joint_Laplace}
Let $(\bX_t)_{t\in\RR_+}$ be a multi-type CBI process with parameters $(d, \bc, \Bbeta, \bB, \nu, \bmu)$ 
 such that the moment condition \eqref{moment_condition_m_new} holds. 
Let $\bx\in\RR_+^d$, $k\in\NN$, $0:=t_0\leq t_1\leq t_2\leq\ldots\leq t_k$ and $\blambda_i\in\RR_+^d$, $i\in\{1,\ldots,k\}$.
Then we have that
 \begin{align*}
   \EE( \ee^{-\sum_{i=1}^k \langle \blambda_i,\bX_{t_i} \rangle }  \mid \bX_0=\bx) 
      = \exp\Bigg\{ -\big\langle \bx, \bv\big(t_1, \by^{(k)}_1\big) \big\rangle 
- \sum_{i=1}^k \int_0^{t_i-t_{i-1}} \psi\big(\bv(u,\by^{(k)}_{i}) \big)\,\dd u
                        \Bigg\},
 \end{align*}
 where 
 \begin{itemize}
 \item for any $\blambda \in \RR_+^d$, the continuously differentiable
       function $\RR_+ \ni t \mapsto \bv(t, \blambda) \in \RR_+^d$ is the unique locally bounded solution
       to the system \eqref{help12},
 \item the vectors $\by^{(k)}_1, \dots, \by^{(k)}_k\in\RR_+^d$ are defined by
       \begin{align}\label{help_y_recursion}
       \by^{(k)}_k:=\blambda_k\qquad \mbox{ and } \qquad  \by^{(k)}_i :=\blambda_{i} + \bv\big(t_{i+1}-t_{i},\by^{(k)}_{i+1}\big),
          \qquad i\in\{1,\ldots,k-1\},      
       \end{align}
 \item $\psi$ is given by \eqref{dCBI_immigration}.
 \end{itemize} 
\end{Pro}

\noindent{\bf Proof.}
We prove the statement by induction on $k$.
The statement trivially holds for $k=1$ due to Theorem \ref{CBI_exists}.
Let $k\in\NN$, $0=t_0\leq t_1\leq t_2\leq\ldots\leq t_k\leq t_{k+1}$ and $\blambda_i\in\RR_+^d$, $i\in\{1,\ldots,k+1\}$.

First, assume that $t_k<t_{k+1}$.
Using that $(\bX_t)_{t\in\RR_+}$ is a time-homogeneous Markov process, 
 by the tower rule for conditional expectation, we have that
 \begin{align*}
   \EE( \ee^{- \sum_{i=1}^{k+1} \langle \blambda_i,\bX_{t_i} \rangle }  \mid \bX_0)
    & = \EE\Big( \EE( \ee^{- \sum_{i=1}^{k+1} \langle \blambda_i,\bX_{t_i} \rangle } \mid \bX_0,\bX_{t_1},\ldots,\bX_{t_k}) \mid \bX_0\Big)\\
    & = \EE\big( \ee^{-\sum_{i=1}^k \langle \blambda_i,\bX_{t_i} \rangle} \EE( \ee^{- \langle \blambda_{k+1}, \bX_{t_{k+1}}\rangle} \mid \bX_0,\bX_{t_1},\ldots,\bX_{t_k}) \mid \bX_0\big) \\
    & = \EE\big( \ee^{-\sum_{i=1}^k \langle \blambda_i,\bX_{t_i} \rangle} \EE( \ee^{-  \langle \blambda_{k+1}, \bX_{t_{k+1}}\rangle } \mid \bX_{t_k} ) \mid \bX_0\big).
 \end{align*}
Hence, using Theorem \ref{CBI_exists}, we get that
 \begin{align}
\label{help_CBI_joint_Laplace11}
 \begin{split}
   &\EE( \ee^{-\sum_{i=1}^{k+1} \langle \blambda_i,\bX_{t_i} \rangle}  \mid \bX_0)\\
   &= \EE\left( \ee^{-\sum_{i=1}^k \langle \blambda_i,\bX_{t_i} \rangle}
               \exp\left\{ -\big \langle \bX_{t_k} ,\bv(t_{k+1}-t_k,\blambda_{k+1})\big\rangle 
                           - \int_0^{t_{k+1}-t_k} \psi(\bv(u,\blambda_{k+1})) \,\dd u \right\} \,\Big\vert\, \bX_0\right) \\
   & = \exp\left\{ - \int_0^{t_{k+1}-t_k} \psi(\bv(u,\blambda_{k+1})) \,\dd u  \right\} \\
   &\phantom{=\;\;}   \times\EE\left( \exp\left\{ - \sum_{i=1}^{k-1} \langle \blambda_i,\bX_{t_i} \rangle 
                            - \big\langle \blambda_k + \bv(t_{k+1}-t_k,\blambda_{k+1}), \bX_{t_k}\big\rangle \right\} \,\Big\vert\, \bX_0\right).
\end{split}
 \end{align}
By the induction hypothesis, using also that $\bv(t,\blambda)\in\RR_+^d$ for any $t\in\RR_+$ and $\blambda\in\RR_+^d$, we have that 
 \begin{align}\label{help_CBI_joint_Laplace2}
 \begin{split}
  &\EE\left( \exp\left\{ - \sum_{i=1}^{k-1} \langle \blambda_i,\bX_{t_i} \rangle 
                            - \big\langle \blambda_k + \bv(t_{k+1}-t_k,\blambda_{k+1}), \bX_{t_k}\big\rangle \right\} \,\Big\vert\, \bX_0\right)\\
  &\qquad = \exp\left\{ - \langle \bX_0, \bv(t_1,\widetilde\by_1) \rangle 
                           - \sum_{i=1}^k \int_0^{t_i-t_{i-1}} \psi(\bv(u,\widetilde\by_{i})) \,\dd u \right\},
  \end{split}                         
 \end{align}
 where  $\widetilde\by_k:=\blambda_k + \bv(t_{k+1}-t_k,\blambda_{k+1})$ and 
 \[
 \widetilde\by_i:= \blambda_{i} + \bv(t_{i+1}-t_{i},\widetilde\by_{i+1}),
    \qquad i\in\{1,\ldots,k-1\}.
 \]
Observe that 
$\by^{(k+1)}_{k+1}=\blambda_{k+1}$ and $\by^{(k+1)}_{k}=\blambda_{k} + \bv(t_{k+1}-t_{k},\blambda_{k+1})=\widetilde\by_k$.
Moreover, by a kind of backward induction argument, we get that
 \[ 
 \by^{(k+1)}_{i}=\widetilde\by_i, \qquad  i\in\{1,\ldots,k\}. 
 \]
Indeed, if we suppose that $ \by^{(k+1)}_{i}=\widetilde\by_i$ for some $i\in\{2,\ldots,k\}$, then
 \[ 
  \by^{(k+1)}_{i-1}=\blambda_{i-1} + \bv(t_{i}-t_{i-1},\by^{(k+1)}_{i})= \blambda_{i-1} + \bv(t_{i}-t_{i-1},\widetilde\by_i)=\widetilde\by_{i-1}.
 \] 

Consequently, by \eqref{help_CBI_joint_Laplace11} and \eqref{help_CBI_joint_Laplace2}, we obtain that 
 \begin{align*}
   \EE( \ee^{-\sum_{i=1}^{k+1} \langle \blambda_i,\bX_{t_i} \rangle}  \mid \bX_0)
   & = \exp\Bigg\{ - \int_0^{t_{k+1}-t_k}  \psi\big(\bv(u, \by^{(k+1)}_{k+1})\big) \,\dd u 
                   - \langle \bX_0, \bv\big(t_1,\by^{(k+1)}_{1}\big) \rangle \\
   &\phantom{= \exp\Bigg\{\; } 
     - \sum_{i=1}^k \int_0^{t_i-t_{i-1}} \psi\big(\bv(u,\by^{(k+1)}_{i})\big) \,\dd u  \Bigg\}\\
   & = \exp\left\{ - \big\langle \bX_0, \bv\big(t_1,\by^{(k+1)}_{1}\big) \big\rangle  
                   - \sum_{i=1}^{k+1} \int_0^{t_i-t_{i-1}} \psi\big(\bv(u,\by^{(k+1)}_{i})\big) \,\dd u  
        \right\},
 \end{align*}
 as desired.
 
Now, assume that $t_k=t_{k+1}$.
Then, by the induction hypothesis, we get that 
 \begin{align}\label{help_AppA_1}
  \begin{split}
  \EE( \ee^{-\sum_{i=1}^{k+1} \langle \blambda_i,\bX_{t_i} \rangle}  \mid \bX_0)
     &=\EE( \ee^{-\sum_{i=1}^{k-1} \langle \blambda_i,\bX_{t_i} \rangle - \langle \blambda_k+\blambda_{k+1},\bX_{t_k} \rangle }  \mid \bX_0)\\
     &= \exp\left\{ - \langle \bX_0, \bv(t_1,\widehat\by_1) \rangle 
                           - \sum_{i=1}^k \int_0^{t_i-t_{i-1}} \psi(\bv(u,\widehat\by_{i})) \,\dd u \right\},
  \end{split}   
 \end{align}
 where $\widehat\by_k:=\blambda_k + \blambda_{k+1}$ and 
 \[
 \widehat\by_i:= \blambda_{i} + \bv(t_{i+1}-t_{i},\widehat\by_{i+1}),
    \qquad i\in\{1,\ldots,k-1\}.
 \]
Observe that $\by^{(k+1)}_{k+1}=\blambda_{k+1}$ and 
 \[
  \by^{(k+1)}_{k}=\blambda_{k} + \bv(t_{k+1}-t_{k},\blambda_{k+1})=\blambda_{k} + \bv(0,\blambda_{k+1})
                 = \blambda_{k} + \blambda_{k+1}
                 = \widehat y_k.
  \]
Moreover, by a kind of backward induction argument, we get that
 \[ 
 \by^{(k+1)}_{i}=\widehat\by_i, \qquad  i\in\{1,\ldots,k\}. 
 \]
Indeed, if we suppose that $\by^{(k+1)}_{i}=\widehat\by_i$ for some $i\in\{2,\ldots,k\}$, then
 \[ 
  \by^{(k+1)}_{i-1}=\blambda_{i-1} + \bv(t_{i}-t_{i-1},\by^{(k+1)}_{i})= \blambda_{i-1} + \bv(t_{i}-t_{i-1},\widehat\by_i)=\widehat\by_{i-1}.
 \] 
Consequently, by \eqref{help_AppA_1}, 
 \begin{align*}
 \EE( \ee^{-\sum_{i=1}^{k+1} \langle \blambda_i,\bX_{t_i} \rangle}  \mid \bX_0)
   & = \exp\left\{ - \langle \bX_0, \bv(t_1,\by^{(k+1)}_1) \rangle 
                           - \sum_{i=1}^k \int_0^{t_i-t_{i-1}} \psi(\bv(u,\by^{(k+1)}_{i})) \,\dd u \right\} \\
   & =  \exp\left\{ - \langle \bX_0, \bv(t_1,\by^{(k+1)}_1) \rangle 
                           - \sum_{i=1}^{k+1} \int_0^{t_i-t_{i-1}} \psi(\bv(u,\by^{(k+1)}_{i})) \,\dd u \right\}                            
 \end{align*} 
 as desired, where at the second equality, we used that 
 $\int_0^{t_{k+1}-t_k} \psi(\bv(u,\by^{(k+1)}_{k+1})) \,\dd u=0$ due to that $t_{k+1}=t_k$.
\proofend

In the next Proposition \ref{Pro_Li_CB_joint_Laplace}, 
 we present a formula for the joint Laplace transform of a single-type CBI process.
For corresponding results on superprocesses (without immigration)
 and single-type CB processes, see Li \cite[Proposition 5.14]{Li} and Li \cite[Proposition 4.1 and Theorem 4.2]{Li3},
 respectively.
Furthermore, we remark that part (i) of Proposition \ref{Pro_Li_CB_joint_Laplace} can be viewed as a special case 
 of Theorem 9.22 in Li \cite{Li}, which gives a characterization of the weighted occupation times 
 of immigration superprocesses in terms of Laplace transforms, 
 since CBI processes are special cases of immigration superprocesses (see Example 12.2 in Li \cite{Li}).
Concerning part (ii) of Proposition \ref{Pro_Li_CB_joint_Laplace}, 
 we cannot address a direct reference, and therefore we give a proof of Proposition \ref{Pro_Li_CB_joint_Laplace},
 and we point out the fact that
 in our proof we do not apply any result proved for the more general class of immigration superprocesses.

\begin{Pro}\label{Pro_Li_CB_joint_Laplace}
Let $(X_t)_{t\in\RR_+}$ be a CBI process with parameters $(1, c, \beta, B, \nu, \mu)$.
Let $\varphi$ and $\psi$ denote the branching and immigration mechanisms of $(X_t)_{t\in\RR_+}$,
respectively.
 \begin{itemize}
  \item[(i)] For any $0\leq r\leq t$ and $\lambda\in\RR_+$, we have that 
             \begin{align}\label{help_CBI_Laplace4}
                \EE(\ee^{-\lambda X_t} \mid X_r=x)
                    = \ee^{-x  g_t(r,\lambda) - \int_r^t \psi\left( g_t(u,\lambda)\right)\,\dd u}, \qquad x\in\RR_+,
             \end{align}
            where $[0,t]\ni s\mapsto g_t(s,\lambda):=v(t-s,\lambda)\in\RR_+$ is the unique locally bounded solution 
            to the equation
            \begin{align}\label{help_g_function_1}
              g_t(s,\lambda) + \int_s^t \varphi(g_t(u,\lambda))\,\dd u =\lambda,\qquad s\in[0,t].
            \end{align}
            Here, for any $\lambda\in\RR_+$, the continuously differentiable
            function $\RR_+ \ni t \mapsto v(t, \lambda)\in \RR_+$ is the unique locally bounded solution of
             \eqref{help12} with $d=1$.   
  \item[(ii)] For any $k\in\NN$, $0\leq r\leq t_1< t_2< \cdots< t_k$ and $\lambda_1,\ldots,\lambda_k\in\RR_+$,
              we have that
              \[
              \EE\left(\ee^{-\sum_{i=1}^k \lambda_i X_{t_i}}  \mid X_r=x\right)
                  = \ee^{-x g(r) - \int_r^{t_k} \psi(g(u))\,\dd u },\qquad r\in[0,t_1],\qquad x\in\RR_+,
              \]  
              where $[0,t_k]\ni s\mapsto g(s)$ is the unique locally bounded non-negative solution
              to the equation
              \begin{align}\label{help_g_function_2}
                g(s) + \int_s^{t_k} \varphi(g(u))\,\dd u
                 = \sum_{i=1}^k \lambda_i \bbone_{[0,t_i]}(s),\qquad s\in[0, t_k].
              \end{align}
 \end{itemize} 
\end{Pro}

\noindent{\bf Proof.}
(i). Using that $(X_t)_{t\in\RR_+}$ is a time-homogeneous Markov process, 
 by Theorem \ref{CBI_exists} with $d=1$, we have that, for any $r\in[0,t]$ and $\lambda\in\RR_+$,
 \begin{align}\label{help_CBI_Laplace3}
    \EE(\ee^{-\lambda X_t} \mid X_r=x)
                    = \ee^{-x v(t-r,\lambda) - \int_0^{t-r} \psi\left(v(u,\lambda)\right)\,\dd u}, \qquad x\in\RR_+,
 \end{align}
 where, for any \ $\lambda \in \RR_+$, \ the continuously differentiable function
 \ $\RR_+ \ni t \mapsto v(t, \lambda)\in\RR_+$ \ is the unique locally bounded solution 
  to \eqref{help12} with $d=1$.
Writing \eqref{help12} with $d=1$ in an integral form, for any \ $\lambda \in \RR_+$, \ we have that 
 \[
  v(t-r,\lambda) = \lambda - \int_0^{t-r} \varphi(v(u,\lambda))\,\dd u, \qquad r\in[0,t],
 \]
 and hence, using that $g_t(s,\lambda) = v(t-s,\lambda)$, $s\in[0,t]$, it implies that 
  \[
    g_t(s,\lambda) = \lambda - \int_0^{t-s} \varphi( g_t(t-u,\lambda))\,\dd u
                          = \lambda - \int_s^t \varphi( g_t(w,\lambda))\,\dd w, \qquad s\in[0,t],
  \]
  where we used the substitution $w=t-u$ at the last equality.
Therefore, the function $[0,t]\ni s\mapsto g_t(s,\lambda)$
 is a locally bounded non-negative solution to equation \eqref{help_g_function_1}.
The uniqueness of a locally bounded non-negative solution to the equation \eqref{help_g_function_1}
 follows by Step 3 of the proof of Theorem 4.2 in Li \cite{Li3} (in a more general setup).
Furthermore, using \eqref{help_CBI_Laplace3}, the definition of $g_t$ 
 and that 
 \[
    \int_0^{t-r} \psi\left(v(u,\lambda)\right)\,\dd u 
        = \int_0^{t-r} \psi\big(g_t(t-u,\lambda)\big)\,\dd u
        = \int_r^t \psi( g_t(w,\lambda))\,\dd w, 
        \qquad r\in[0,t],
 \] 
 we obtain \eqref{help_CBI_Laplace4}, as desired.

(ii).
We prove part (ii) by induction on $k$. 
For $k=1$, the statement is evident from part (i). Let $k\in\NN$, $0\leq r\leq t_1<t_2<\ldots<t_k<t_{k+1}$, and $\lambda_i\in\RR_+$, $i\in\{1,\ldots,k+1\}$. 
Using that $(X_t)_{t\in\RR_+}$ is a time-homogeneous Markov process, 
 by the tower rule for conditional expectations and part (i), we have that
  \begin{align*}
   \EE( \ee^{- \sum_{i=1}^{k+1} \lambda_iX_{t_i} }  \mid X_r)
    & = \EE\Big( \EE( \ee^{- \sum_{i=1}^{k+1} \lambda_iX_{t_i}} \mid X_r,X_{t_1},\ldots,X_{t_k}) \mid X_r\Big)\\
    & = \EE\big( \ee^{-\sum_{i=1}^k \lambda_iX_{t_i}} \EE( \ee^{- \lambda_{k+1}X_{t_{k+1}}} \mid X_r,X_{t_1},\ldots,X_{t_k}) \mid X_r\big) \\
    & = \EE\big( \ee^{-\sum_{i=1}^k \lambda_iX_{t_i} } \EE( \ee^{-  \lambda_{k+1}X_{t_{k+1}}} \mid X_{t_k} ) \mid X_r\big)\\
    &=\EE\big( \ee^{-\sum_{i=1}^k \lambda_iX_{t_i} }  \ee^{-X_{t_k}h(t_k,\lambda_{k+1})  - \int_{t_k}^{t_{k+1}} \psi(h(u,\lambda_{k+1}))\,\dd u }\mid X_r\big),
 \end{align*}
 where $[0,t_{k+1}]\ni s\mapsto h(s,\lambda_{k+1}):=v(t_{k+1}- s,\lambda_{k+1})\in\RR_+$ is 
 the unique locally bounded solution to the equation
 \begin{align}\label{help_g_function_11}
   h(s,\lambda_{k+1}) + \int_s^{t_{k+1}} \varphi(h(u,\lambda_{k+1}))\,\dd u =\lambda_{k+1},\qquad s\in[0,t_{k+1}].
 \end{align}
Hence, by the induction hypothesis, we have that 
\begin{align}\label{help_g_function_4}
\begin{split}
   \EE( \ee^{- \sum_{i=1}^{k+1} \lambda_iX_{t_i} }  \mid X_r)
   &= \ee^{- \int_{t_k}^{t_{k+1}} \psi(h(u,\lambda_{k+1}))\,\dd u}
            \EE\big( \ee^{-\sum_{i=1}^{k-1} \lambda_iX_{t_i} - (\lambda_k+h(t_k,\lambda_{k+1})) X_{t_k} }  \mid X_r\big) \\
   &= \ee^{ - \int_{t_k}^{t_{k+1}} \psi(h(u,\lambda_{k+1}))\,\dd u } \ee^{-X_r \widetilde g(r) - \int_r^{t_k} \psi(\widetilde g(u))\,\dd u },
\end{split}   
 \end{align}
 where $[0,t_k]\ni s\mapsto \widetilde g(s)$ is  the unique locally bounded non-negative solution
              to the equation
              \begin{align}\label{help_g_function_21}
                \widetilde g(s) + {\int_s^{t_k}} \varphi(\widetilde g(u))\,\dd u
                 = \sum_{i=1}^{k-1} \lambda_i{\bbone_{[0,t_i]}(s)}+(\lambda_k+h(t_k,\lambda_{k+1}))\bbone_{[0,t_k]}(s),\qquad s\in[0,{t_k}].
              \end{align}
Let $g:[0,t_{k+1}]\to\RR_+$ be defined by
 \[
         g(s):=\begin{cases}
                      \widetilde g(s) & \text{if $s\in[0,t_k]$,}\\
                       h(s,\lambda_{k+1}) & \text{if $s\in(t_k,t_{k+1}]$.}   
                     \end{cases} 
 \]
Then, in particular, we have $g(r) = \widetilde g(r)$, since $r\in[0,t_k]$. 
Therefore, by \eqref{help_g_function_4}, we get that
\begin{align*}
   \EE( \ee^{- \sum_{i=1}^{k+1} \lambda_iX_{t_i} }  \mid X_r)
    &= \ee^{ -X_r g(r) - \int_r^{t_{k+1}} \psi(g(u))\,\dd u}.
 \end{align*}
Next, we prove that the function $[0,{t_{k+1}}]\ni s\mapsto g(s)$ 
 is the unique locally bounded non-negative solution to the equation
 \begin{align}\label{help_DE_g}
   g(s) + {\int_s^{t_{k+1}}} \varphi(g(u))\,\dd u
                 = \sum_{i=1}^{k+1} \lambda_i{\bbone_{[0,t_i]}(s)},\qquad s\in[0,{t_{k+1}}].
 \end{align}
If $s\in[0,t_k]$, then, by splitting the integral at the time point $t_k$ into two parts,
 using the definition of $g$, and applying \eqref{help_g_function_21}, and \eqref{help_g_function_11} at the point $t_k$, 
 we have that 
 \begin{align*}
   g(s) + \int_s^{t_{k+1}} \varphi(g(u))\,\dd u
   & = g(s) + \int_s^{t_k} \varphi(g(u))\,\dd u
       + \int_{t_k}^{t_{k+1}} \varphi(  g(u) )\,\dd u \\
   & = \widetilde g(s) + \int_s^{t_k} \varphi(\widetilde g(u))\,\dd u
       + \int_{t_k}^{t_{k+1}} \varphi( h(u,\lambda_{k+1}) )\,\dd u \\    
   & = \sum_{i=1}^{k-1} \lambda_i{\bbone_{[0,t_i]}(s)}
       + (\lambda_k + h(t_k,\lambda_{k+1}) ) \bbone_{[0,t_k]}(s)
       + \lambda_{k+1} - h(t_k,\lambda_{k+1})\\
   &= \sum_{i=1}^{k-1} \lambda_i{\bbone_{[0,t_i]}(s)}  + \lambda_k + \lambda_{k+1}
    = \sum_{i=1}^{k+1} \lambda_i{\bbone_{[0,t_i]}(s)}.
 \end{align*} 
If $s\in(t_k,t_{k+1}]$, then, by the definition of $g$ and \eqref{help_g_function_11}, we have that 
 \begin{align*}
   g(s) + \int_s^{t_{k+1}} \varphi(g(u))\,\dd u
      = h(s,\lambda_{k+1}) + \int_{s}^{t_{k+1}} \varphi( h(u,\lambda_{k+1}) )\,\dd u 
      = \lambda_{k+1} = \sum_{i=1}^{k+1} \lambda_i{\bbone_{[0,t_i]}(s)}. 
  \end{align*}
This argument shows that $g$ is a solution of \eqref{help_DE_g}. 
Further, since $\widetilde g$ and $h$ are non-negative and locally bounded, we have that $g$ is non-negative and locally bounded as well. 
Therefore $g$ is indeed a locally bounded non-negative solution to the equation \eqref{help_DE_g}.
The uniqueness of a locally bounded non-negative solution to the equation \eqref{help_DE_g} 
 follows by Step 3 of the proof of Theorem 4.2 in Li \cite{Li3} (in a more general setup).
\proofend

In the next proposition, we relate the function $v$ 
 that appears in Proposition \ref{Pro_CBI_joint_Laplace} with $d=1$, $\beta=0$ and $\nu=0$ and the function $g$ defined in part (ii) of Proposition \ref{Pro_Li_CB_joint_Laplace}
 with $\beta=0$, $\nu=0$ and $r=0$. We note that the proof  does not use the theory of CB processes. 

\begin{Pro}\label{Relation_g_v}
Let $(X_t)_{t\in\RR_+}$ be a CB process with parameters $(1, c, 0, B, 0, \mu)$.
Let $\varphi$ denote the branching mechanism of $(X_t)_{t\in\RR_+}$.
For all $k\in\NN$, $0=:t_0<t_1<\ldots<t_k$
 and $\lambda_1,\ldots,\lambda_k\in\RR_+$, the unique locally bounded non-negative solution $g:[0, t_k]\rightarrow \mathbb{R}$ to the equation 
  \eqref{help_g_function_2} takes the form 
  \begin{align}\label{help_g_function_3}
    g(s) = v(t_1-s,y^{(k)}_1 ) \bbone_{[0,t_1]}(s) + \sum_{i=2}^k v(t_i-s,y^{(k)}_{i}) \bbone_{(t_{i-1},t_i]}(s),
         \qquad s\in[0,t_k],
 \end{align}
 where 
 \begin{itemize}    
  \item for any $\lambda \in \RR_+$, the continuously differentiable
       function $\RR_+ \ni t \mapsto v(t, \lambda) \in \RR_+$ is the unique locally bounded solution
       to \eqref{help12},
  \item $y^{(k)}_k=\lambda_k$, and $y^{(k)}_i$, $i\in\{1,\ldots,k-1\}$, are given by \eqref{help_y_recursion}.
 \end{itemize}
\end{Pro}

\noindent{\bf Proof.} 
We first observe that, as a consequence of \eqref{help_g_function_2}, for all $s,t\in[0,t_k]$ with $s\leq r$, it holds that
\begin{equation}
\label{eq:gsr}
g(s)-g(r)+\int_s^{r} \varphi(g(u))\,\dd u =\sum_{i=1}^k\lambda_i \bbone_{ \{s\leq t_i <r\}}. 
\end{equation}
 Let $k\in\NN$, $0=t_0<t_1<\ldots<t_k$ and $\lambda_1,\ldots,\lambda_k\in\RR_+$ be fixed.
Let us prove \eqref{help_g_function_3} by a kind of backward induction argument. First, we consider the case $s\in(t_{k-1},t_k]$.
Note that the equation \eqref{help_g_function_2} on the interval $(t_{k-1},t_k]$ takes the form
  \begin{equation}\label{eq:g_help_1}
    g(s) + \int_s^{t_k} \varphi(g(u))\,\dd u = \lambda_k ,\qquad s\in(t_{k-1},t_k].
  \end{equation} 
In particular, $g(t_k)=\lambda_k=y^{(k)}_k=v(0,y^{(k)}_k)$,
 where the last equality follows by \eqref{help12}.
Since, on the interval $(t_{k-1},t_k]$, the equation \eqref{eq:g_help_1} coincides with the 
 equation \eqref{help_g_function_1} with $t:=t_k$ and $\lambda:=\lambda_k$,  and 
 $\eqref{help_g_function_1}$ has a unique locally bounded non-negative solution 
 $[0,t_k]\ni s\mapsto g_{t_k}(s,\lambda_k)=v(t_k-s,\lambda_k)\in\RR_+$,
 we get that the two functions $g$ and $[0,t_k]\ni s\mapsto g_{t_k}(s,\lambda_k)$ must coincide on $(t_{k-1},t_k]$, i.e.,
 \begin{equation}\label{eq: g_help_2}
  g(s) =  g_{t_k}(s,\lambda_k) = v(t_k-s,\lambda_k) = v\big(t_k-s,y^{(k)}_k\big), \qquad s\in (t_{k-1},t_k].
 \end{equation}
 
Moreover, using first \eqref{eq:gsr} with $s=t_{k-1}$ and $r=t_k$ together with $g(t_k)=\lambda_k$, the equation 
 \eqref{eq: g_help_2}, and \eqref{help_g_function_1} with $s=t_{k-1}$, $t=t_k$ and $\lambda=\lambda_k$, we get that
\begin{align*}
g(t_{k-1})&=\lambda_{k-1}+\lambda_k-\int_{t_{k-1}}^{t_k} \varphi(g(u))\,\dd u =\lambda_{k-1}+\lambda_k-\int_{t_{k-1}}^{t_k} \varphi(g_{t_k}(u, \lambda_k))\,\dd u \\
&=\lambda_{k-1}+g_{t_k}(t_{k-1},\lambda_k)=\lambda_{k-1}+v(t_k-t_{k-1},\lambda_k)= y^{(k)}_{k-1}= v(0,y^{(k)}_{k-1}),
\end{align*}
where, at the fourth equality, we used the definition of $g_{t_k}$, 
at the fifth equality, the definition of $y^{(k)}_{k-1}$,
 and, at the last equality, \eqref{help12}. 
Therefore, \eqref{help_g_function_3} holds at the point $t_{k-1}$, which together with \eqref{eq: g_help_2} yield that 
 \eqref{help_g_function_3} holds for $s\in[t_{k-1},t_k]$.

Now suppose that \eqref{help_g_function_3} holds for any $s\in[t_m,t_k]$ for some $m\in\{1,\dots,k-1\}$. 
In particular, we have 
 \begin{equation} \label{eq:g(t_m)}
 g(t_m)=v(0,y^{(k)}_{m})=y^{(k)}_{m}. 
 \end{equation}
Next, we check that \eqref{help_g_function_3} holds for $s\in[t_{m-1},t_m)$ as well.
We do it in two steps, first we verify that it holds for $s\in(t_{m-1},t_m)$, and then for $s=t_{m-1}$ as well.
Using \eqref{eq:gsr} with $s\in(t_{m-1},t_m)$ and $r=t_m$  and \eqref{eq:g(t_m)}, we obtain that
\begin{equation}\label{eq:g_help_3}
 g(s)+\int_s^{t_m} \varphi(g(u))=g(t_m)= y^{(k)}_{m},  \qquad  s\in (t_{m-1},t_m].
\end{equation}
On the interval $(t_{m-1},t_m]$, the equation \eqref{eq:g_help_3} coincides with 
 equation \eqref{help_g_function_1} with $t:=t_{m}$ and $\lambda:=y^{(k)}_{m}$, 
 which has a unique locally bounded non-negative solution $[0,t_m]\ni s\mapsto g_{t_m}(s,y^{(k)}_{m}):=v(t_m- s,y^{(k)}_{m})\in\RR_+$. 
Therefore, it holds that
\begin{equation} \label{eq:h_help_4}
 g(s)= g_{t_m}(s,y^{(k)}_{m})=v(t_m-s,y^{(k)}_{m}), \qquad s\in (t_{m-1},t_m].
\end{equation}

Finally, let us prove that \eqref{help_g_function_3} holds for $s=t_{m-1}$ as well.
We define $\lambda_0:=0$.
Using \eqref{eq:gsr} with $s=t_{m-1}$ and $r=t_m$, together with \eqref{eq:g(t_m)}, \eqref{eq:h_help_4}
 and \eqref{help_g_function_1} with $s=t_{m-1}$, $t=t_m$ and $\lambda=y^{(k)}_{m}$, we obtain that
\begin{align*}
g(t_{m-1})&=\lambda_{m-1}+g(t_m)-\int_{t_{m-1}}^{t_m} \varphi(g(u))\,\dd u =\lambda_{m-1}
   + y^{(k)}_{m}-\int_{t_{m-1}}^{t_m} g_{t_m}(u,y^{(k)}_{m}) \,\dd u \\
&=\lambda_{m-1}+g_{t_m}(t_{m-1},y^{(k)}_{m})=\lambda_{m-1}+v(t_m-t_{m-1},y^{(k)}_{m}),
\end{align*}
where, at the last equality, we used the definition of $g_{t_m}$.
Therefore, using the definition of $y^{(k)}_{m-1}$ (see \eqref{help_y_recursion}) and \eqref{help12}, we get that
\[
 g(t_{m-1})= 
 \begin{cases}
  v(t_1, y^{(k)}_{1})& \text{if  $m=1$,}\\
  y^{(k)}_{m-1}=v(0, y^{(k)}_{m-1})& \text{if $m\in\{2,\ldots,k-1\}$.}  
 \end{cases}
\]
 Thus \eqref{help_g_function_3} holds for all $s\in[t_{m-1},t_m)$ as well.
Therefore, \eqref{help_g_function_3} holds for all $s\in[0,t_k]$, as desired.
\proofend

\begin{Rem}
(i).
The fact that the equation \eqref{help_g_function_2} has a unique locally bounded non-negative solution
 also follows from Step 3 of the proof of Theorem 4.2 in Li \cite{Li3}, 
 where the author, in a more general setup (replacing $\sum_{i=1}^k \lambda_i \bbone_{[0,t_i]}$
 by a more general function), proved this fact by contradiction
 using Gr\"onwall's inequality.
In Proposition \ref{Relation_g_v}, we derived an explicit formula 
 for the unique solution in question in our special case.
In the general setup, no explicit formula has been derived in Li \cite{Li3}.

(ii).
For any $t\in\RR_+$, the function $g_t$ given by \eqref{help_g_function_1} is continuous. 
However, when $\lambda_i>0$, $i\in\{1,\dots,k\}$, the function $g$ given by \eqref{help_g_function_2} 
 (for an explicit formula, see \eqref{help_g_function_3})
 is continuous only on $C:=[0,t_1)\cup(t_1,t_2)\cup\cdots\cup (t_{k-1},t_k]$. 
Indeed, by \eqref{help_g_function_3} and the fact that the function $\RR_+\ni t\mapsto v(t,\lambda)\in\RR_+$ 
 is continuous for any $\lambda\in\RR_+$, we can see that $g$ is continuous on $C$. 
Now consider $i\in\{1,\dots,k-1\}$ and a sequence $\{s_n,n\in\NN\}$ with $t_i<s_n<t_{i+1}$ and $s_n\downarrow t_i$ 
 as $n\rightarrow \infty$. 
Using \eqref{help_g_function_3} by choosing $s=t_i$ and $s=s_n,n\in\NN$, respectively, and \eqref{help12},
 we have that $g(t_i)=v(0,y^{(k)}_{i})=y^{(k)}_{i}$ and 
 \begin{align*}
   \lim_{n\rightarrow\infty} g(s_n)
      = \lim_{n\rightarrow\infty} v(t_{i+1}-s_n,y^{(k)}_{i+1})
      = v(t_{i+1}-t_i,y^{(k)}_{i+1})
      = y^{(k)}_i-\lambda_i<g(t_i),
 \end{align*}
 where the third equality follows by \eqref{help_y_recursion}.
Hence, if $\lambda_i>0$, $i\in\{1,\ldots,k\}$, then $g$ is not continuous at the points $t_i$, $i\in\{1,\dots,k-1\}$.
\proofend
\end{Rem}

In the next remark, we check that Proposition \ref{Pro_CBI_joint_Laplace}  with $d=1$ and 
 part (ii) of Proposition \ref{Pro_Li_CB_joint_Laplace}  with $r=0$ are in accordance.

\begin{Rem}\label{Rem_AppA_2}
Let $(X_t)_{t\in\RR_+}$ be a CBI process with parameters $(1, c, \beta, B, \nu, \mu)$.
Let $x\in\RR_+$, $k\in\NN$, $0=:t_{0}< t_1< t_{2}<\cdots< t_k$ 
 and $\lambda_1,\ldots,\lambda_k\in\RR_+$.
By a direct calculation, we check that the formulae for 
 $\EE\big(\ee^{-\sum_{i=1}^k \lambda_i X_{t_i}}  \mid X_0=x\big)$
 given by Proposition \ref{Pro_CBI_joint_Laplace} with $d=1$
 and by part (ii) of Proposition \ref{Pro_Li_CB_joint_Laplace} with $r=0$  indeed coincide.
On the one hand, Proposition \ref{Pro_CBI_joint_Laplace} with $d=1$ implies that 
 \begin{align} \label{help_formula_1_CBI_Laplace}
 \EE\left( \ee^{-\sum_{i=1}^k \lambda_i X_{t_i} }  \mid X_0=x \right) 
      = \exp\Bigg\{ - x v\big(t_1, y^{(k)}_1\big)
                        - \sum_{i=1}^k \int_0^{t_i-t_{i-1}} \psi\big( v(u,y^{(k)}_{i}) \big) \,\dd u
                        \Bigg\},
 \end{align}
 and, on the other hand, by part (ii) of Proposition \ref{Pro_Li_CB_joint_Laplace} with $r=0$, we have
 \begin{align} \label{help_formula_2_CBI_Laplace}
 \EE\left(\ee^{-\sum_{i=1}^k \lambda_i X_{t_i}}  \mid X_0=x\right)
                  = \ee^{-x g(0) - \int_0^{t_k} \psi(g(u))\,\dd u }.
\end{align}
In what follows, by a direct calculation, we verify that the right hand sides of \eqref{help_formula_1_CBI_Laplace}
 and \eqref{help_formula_2_CBI_Laplace} indeed coincide.
Recall that, by \eqref{help_g_function_3} with $s=0$, we have that
  $g(0) =  v(t_1,y^{(k)}_1)$.
Furthermore, for each $i\in\{1,\ldots,k\}$, by the substitution $u=t_i-s$, we have that 
 \[
    \int_0^{t_i-t_{i-1}} \psi\big( v(u,y^{(k)}_{i}) \big)\,\dd u 
      = \int_{t_{i-1}}^{t_i} \psi\big( v(t_i-s,y^{(k)}_{i}) \big)\,\dd s
      = \int_{t_{i-1}}^{t_i} \psi(g(s))\,\dd s, 
 \]
 where the last equality follows by \eqref{help_g_function_3}.
 By summing up these equations and taking into account that $t_{0}=0$, we get that
 \[
  \sum_{i=1}^k\int_0^{t_i-t_{i-1}} \psi\big( v(u,y^{(k)}_{i}) \big)\,\dd u   
      = \sum_{i=1}^k \int_{t_{i-1}}^{t_i} \psi(g(s))\,\dd s
      = \int_0^{t_k} \psi(g(s))\,\dd s,
 \]
 as desired.
\proofend
\end{Rem}

\section{Proof of Lemma \ref{Lem_vjumps}}\label{Appendix_Lemma_v}

This appendix is devoted to the proof Lemma \ref{Lem_vjumps}.

\smallskip

(i).
Let $x\in\RR_+$ and $(X_t)_{t\in\RR_+}$ be a CB process with parameters $(1, c, 0, B, 0,\mu)$ such that $X_0=x$.
By part (iv) of Corollary \ref{Cor_Laplace}, we have
\begin{equation}\label{eq:aux}
\EE\left( \ee^{-\lambda_0 X_t - \sum_{\ell=1}^k \lambda_\ell J_t(A_\ell)}  \mid X_0 = x\right)= \exp\left\{ -x v^{(k)}(t,\lambda_0,\lambda_1,\ldots,\lambda_k) 
               \right\}
\end{equation}
 for all $t,\lambda_\ell \in\RR_+$, $\ell\in\{0,1,\ldots,k\}$.
Let $i\in\{0,1,\dots,k\}$ and $t,\lambda_j \in\RR_+$, $j\in \{0,1,\dots,k\}\setminus\{i\}$ be fixed.
Then, using the monotonicity of expectation, the function $\RR_+\ni\lambda_i \mapsto \EE\left( \ee^{-\lambda_0 X_t - \sum_{\ell=1}^k \lambda_\ell J_t(A_\ell)}  \mid X_0 = x\right)$ is monotone decreasing. 
Hence, by \eqref{eq:aux}, the function $\RR_+\ni\lambda_i \mapsto v^{(k)}(t,\lambda_0,\lambda_1,\ldots,\lambda_k)$ 
 is monotone increasing, as desired.

(ii). The proof is divided into four steps.
First, we will show that there exists a unique locally bounded non-negative solution to 
 the differential equation \eqref{eq_vAlambda}. 
Then we will prove that $v^{(k)}(t,\eta_s,s\lambda_1,\dots,s\lambda_k)$ converges to a finite limit
 in a monotone increasing way as $s\uparrow \infty$ for all $t,\lambda_0,\lambda_1,\ldots,\lambda_k\in\RR_+$ 
 and $\eta_s\uparrow\lambda_0$ as $s \uparrow \infty$.
Finally, we will verify that the limit in question is $v^{(A_\blambda)}(t,\lambda_0)$.
 and it does not depend on the choice of the function $\RR_+\ni s\mapsto \eta_s$ 
 (only on its limit $\lambda_0$ as $s\uparrow \infty$).

{\sl Step 1.} Let $\blambda=(\lambda_1,\ldots,\lambda_k)\in\RR_+^k$ be arbitrarily fixed.
Let $(X^{(A_\blambda)}_t)_{t\in\RR_+}$ be a CB process
 with parameters \ $(1, c,0, B^{(A_\blambda)}, 0, \mu^{(A_\blambda^{-1})})$,
 where $A_\blambda$ is defined in \eqref{Alambda}, and
  $B^{(A_\blambda)}$ and $\mu^{(A_\blambda^{-1})}$ are given by \eqref{mod_parameters} with $d=1$ and $A=A_\blambda$. 
By Remark \ref{rem:VarphiA}, the branching mechanism of $(X^{(A_\blambda)}_t)_{t\in\RR_+}$ is
 $\varphi^{(A_\blambda)}$ given by \eqref{eq:psiA} with $A=A_\blambda$.
Then, by Corollary 5.17 and Theorem 5.18 in Li \cite{Li}, we have
\begin{equation*}
   \EE\left( \ee^{-\lambda_0 X^{(A_\blambda)}_t - \mu(A_\blambda)\int_{0}^{t}X^{(A_\blambda)}_s\dd s}  \mid X^{(A_\blambda)}_0 = x\right)=\exp\Big\{ -x v^{(A_\blambda)}(t,\lambda_0) \Big\},
      \qquad t,\lambda_0,x\in\RR_+,
\end{equation*}
 where, {for all $\lambda_0\in\RR_+$, the continuously differentiable function} $\RR_+\ni t \mapsto v^{(A_\blambda)}(t,\lambda_0){\in\RR_+}$
 is the unique locally bounded non-negative solution to the differential equation
 \begin{align*}%\label{help_DE_uj_b1}
 \partial_1 v^{(A_\blambda)}(t,\lambda_0) =\mu(A_\blambda)- \varphi^{(A_\blambda)}(v^{(A_\blambda)}(t,\lambda_0)),\qquad  t\in\RR_+,\qquad  v^{(A_\blambda)}(0,\lambda_0) = \lambda_0.
 \end{align*}
This differential equation is nothing else but the differential equation \eqref{eq_vAlambda},
 and therefore the differential equation \eqref{eq_vAlambda} has a unique locally bounded non-negative solution.

{\sl Step 2.}
Let $\lambda_0\in\RR_+$ and $\blambda=(\lambda_1,\ldots,\lambda_k)\in\RR_+^k$ be arbitrarily fixed.
Let $\RR_+\ni s\mapsto \eta_s$ be a function such that $\eta_s\uparrow\lambda_0$ as $s \uparrow \infty$.
By part (i), for all $t\in\RR_+$, the function $\RR_+\ni s \mapsto v^{(k)}(t,\eta_s,s\lambda_1,\dots,s\lambda_k)$ is monotone increasing. 
Indeed, if $s\leq s^*$, $s,s^*\in\RR_+$, then, by part (i), for all $t\in\RR_+$, we have that 
 \begin{align*}
  & v^{(k)}(t,\eta_s,s\lambda_1,\dots,s\lambda_k) 
    \leq v^{(k)}(t, \eta_{s^*},s\lambda_1,\dots,s\lambda_k)
     \leq v^{(k)}(t,\eta_{s^*},s^*\lambda_1,s\lambda_2,\dots,s\lambda_k)\\
  & \leq v^{(k)}(t,\eta_{s^*},s^*\lambda_1,s^*\lambda_2,s\lambda_3,\dots,s\lambda_k)\leq \cdots 
     \leq v^{(k)}(t,\eta_{s^*},s^*\lambda_1,s^*\lambda_2,s^*\lambda_3,\dots,s^*\lambda_k),
 \end{align*}
 as desired.
Then, for all $t\in\RR_+$, the limit
 \begin{align}\label{help_def_v}
    \tv^{(A_{\blambda})}(t,\lambda_0) :=  \lim_{s\uparrow\infty} v^{(k)}(t,\eta_s,s\lambda_1,\dots,s\lambda_k)
 \end{align}
 exists in $[0,\infty]$.
We call the attention that the limit $\tv^{(A_{\blambda})}(t,\lambda_0)$ in \eqref{help_def_v}, in principle,
 could depend on the function $\RR_+\ni s \mapsto \eta_s$ (not only on its limit $\lambda_0$ as $s\uparrow \infty$),
 but we do not denote this possible dependence, since, at the end of Step 4,
 we prove that it is not the case.
Observe that, by \eqref{def_v_k+1_CBI}, we have
$\tv^{(A_{\blambda})}(0,\lambda_0)=\lim_{ s\uparrow\infty} \eta_s = \lambda_0$.

{\sl Step 3.} We prove that $ \tv^{(A_{\blambda})}(t,\lambda_0) <\infty$ 
 for all $t,\lambda_0\in\RR_+$, $\eta_s\uparrow \lambda_0$ as $s\uparrow \infty$, and $\blambda = (\lambda_1,\ldots,\lambda_k)\in\RR_+^k$.
Let $\lambda_0\in\RR_+$, $\eta_s\uparrow \lambda_0$ as $s\uparrow \infty$, and $\blambda=(\lambda_1,\ldots,\lambda_k)\in\RR_+^k$ be arbitrarily fixed.
Recall that the function $\varphi_1^{(k)}:\RR_+^{k+1}\to\RR$ is given by \eqref{2CBI_branching}.
By \eqref{2CBI_branching}, the dominated convergence theorem and \eqref{eq:psiA} with $A=A_\blambda$, we have that 
\begin{align}\label{eq:banchinglimit}
\begin{split}
 \widetilde{\varphi}^{(\blambda)}(\lambda_0)
:&=\lim_{s\to\infty }\varphi_1^{(k)}(\lambda_0,s\lambda_1,\dots,s\lambda_k)\\
     & = c\lambda_0^2 - B\lambda_0
         + \lim_{s\to\infty }
           \int_{A_\blambda} \bigl(  \ee^{- \lambda_0 z -  \sum_{i=1}^k s\lambda_i\bbone_{A_{i}}(z)} - 1
                             + \lambda_0 (1 \land z) \bigr) \, \mu(\dd z) \\
     &\phantom{=\;} + \int_{\cU_1\setminus A_\blambda} \bigl( \ee^{- \lambda_0 z } - 1 + \lambda_0 (1 \land z) \bigr) \, \mu(\dd z)\\
     & = c\lambda_0^2 - \left(B - \int_{A_\blambda} (1\wedge z)\,\mu(\dd z)\right)\lambda_0 \\
     &\phantom{=\;}
        + \int_{\cU_1\setminus A_\blambda} \bigl( \ee^{- \lambda_0 z } - 1 + \lambda_0 (1 \land z) \bigr) \, \mu(\dd z)
         -\mu(A_\blambda)\\
     & = \varphi^{(A_\blambda)}(\lambda_0) - \mu(A_\blambda),
 \end{split}
 \end{align}
 where, at the last but one equality, we can indeed interchange the limit as $s\to\infty$
 and the integration on the set $A_\blambda$, since we can use the dominated convergence theorem 
 due to the facts that  
 \[
   \Big\vert  \ee^{- \lambda_0 z -  \sum_{i=1}^k s\lambda_i \bbone_{A_{i}}(z)} - 1 + \lambda_0(1\wedge z)\Big\vert
      \leq 2 + \lambda_0(1\wedge z) \leq 2 + \lambda_0,
      \qquad  s,z\in\RR_+,
 \]
 and
 \[
  \int_{A_\blambda} (2+\lambda_0)\,\mu(\dd z) = (2+\lambda_0) \mu(A_\blambda)<\infty.
 \]

Additionally, we check that
 \begin{align}\label{eq:help_comparison_varphi_1}
\varphi_1^{(k)}(\lambda_0,s\lambda_1,\dots,s\lambda_k)\geq\widetilde{\varphi}^{(\blambda)}(\lambda_0), \qquad s\in\RR_+.
 \end{align}
Using \eqref{2CBI_branching} and \eqref{eq:banchinglimit}, we have that \eqref{eq:help_comparison_varphi_1} is equivalent to
 \begin{align*}
  &\int_{\cU_1} \bigl( \ee^{- \lambda_0 z -  \sum_{i=1}^k s\lambda_i\bbone_{A_{i}}(z)} - 1
 + \lambda_0 (1 \land z) \bigr) \, \mu(\dd z)\\
  &\qquad \geq \int_{A_\blambda} (-1 + \lambda_0(1 \land z)) \, \mu(\dd z) 
 + \int_{\cU_1\setminus A_\blambda} \bigl( \ee^{- \lambda_0 z} - 1 + \lambda_0 (1 \land z) \bigr) \, \mu(\dd z),
  \qquad s\in\RR_+.
 \end{align*}
By splitting the domain of integration of the integral on the left-hand side of the inequality above
 into two parts, namely, by integrating on $A_\blambda$ and $\cU_1\setminus A_\blambda$,
 and, using that $\mu(A_\blambda)<\infty$, we have that \eqref{eq:help_comparison_varphi_1} is equivalent to
 \[
  \int_{A_\blambda} \ee^{- \lambda_0 z -  \sum_{i=1}^k s\lambda_i\bbone_{A_{i}}(z)}\,\mu(\dd z)\geq 0,
 \]
 which holds trivially.

For convenience in the subsequent derivations,  using \eqref{eq:banchinglimit}, for all $\lambda_0\in\RR_+$,
 we can express the differential equation \eqref{eq_vAlambda}
 in the following form:
 \begin{align}\label{eq:help_DECI}
  \partial_1 v^{(A_\blambda)}(t,\lambda_0) =-\widetilde{\varphi}^{(\blambda)}(v^{(A_\blambda)}(t,\lambda_0)),\qquad
 t\in\RR_+, \qquad  v^{(A_\blambda)}(0,\lambda_0) = \lambda_0.
 \end{align}
Recall also that, by \eqref{def_v_k+1_CBI}, for all $\lambda_0,s\in\RR_+$ and $\blambda\in\RR_+^k$, we have that
 \begin{align}\label{eq:help_DEVI}
\begin{split}
&\partial_1 v^{(k)}(t,\lambda_0,s\lambda_1,\ldots,s\lambda_k)
=  - \varphi_1^{(k)}\big( v^{(k)}(t,\lambda_0,s\lambda_1,\ldots,s\lambda_k),s\lambda_1,\ldots,s\lambda_k\big),\\
&v^{(k)}(0, \lambda_0,s\lambda_1,\ldots,s\lambda_k) = \lambda_0.
\end{split}
\end{align} 

We will compare the solutions of the differential equations \eqref{eq:help_DECI} and \eqref{eq:help_DEVI}
 by using the following version of comparison theorem for differential equations (see, e.g., Volkmann \cite{Vol}
 or Lemma C.3 in Filipovi\'{c} et al.\ \cite{FilMaySch}): if \ $S:\RR_+\times\RR\to\RR$ \ is a continuous function
 which is locally Lipschitz continuous in its second variable and \ $p,q:\RR_+\to\RR$ \ are differentiable functions
 satisfying
 \begin{align*}
   & p'(t) \leq S(t,p(t)),\quad t\in\RR_+,\\
   & q'(t) = S(t,q(t)),\quad t\in\RR_+,\\
   & p(0)\leq q(0),
 \end{align*}
 then \ $p(t)\leq q(t)$ \ for all \ $t\in\RR_+$.
For all $\lambda_0,s\in\RR_+$ and $\blambda=(\lambda_1,\ldots,\lambda_k)\in\RR_+^k$, 
 let us choose $p(t):=v^{(k)}(t,\lambda_0,s\lambda_1,\ldots,s\lambda_k)$, $t\in\RR_+$, $q(t):=v^{(A_\blambda)}(t,\lambda_0)$, $t\in\RR_+$,
 and 
\[
  S(t,x):=- \widetilde{\varphi}^{(\blambda)}(x)\bbone_{\{x\geq 0\}}
           + \mu(A_\blambda)\bbone_{\{x< 0\}}, \qquad t\in\RR_+,x\in\RR,
\]
where $\widetilde{\varphi}^{(\blambda)}$ is given by \eqref{eq:banchinglimit}.
Then $S$ is continuous (since $\widetilde{\varphi}^{(\blambda)}(0) = \varphi^{(A_\blambda)}(0)-\mu(A_\blambda)
   = \varphi(0) - \mu(A_\blambda) = -\mu(A_\blambda)$),
 and we check that it is locally Lipschitz continuous in its second variable
 (because of the special form of $S$, in fact, it also implies the continuity of $S$).
Since $\int_{\cU_1} (z\wedge z^2)\,\mu(\dd z)<\infty$ (see \eqref{help_moment_mu})
 and $\int_{\cU_1} (1\wedge z^2)\,\mu(\dd z)<\infty$ (see the paragraph after Definition \ref{Def_CBI}),
 by Proposition 1.48 and the paragraph at the end of page 26 in Li \cite{Li}, we have that
 the function
 \begin{align}\label{eq:help_loc_Lip11}
  \RR_+\ni x\mapsto c x^2 - \left(B - \int_{A_\blambda} (1\wedge z)\,\mu(\dd z)\right)x 
      + \int_{\cU_1\setminus A_\blambda} \Bigl( \ee^{- x z } - 1 +  xz\bone_{(0,1]}(z) \Bigr) \, \mu(\dd z)
  \end{align}
 is locally Lipschitz continuous at any point of $\RR_+$.
Note that
 \begin{align}\label{eq:help_loc_Lip21}
 \begin{split}
 &\int_{\cU_1\setminus A_\blambda} \bigl( \ee^{- x z } - 1 + x (1 \land z) \bigr) \, \mu(\dd z)\\
 &\qquad  = \int_{\cU_1\setminus A_\blambda} \Bigl( \ee^{- x z } - 1 +  xz\bone_{(0,1]}(z)  \Bigr) \, \mu(\dd z)
       - x\int_{\cU_1\setminus A_\blambda}\left( z\bone_{(0,1]}  - (1\wedge z)\right)\mu(\dd z),
 \end{split}      
 \end{align}
 where
\begin{align*}
  \int_{\cU_1\setminus A_\blambda}\left( z\bone_{(0,1]}  - (1\wedge z)\right)\mu(\dd z)
  & = \int_{(\cU_1\setminus A_\blambda)\cap(0,1]}(z-z)\,\mu(\dd z)
      + \int_{(\cU_1\setminus A_\blambda)\cap(1,\infty)}(-1)\,\mu(\dd z)\\
  & = -\mu\big( (\cU_1\setminus A_\blambda)\cap(1,\infty)\big),
 \end{align*}
 which is finite.
Indeed, by \eqref{help_moment_mu}, we get that
 \[
   0\leq \mu\big( (\cU_1\setminus A_\blambda)\cap(1,\infty)\big)
    \leq \mu( (1,\infty) )
    =  \int_{(1,\infty)} 1\,\mu(\dd z)
    \leq \int_{(1,\infty)} (z\wedge z^2) \,\mu(\dd z)
    <\infty.
 \]
Consequently, taking into account \eqref{eq:help_loc_Lip21} and that the function given in \eqref{eq:help_loc_Lip11} 
 is locally Lipschitz continuous at any point of $\RR_+$, we have that, for all $t\in\RR_+$, the function $\RR_+\ni x \mapsto S(t,x)$ is locally Lipschitz continuous
  at any point of $\RR_+$.
Furthermore, since $S(t,x) = \mu(A_\blambda)$ for $t\in\RR_+$ and $x<0$, 
 we readily have that,  for all $t\in\RR_+$, the function $(-\infty,0)\ni x \mapsto S(t,x)$ is locally Lipschitz continuous 
 at any point of $(-\infty,0)$.
It remains to verify that the function $\RR\ni x \mapsto S(t,x)$ is locally Lipschitz continuous at $0$.
This follows from the fact that, for all $t\in\RR_+$, it holds that
 \[
   \vert S(t,x) - S(t,0)\vert
     =\begin{cases}
         \vert S(t,0) - S(t,0)\vert = 0 & \text{for $x<0$,} \\
         \leq K x = K\vert x\vert & \text{for $x\in [0,\vare)$,}
       \end{cases}
 \]
 with some appropriate $K>0$ and $\vare>0$ (the existence of such $K$ and $\vare$ follows from the fact that $\RR_+\ni x \mapsto S(t,x)$ 
 is locally Lipschitz continuous at any point of $\RR_+$).
Consequently, using also that local Lipschitz continuity at every point of $\RR$ implies local
 Lipschitz continuity on compact sets of $\RR$, we obtain that
 $S$ is locally Lipschitz continuous in its second variable.

By \eqref{eq:help_DECI}, \eqref{eq:help_DEVI}, and the non-negativity of $v^{(A_\blambda)}(t,\lambda_0)$, $t\in\RR_+$,
 we have that $p(0)=q(0) = \lambda_0$ and
 \[
    q'(t) = S(t,q(t)) = - \widetilde{\varphi}^{(\blambda)}(q(t)), \qquad t\in\RR_+.
 \]
By \eqref{eq:help_comparison_varphi_1} and \eqref{eq:help_DEVI}, we get that
 \begin{align*}
   p'(t) &= \partial_1 v^{(k)}(t,\lambda_0,s\lambda_1,\ldots,s\lambda_k)
=  - \varphi_1^{(k)}\big( v^{(k)}(t,\lambda_0,s\lambda_1,\ldots,s\lambda_k),s\lambda_1,\ldots,s\lambda_k\big)\\
         &\leq  - \widetilde{\varphi}^{(\blambda)}( v^{(k)}(t,\lambda_0,s\lambda_1,\ldots,s\lambda_k)\big))
         = S(t,p(t)), \qquad t\in\RR_+.
 \end{align*}
Consequently, one can indeed apply the above recalled comparison theorem in order to obtain that
 \begin{align}\label{eq:help_comparison_vi_vci}
 v^{(k)}(t,\lambda_0,s\lambda_1,\ldots,s\lambda_k)\leq v^{(A_\blambda)}(t,\lambda_0), \qquad  t,s\in\RR_+.
 \end{align}

Since $\eta_s\uparrow\lambda_0$ as $s \uparrow \infty$, by part (i) and \eqref{eq:help_comparison_vi_vci}, we obtain that 
\begin{align}
\label{eq:bounded_s}
v^{(k)}(t,\eta_s,s\lambda_1,\ldots,s\lambda_k)\leq v^{(k)}(t,\lambda_0,s\lambda_1,\ldots,s\lambda_k)\leq v^{(A_{\blambda})}(t,\lambda_0),
  \qquad t,s\in\RR_+.
\end{align}
Therefore, by taking the limit as  $s\uparrow\infty$ and  using \eqref{help_def_v}, we have
\[
  \tv^{(A_{\blambda})}(t,\lambda_0) =  \lim_{ s\uparrow\infty} v^{(k)}(t, \eta_s,s\lambda_1,\dots,s\lambda_k)
       \leq v^{(A_\blambda)}(t,\lambda_0)
\]
 for all $t,\lambda_0\in\RR_+$, $\eta_s\uparrow \lambda_0$ as $s\uparrow \infty$ and $\blambda=(\lambda_1,\ldots,\lambda_k)\in\RR_+^k$.
Since $v^{(A_\blambda)}(t,\lambda_0)<\infty$ for all $t,\lambda_0\in\RR_+$ and $\blambda\in\RR_+^k$,
 we get that $ \tv^{(A_{\blambda})}(t,\lambda_0)<\infty$ for all $t,\lambda_0\in\RR_+$, $\eta_s\uparrow \lambda_0$ as $s\uparrow \infty$
 and $\blambda\in\RR_+^k$, as desired.

{\sl Step 4.}  Finally, for all $\lambda_0\in\RR_+$, $\eta_s\uparrow \lambda_0$ as $s\uparrow \infty$ and $\blambda=(\lambda_1,\ldots,\lambda_k)\in\RR_+^k$,
 we prove that the function $\RR_+\ni t\mapsto\tv^{(A_{\blambda})}(t,\lambda_0)\in\RR_+$ (introduced in \eqref{help_def_v})
 is a locally bounded non-negative solution to the differential equation \eqref{eq_vAlambda}. 
Note that, by Steps 2 and 3, this function indeed takes values in $\RR_+$, i.e., 
 $\tv^{(A_{\blambda})}(t,\lambda_0)$ is non-negative and finite.
Let $\lambda_0\in\RR_+$, $\eta_s\uparrow \lambda_0$ as $s\uparrow \infty$ and $\blambda=(\lambda_1,\ldots,\lambda_k)\in\RR_+^k$ be arbitrarily fixed.
By \eqref{eq:help_DEVI}, we get that, for all $ s,t\in\RR_+$,
  \begin{align}\label{eq:DE_2_integral_form_vis}
   v^{(k)}(t,\eta_s,s\lambda_1,\ldots,s\lambda_k)
     = \eta_s - \int_0^t  \varphi_1^{(k)}\big( v^{(k)}(r,\eta_s,s\lambda_1,\ldots,s\lambda_k),s\lambda_1,\ldots,s\lambda_k\big)\,\dd r.
 \end{align}
Recall that $ \tv^{(A_{\blambda})}(0,\lambda_0) =\lambda_0$ (see Step 2).
 In what follows, using \eqref{help_def_v}, \eqref{eq:DE_2_integral_form_vis}, the dominated convergence theorem 
 and \eqref{eq:banchinglimit}, we verify that
 \begin{align}\label{eq:help210_limit_lambda_0}
 \begin{split}
  \tv^{(A_{\blambda})}(t,\lambda_0)
   & =\lambda_0 - \lim_{ s\uparrow\infty} \int_0^t \varphi_1^{(k)}\big( v^{(k)}(r, \eta_s,s\lambda_1,\ldots,s\lambda_k),s\lambda_1,\ldots,s\lambda_k\big)\,\dd r\\
   & =\lambda_0+\mu(A_{\blambda})t - \int_0^t \varphi^{(A_{\blambda})}( \tv^{(A_{\blambda})}(r,\lambda_0)) \,\dd r,\qquad t\in\RR_+,
  \end{split}
 \end{align}
where the function $\varphi^{(A_\blambda)}:\RR_+\to\RR$ is given by \eqref{eq:psiA} with $A =A_\blambda$. 
The first equality in \eqref{eq:help210_limit_lambda_0} follows from  \eqref{help_def_v} and \eqref{eq:DE_2_integral_form_vis}.
Now, we turn to prove the second equality in \eqref{eq:help210_limit_lambda_0}.
Using \eqref{2CBI_branching} we have the decomposition
 \begin{align}\label{help_AppB_1}
 \begin{split} 
  &\int_0^t \varphi_1^{(k)}\big( v^{(k)}(r, \eta_s,s\lambda_1,\ldots,s\lambda_k),s\lambda_1,\ldots,s\lambda_k\big)\,\dd r\\
  &\qquad  = S_1(s)+S_2(s)+S_3(s)+S_4(s),\qquad s,t\in\RR_+,
  \end{split}
 \end{align}
 where
	\begin{align*}
      S_1(s)&:=\int_{0}^{t}-Bv^{(k)}(r,\eta_s,s\lambda_1,\ldots,s\lambda_k)\,\dd r,\\
		S_2(s)&:=\int_{0}^{t}c (v^{(k)}(r,\eta_s,s\lambda_1,\ldots,s\lambda_k))^2\,\dd r,\\
		S_3(s)&:=\int_{0}^{t} \int_{A_{\blambda}} \Big( {\ee}^{-v^{(k)}(r,\eta_s,s\blambda)z  -  \sum_{i=1}^k s\lambda_i\bbone_{A_{i}}(z)} - 1 +v^{(k)}(r,\eta_s,s\blambda)(1 \land z)\Big)\mu(\dd z)\,\dd r,\\
	S_4(s)&:=\int_{0}^{t} \int_{\cU_1\setminus A_{\blambda}} \Big( {\ee}^{-v^{(k)}(r,\eta_s,s\blambda)z } - 1 + v^{(k)}(r,\eta_s,s\blambda)(1 \land z)\Big)\mu(\dd z)\,\dd r.
    \end{align*}
Here $S_i(s)$, $i\in\{1,2,3,4\}$, may depend on $t$, $(\eta_s)_{s\in\RR_+}$ (in particular, its limit $\lambda_0$ as $s\to\infty$)
 and $\blambda$ as well, but we do not denote these dependence,
 since we investigate the limit as $s$ tends to $\infty$
 (so we display the dependence only on $s$).
First, we check that, for all $t\in\RR_+$, it holds that
 \begin{align}\label{eq:help_S1_1}
     \lim_{s\uparrow\infty} S_1(s) = \int_{0}^{t} \lim_{s\uparrow \infty}\big(-Bv^{(k)}(r,\eta_s,s\lambda_1,\ldots,s\lambda_k)\big)\,\dd r
=\int_{0}^{t} 
- B \tv^{(A_{\blambda})}(r,\lambda_0) \,\dd r.
 \end{align} 
By \eqref{eq:bounded_s} and the nonnegativity of $v^{(k)}(r,\eta_s,s\lambda_1,\ldots,s\lambda_k)$ and $v^{(A_{\blambda})}(s,\lambda_0)$,
 we get that
\begin{equation*}\label{eq:S1}
  |-B v^{(k)}(r,\eta_s,s\lambda_1,\ldots,s\lambda_k)|\leq|B|v^{(A_{\blambda})}(r,\lambda_0), \qquad r,s\in\RR_+.
\end{equation*}
Since the function $\RR_+\ni r \mapsto v^{(A_{\blambda})}(r,\lambda_0)$ is locally bounded,
 we have that $[0,t]\ni r\mapsto v^{(A_{\blambda})}(r,\lambda_0)$ is bounded for any $t\in\RR_+$, and therefore
 $$
 \int_{0}^{t}|B|v^{(A_{\blambda})}(r,\lambda_0)\dd r<\infty, \qquad t\in\RR_+.
 $$
Therefore, one can apply the dominated convergence theorem and \eqref{help_def_v} to obtain \eqref{eq:help_S1_1}.

Now, we check that, for all $t\in\RR_+$, it holds that
 \begin{align}\label{eq:help_S2_1}
     \lim_{s\uparrow\infty} S_2(s) = \int_{0}^{t} \lim_{s\uparrow\infty} \big(c (v^{(k)}(r,\eta_s,s\lambda_1,\ldots,s\lambda_k))^2\big) \,\dd r= \int_{0}^{t} 
c ( \tv^{(A_{\blambda})}(r,\lambda_0) )^2    \,\dd r  .
 \end{align} 
Similarly as in case of $S_1(s)$, we get that
  $$
     |c (v^{(k)}(r,\eta_s,s\lambda_1,\ldots,s\lambda_k))^2|\leq c (v^{(A_{\blambda})}(r,\lambda_0))^2, \qquad r,s\in\RR_+.
  $$
 where
\begin{equation*}\label{eq:S2}
   \int_{0}^{t}c (v^{(A_{\blambda})}(r,\lambda_0))^2\dd r<\infty, \qquad t\in\RR_+.
\end{equation*}
Therefore, one can apply the dominated convergence theorem and \eqref{help_def_v}
 to obtain \eqref{eq:help_S2_1}. 

Now, we check that, for all $t\in\RR_+$, it holds that
 \begin{align}\label{eq:help_S3_1}
 \begin{split}
&\lim_{s\uparrow \infty} S_3(s) \\
&= \int_{0}^{t} \int_{A_{\blambda}} 
\lim_{s\uparrow \infty} \Big( {\ee}^{-v^{(k)}(r,\eta_s,s\blambda)z  -  \sum_{i=1}^k s\lambda_i\bbone_{A_{i}}(z) } - 1 +v^{(k)}(r,\eta_s,s\blambda)(1 \land z)\Big)\mu(\dd z)\dd r\\
&= -t\mu(A_{\blambda})+\int_{0}^{t} \int_{A_{\blambda}} 
 \tv^{(A_{\blambda})}(r,\lambda_0)(1 \land z)\mu(\dd z)\dd r.
 \end{split}                      
 \end{align}
Using again \eqref{eq:bounded_s}, the non-negativity of $v^{(k)}(r,\eta_s,s\blambda)$ and $\tv^{(A_{\blambda})}(r,\lambda_0)$
 and that $1\wedge z\in(0,1]$, $z\in\cU_1$, for all $z\in A_{\blambda}$, we have
 $$
  | {\ee}^{-v^{(k)}(r,\eta_s,s\blambda)z  -  \sum_{i=1}^k s\lambda_i\bbone_{A_{i}}(z) } - 1 +v^{(k)}(r,\eta_s,s\blambda)(1 \land z)|\leq  2 + v^{(A_{\blambda})}(r,\lambda_0).
 $$
Since $\RR_+\ni r \mapsto v^{(A_{\blambda})}(r,\lambda_0)$ is locally bounded and $\mu(A_{\blambda})<\infty$, we have that 
\begin{equation*}\label{eq:S3}
  \int_{0}^{t} \int_{A_{\blambda}} \left( 2 + v^{(A_{\blambda})}(r,\lambda_0)\right)\mu(\dd z)\dd r
    = 2\mu(A_{\blambda})t + \mu(A_{\blambda}) \int_{0}^{t} v^{(A_{\blambda})}(r,\lambda_0) \,\dd r<\infty,
    \qquad t\in\RR_+.
 \end{equation*}
Note also that $ \sum_{i=1}^k s\lambda_i\bbone_{A_{i}}(z) \to \infty$ as $s\to\infty$ for all $z\in A_{\blambda}$.
Therefore, one can apply the dominated convergence theorem and \eqref{help_def_v} to obtain \eqref{eq:help_S3_1}.

Next, we check that, for all $t\in\RR_+$, it holds that
\begin{align}\label{eq:help_S4_1}
 \begin{split}
\lim_{s\uparrow \infty} S_4(s) &= \int_{0}^{t} \int_{\cU_1\setminus A_{\blambda}}  \lim_{s\uparrow\infty} 
\Big( {\ee}^{-v^{(k)}(r,\eta_s,s\blambda)z } - 1 + v^{(k)}(r,\eta_s,s\blambda)(1 \land z)\Big)\mu(\dd z)\dd r\\
&=\int_0^t\int_{\cU_1\setminus A_{\blambda}} 
\Big( {\ee}^{- \tv^{(A_{\blambda})}(r,\lambda_0)z } - 1 +  \tv^{(A_{\blambda})}(r,\lambda_0)(1 \land z)\Big)\mu(\dd z)\dd r.
 \end{split}
 \end{align}
For $z\in (\cU_1\setminus A_{\blambda})\cap (0,1]$,  using \eqref{eq:bounded_s},
 the inequality $0\leq\ee^{-x} - 1 + x\leq \frac{1}{2}x^2$, $x\in\RR_+$, 
 and that $v^{(k)}(r,\eta_s,s\blambda)$ is non-negative, we have
 \[
 |{\ee}^{-v^{(k)}(r,\eta_s,s\blambda)z } - 1 + v^{(k)}(r,\eta_s,s\blambda)z|
 \leq \frac{1}{2} (v^{(k)}(r,\eta_s,s\blambda))^2 z^2 \leq  (v^{(A_{\blambda})}(r,\lambda_0))^2z^2.
 \]
For $z\in (\cU_1\setminus A_{\blambda})\cap (1,\infty)$, using \eqref{eq:bounded_s}, the inequality $0\leq 1-\ee^{-x}\leq x$, $x\in\RR_+$,
 and that $v^{(k)}(r,\eta_s,s\blambda)$ is non-negative, we have
\begin{align*}
|{\ee}^{-v^{(k)}(r,\eta_s,s\blambda)z } - 1 + v^{(k)}(r,\eta_s,s\blambda)|\leq 2 v^{(k)}(r,\eta_s,s\blambda)z
\leq  2 v^{(A_{\blambda})}(r,\lambda_0)z.
\end{align*}
Since $\RR_+\ni r \mapsto v^{(A_{\blambda})}(r,\lambda_0)$ is locally bounded and \eqref{help_moment_mu} holds, 
 we obtain that
 \begin{equation*}\label{eq:S4}
\begin{split}
 &\int_{0}^{t} \int_{(\cU_1\setminus A_{\blambda})\cap (0,1]}(v^{(A_{\blambda})}(r,\lambda_0))^2z^2\,\mu(\dd z)\dd r
   +\int_{0}^{t} \int_{(\cU_1\setminus A_{\blambda})\cap (1,\infty)} 2 v^{(A_{\blambda})}(r,\lambda_0)z\,\mu(\dd z)\dd r\\
 &\leq \int_{0}^{t} (v^{(A_{\blambda})}(r,\lambda_0))^2\,\dd r
     \int_{ (0,1]} z^2\,\mu(\dd z)
     + 2 \int_{0}^{t} v^{(A_{\blambda})}(r,\lambda_0)\,\dd r 
         \int_{(1,\infty)}  z\,\mu(\dd z)<\infty,
         \qquad t\in\RR_+.
\end{split}
 \end{equation*}
Therefore, one can apply the dominated convergence theorem and \eqref{help_def_v} to obtain \eqref{eq:help_S4_1}.
Hence, by \eqref{help_AppB_1}--\eqref{eq:help_S4_1}, we get that 
 \begin{align*}
  &\lim_{s\uparrow \infty} \int_0^t \varphi_1^{(k)}\big( v^{(k)}(r,\eta_s,s\lambda_1,\ldots,s\lambda_k),s\lambda_1,\ldots,s\lambda_k\big)\,\dd r\\
  &\qquad = -\mu(A_{\blambda})t - \int_{0}^{t} B \tv^{(A_{\blambda})}(r,\lambda_0) \,\dd r
       + \int_{0}^{t} c (\tv^{(A_{\blambda})}(r,\lambda_0))^2 \,\dd r\\
  &\phantom{\qquad=\;} + \int_{0}^{t} \int_{A_{\blambda}} 
            \tv^{(A_{\blambda})}(r,\lambda_0) (1 \land z)\mu(\dd z)\dd r\\
  &\phantom{\qquad=\;} + \int_{0}^{t} \int_{\cU_1\setminus A_{\blambda}} 
           \Big( {\ee}^{-\tv^{(A_{\blambda})}(r,\lambda_0)z } - 1 + \tv^{(A_{\blambda})}(r,\lambda_0)(1 \land z)\Big)\mu(\dd z)\dd r\\
  &\qquad = -\mu(A_{\blambda})t  + \int_{0}^{t} \varphi^{(A_{\blambda})}(\tv^{(A_{\blambda})}(r,\lambda_0))\,\dd r,
  \qquad t\in\RR_+,
 \end{align*}
 where the second equality follows by the definition of $\varphi^{(A_{\blambda})}$ 
 (see \eqref{eq:psiA} with $A=A_\blambda$).
This implies \eqref{eq:help210_limit_lambda_0}, as desired.
In particular, for all $\lambda_0\in\RR_+$, $\eta_s\uparrow \lambda_0$ as $s\uparrow \infty$, and $\blambda\in\RR_+^k$,
 we also have that $\RR_+\ni t\mapsto \tv^{(A_{\blambda})}(t,\lambda_0)$ is continuously differentiable, and hence it is locally bounded.

All in all, for all $\lambda_0\in\RR_+$, $\eta_s\uparrow \lambda_0$ as $s\uparrow \infty$, and $\blambda\in\RR_+^k$,
 the function $\RR_+\ni t\mapsto \tv^{(A_{\blambda})}(t,\lambda_0)$ is a locally bounded non-negative
 solution to the differential equation \eqref{eq_vAlambda}.
 
Since the differential equation \eqref{eq_vAlambda} has a unique locally bounded non-negative solution (see Step 1),
 we have that, for all $\lambda_0\in\RR_+$, $\eta_s\uparrow \lambda_0$ as $s\uparrow \infty$, and $\blambda\in\RR_+^k$,
 the function $\RR_+\ni t\mapsto \tv^{(A_{\blambda})}(t,\lambda_0)\in\RR_+$ (introduced in \eqref{help_def_v}) is 
 the unique locally bounded non-negative solution to the differential equation \eqref{eq_vAlambda}.

Furthermore, since the differential equation \eqref{eq_vAlambda} depends on $\lambda_0\in\RR_+$
 and $\blambda\in\RR_+^k$, but not on the function $\RR_+\ni s\mapsto \eta_s$,
 we obtain that $\tv^{(A_{\blambda})}(t,\lambda_0)$ depends on $\lambda_0$ and $\blambda$, but not on the choice
 of the function $\RR_+\ni s\mapsto \eta_s$ (only on its limit $\lambda_0$ as $s\uparrow \infty$)
 justifying our notation in \eqref{help_def_v}.

\section*{Acknowledgements}
We would like to thank Professor Zenghu Li for his constant help in discussing our questions.

%\addcontentsline{toc}{section}{References}


\begin{thebibliography}{99}

\bibitem{App}
\textsc{Applebaum, D.} (2009).
\textit{L\'evy Processes and Stochastic Calculus, 2nd ed.}
Cambridge University Press, Cambridge.

 
\bibitem{BarBenKebPap}
\textsc{Barczy, M., Ben Alaya, M., Kebaier, A.} and \textsc{Pap, G.} (2019).
Asymptotic properties of maximum likelihood estimator for the growth rate of a stable CIR process based on continuous time observations.
\textit{Statistics}
\textbf{53(3)} 533--568. 

\bibitem{BarLiPap2}
\textsc{Barczy, M., Li, Z.} and \textsc{Pap, G.} (2015).
Stochastic differential equation with jumps for multi-type continuous state
 and continuous time branching processes with immigration.
\textit{ALEA. Latin American Journal of Probability and Mathematical Statistics}
\textbf{12(1)} 129--169.

\bibitem{BarPal}
\textsc{Barczy, M.} and \textsc{Palau, S.} (2024).
Distributional properties of jumps of multi-type CBI processes.
\textit{Electronic Journal of Probability}
\textbf{29} 1--39.

\bibitem{DufFilSch}
\textsc{Duffie, D., Filipovi\'{c}, D.} and \textsc{Schachermayer, W.} (2003).
Affine processes and applications in finance.
\textit{Annals of Applied Probability}
\textbf{13} 984--1053.

\bibitem{Fel}
\textsc{Feller, W.} (1950).
\textit{Diffusion processes in genetics}.
In: Proceedings of the Second Berkeley Symposium on Mathematical Statistics and Probability,
University of California Press, Berkeley and Los Angeles, 227--246.

\bibitem{FilMaySch}
\textsc{Filipovi\'{c}, D., Mayerhofer, E.} and \textsc{Schneider, P.} (2013).
Density approximations for multivariate affine jump-diffusion processes.
\textit{Journal of Econometrics}
\textbf{176(2)} 93--11.
 
\bibitem{FuLi}
\textsc{Fu, Z.} and \textsc{Li, Z.} (2010).
Stochastic equations of non-negative processes with jumps.
\textit{Stochastic Processes and their Applications}
\textbf{120} 306--330.

\bibitem{HeLi}
\textsc{He, X.} and \textsc{Li, Z.} (2016).
Distributions of jumps in a continuous-state branching process with immigration.
\textit{Journal of Applied Probability}
\textbf{53} 1166--1177.

\bibitem{JiaoMaScotti}
\textsc{Jiao, Y., Ma, C.} and \textsc{Scotti, S.} (2017).
Alpha-CIR model with branching processes in sovereign interest rate modeling.
\textit{Finance and Stochastics}
\textbf{21(3)} 789--813.

\bibitem{K2}
\textsc{Kallenberg, O.} (2017).
\textit{Random Measures, Theory and Applications}.
Springer, Cham.

\bibitem{KelSchTei}
\textsc{Keller-Ressel, M.}, \textsc{Schachermayer, W.} and \textsc{Teichmann, J.} (2011).
Affine processes are regular.
\textit{Probabaility Theory and Related Fields}
\textbf{151} 591--611.

\bibitem{Kin}
\textsc{Kingman, J. F. C.} (1993).
\textit{Poisson Processes}.
Oxford University Press, New York.

\bibitem{Li}
\textsc{Li, Z.} (2022).
\textit{Measure-Valued Branching Markov Processes, 2nd ed. }
  Springer-Verlag GmbH, Berlin.

\bibitem{Li3}
\textsc{Li, Z.} (2020).
Continuous-state branching processes with immigration.
 From probability to finance--lecture notes of BICMR Summer School on Financial Mathematics.
\textit{Mathematical Lectures from Peking University,}
\textbf{1--69.}
Springer, Singapore.


\bibitem{Vol}
\textsc{Volkmann, P.} (1972).
Gew\"ohnliche Differentialungleichungen mit quasimonoton wachsenden
 Funktionen in topologischen Vektorr\"aumen.
\textit{Mathematische Zeitschrift}
\textbf{127} 157--164.

\end{thebibliography}
\end{document}